\documentclass[preprint,review,10pt]{elsarticle}

\usepackage{lineno,hyperref}
\usepackage{graphicx}
\usepackage{subfigure}
\usepackage{amsmath}
\usepackage{geometry}
\usepackage{booktabs}
\usepackage{threeparttable}
\usepackage{textcomp}
\usepackage{multirow}
\usepackage{color,xcolor}
\usepackage{amssymb}
\usepackage{changepage}
\usepackage[latin1]{inputenc}

\modulolinenumbers[50]
\journal{The Report}

\begin{document}

\begin{frontmatter}

\title{Nonlinear Analysis in $p$-Vector Spaces for Singe-Valued 1-Set Contractive Mappings}


\author{George Xianzhi YUAN}
\address{Chengdu University, Chengdu 610601, China}
\address{Mathematics, Sichuan University, Chengdu 610065, China}
\address{Sun Yat-Sen University, Guangzhou 510275, China and}
\address{East China University of Science and Technology, Shanghai 200237, China}
\address{\bf george\_yuan99@yahoo.com}

\begin{abstract}
The goal of this paper is to develop some fundamental and important nonlinear analysis for single-valued mappings under the framework of $p$-vector spaces, in particular, for locally $p$-convex spaces for $0 < p \leq 1$. More precisely, based on the fixed point theorem of single-valued continuous condensing mapping in locally $p$-convex spaces as the starting point, we first establish best approximation results for (single-valued) continuous condensing mappings which are then used to develop new results for three classes of nonlinear mappings consisting of 1) condensing; 2) 1-set contractive; and 3) semiclosed 1-set contractive mappings in locally $p$-convex spaces. Next they are used to establish general principle for nonlinear alternative, Leray - Schauder alternative, fixed points for non-self mappings with different boundary conditions for nonlinear mappings from locally $p$-convex spaces, to nonexpansive mappings in  uniformly convex Banach spaces, or locally convex spaces with Opial condition. The results given by this paper not only include the corresponding ones in the existing literature as special cases, but also expected to be useful tools for the development of new theory in nonlinear functional analysis and applications to the study of related nonlinear problems arising from practice under the general framework of $p$-vector spaces for $0< p \leq 1$.

Finally, the  work presented by this paper focuses on the development of nonlinear analysis for single-valued (instead of set-valued) mappings for locally $p$-convex spaces, essentially, is indeed the continuation of the associated work given recently by Yuan \cite{yuan2022}  therein, the attention is given to the study of nonlinear analysis for set-valued mappings in locally $p$-convex spaces for  $0 < p  \leq 1$.

\end{abstract}

\vskip.1in
\begin{keyword}
Nonlinear analysis, $p$-convex, Fixed points, Measure of noncompactness, Condensing mapping, 1-set contractive mapping, Semiclosed mapping,
Nonexpansive mapping, Best approximation, Nonlinear alternative, Leray - Schauder alternative, Demiclosed principle, Opial condition,  $p$-inward and $p$-outward set, $p$-vector space, locally $p$-convex space, Uniform convex space.

\vskip.1in
\noindent
AMS Classification: 47H04, 47H10, 46A16, 46A55, 49J27, 49J35,  52A07,  54C60, 54H25, 55M20 

\end{keyword}

\end{frontmatter}

\newpage
\linenumbers
\section{Introduction}
\label{intro}

\vskip.1in
It is known that the class of $p$-semi-norm spaces $(0 < p \leq 1)$ is an important generalization of usual normed spaces with  rich topological and geometrical structures, and related study has received a lot of attention, e.g., see work by Alghamdi et al.\cite{alghamdiporegan}, Balachandran \cite{balachandra}, Bayoumi \cite{bayoumi}, Bayoumi et al.\cite{bayoumi2015}, Bernu\'{e}es and Pena \cite{bernues1997}, Ding \cite{ding}, Ennassik and Taoudi \cite{ennassik2021}, Ennassik et al.\cite{ennassik},
Gal and Goldstein \cite{gal}, Gholizadeh et al.\cite{gholizadeh}, Jarchow \cite{jarchow}, Kalton \cite{kalton1977}-\cite{kalton1977b}, Kalton et al.\cite{kalton1984}, Machrafi and Oubbi \cite{machrafioubbi}, Park \cite{park2016}, Qiu and Rolewicz \cite{qiu}, Rolewicz \cite{rolewicz}, Silva et al.\cite{silva}, Simons \cite{simons}, Tabor et al.\cite{tabor}, Tan \cite{tan}, Wang \cite{wang}, Xiao and Lu \cite{xiaolu}, Xiao and Zhu \cite{xiaozhu2011}-\cite{xiaozhu2018}, Yuan \cite{yuan2022}, and many others. However, to the best of our knowledge, the corresponding basic tools and associated results in the category of nonlinear functional analysis for $p$-vector spaces have not been well developed, in particular for the three classes of (single-valued) continuous nonlinear mappings which are: 1) condensing; 2) 1-set contractive; and 3) semiclosed 1-set contractive operators under locally $p$-convex spaces. Our goal in this paper is to develop some fundamental and important nonlinear analysis for single-valued mappings under the framework of $p$-vector spaces, in particular, for locally $p$-convex spaces for $0 < p \leq 1$. More precisely, based on the fixed point theorem of single-valued continuous condensing mapping in locally $p$-convex spaces as the starting point, we first establish best approximation results for (single-valued) continuous condensing mappings which are then used to develop new results for three classes of nonlinear mappings, which are 1): condensing; 2): 1-set contractive; and 3): semiclosed 1-set contractive in locally $p$-convex spaces. Then these new results are used to establish general principle for nonlinear alternative, Leray - Schauder alternative, fixed points for non-self mappings with different boundary conditions for nonlinear mappings from locally $p$-convex spaces, to  nonexpansive mappings in  uniformly convex Banach spaces, or locally convex spaces with Opial condition. The results given by this paper not only include the corresponding results in the existing literature as special cases, but also expected to be useful tools for the development of new theory in nonlinear functional analysis and applications to the study of related nonlinear problems arising from practice under the general framework of $p$-vector spaces for $0< p \leq 1$.

In addition, we like to point out that the work presented by this paper focuses on the development of nonlinear analysis for single-valued (instead of set-valued) mappings for locally $p$-convex spaces, essentially, is very important, and also the continuation of the work given recently by Yuan \cite{yuan2022} therein, the attention was given to establish new results on fixed points, principle of nonlinear alternative for nonlinear mappings mainly on set-valued (instead of single-valued) mappings developed in locally $p$-convex spaces for $0 < p  \leq 1$.
Though some new results for set-valued mappings in locally $p$-convex spaces have been developed (see Gholizadeh et al.\cite{gholizadeh}, Park \cite{park2016}, Qiu and Rolewicz \cite{qiu}, Xiao and Zhu \cite{xiaozhu2011}-\cite{xiaozhu2018}, Yuan \cite{yuan2022} and others), we still like to emphasize that results obtained for set-valued mappings  for $p$-vector spaces may face some challenging in dealing with true nonlinear problems.  One example is that the assumption used for $``$set-valued mappings with closed $p$-convex values" seems too strong as it always means that the  zero element is a trivial fixed point of the set-valued mappings, and this was also discussed in P.40-41 by Yuan \cite{yuan2022} for $0 < p \leq 1$.

For the development since 1920s on the development, and in particular, how the fixed points for non-self mappings, best approximation method and related to the study on some key aspects of nonlinear analysis related to Birkhoff-Kellogg problems, nonlinear alternative, Leray - Schauder alternative, KKM principle, best approximation, and related topics, readers can find some most important contributions by Birkhoff and Kellogg \cite{bk} in 1920's, Leray and Schauder \cite{lerayschauder} in 1934's, Fan \cite{fan1969} in 1969; plus the related comprehensive references given by Agarwal et al.\cite{Agarwal}, Bernstein \cite{berbstein}, Chang et al.\cite{chang1993}, Granas and Dugundji \cite{granas}, Isac \cite{isac},  Park \cite{park1997}, Singh et al.\cite{singh1997}, Zeidler \cite{Zeidler}; and also see work contributed by
Agarwal and O'Regan \cite{agarwalorgan2003}-\cite{agarwalorgan2004}, Furi and Pera \cite{furipera}, Park \cite{park1997}, O'Regan \cite{oregan2019}, O'Regan and Precup \cite{OP}), Poincare \cite{poincare2},  
Rothe \cite{rothe1986}-\cite{rothe1981}, Yuan \cite{yuan1998}-\cite{yuan2022}, Zeidler \cite{Zeidler}.

It is well-known that the best approximation is one of very important aspects for the study of nonlinear problems related to the problems on their solvability for partial differential equations, dynamic systems, optimization, mathematical program, operation research; and in particularly, the one approach well accepted for study of nonlinear problems in optimization, complementarity problems and of variational inequalities problems and so on, strongly based on today called Fan's best approximation theorem given by Fan \cite{fan1952}-\cite{fan1972} in 1969 which acts as a very powerful tool in nonlinear analysis, and see the book of Singh et al.\cite{singh1997} for the related discussion and study on the fixed point theory and best approximation with the KKM-map principle), among them, the related tools are Rothe type and principle of Leray-Schauder alterative in topological vector spaces (TVS), and local topological vector spaces (LCS) which are comprehensively studied by Chang et al.\cite{chang1993},
Chang et al.\cite{chang1996}-\cite{chang2001}, 
Carbone and Conti \cite{carbone1991}, 
Ennassik and Taoudi \cite{ennassik2021}, Ennassik et al.\cite{ennassik},
Isac \cite{isac}, Granas and Dugundji \cite{granas}, Kirk and Shahzad \cite{kirk2014}, Liu \cite{liu2001},
Park \cite{park100}, Rothe \cite{rothe1986}-\cite{rothe1981}, Shahzad \cite{shahzad2006}-\cite{shahzad2004}, Xu \cite{xu1991}, Yuan \cite{yuan1998}-\cite{yuan2022}, Zeidler \cite{Zeidler},  
and references therein. 

On the other hand, since the celebrated so-called KKM principle established in 1929 in \cite{kkm}, was based on the celebrated Sperner combinatorial lemma and first applied to a simple proof of the Brouwer fixed point theorem. Later it became
clear that these three theorems are mutually equivalent and they were regarded as a sort of mathematical
trinity (Park \cite{park100}). Since Fan extended the classical KKM theorem to infinite-dimensional spaces in 1961
by Fan \cite{fan1960}-\cite{fan1972}, there have been a number of generalizations and applications in numerous
areas of nonlinear analysis, and fixed points in TVS and LCS as developed by Browder \cite{bro1965}-\cite{bro1983} and related references therein. Among them, Schauder's fixed point theorem \cite{schauder} in normed spaces is one of the powerful
tools in dealing with nonlinear problems in analysis. Most notably, it has played a major role in the
development of fixed point theory and related nonlinear analysis and mathematical theory of partial and
differential equations and others. A generalization of Schauder's theorem from normed space to general
topological vector spaces is an old conjecture in fixed point theory which is explained by the Problem 54 of the book $``$The Scottish Book" by Mauldin \cite{Mauldin} as stated as Schauder's conjecture: $``$Every nonempty compact convex set in a topological vector space has the fixed point property, or in its analytic statement, does a continuous function defined on a compact convex subset of a topological vector space to itself have a fixed point?" Recently, this question has been recently answered by the work of Ennassik and Taoudi \cite{ennassik2021}
by using $p$-seminorm methods under locally $p$-convex spaces! See also related work in this direction given by  Askoura and Godet-Thobie \cite{Askoura}, Cauty \cite{cauty2005}-\cite{cauty2007}, Chang \cite{chang1997}, Chang et al.\cite{chang1993}, Chen \cite{chen},
Dobrowolski \cite{dobrowolski2003}, Gholizadeh et al.\cite{gholizadeh}, G\'{o}rniewicz \cite{gorniewicz},  G\'{o}rniewicz et al.\cite{gorniewiczetal}, Isac \cite{isac}, Li \cite{li}, Li et al.\cite{lixuduan2006}, Liu \cite{liu2001},  Nhu \cite{nhu1996}, Okon \cite{okon}, Park \cite{park2016}-\cite{park2022}, Reich \cite{reich}, Smart \cite{smart}, Weber \cite{weber1992}-\cite{weber1991}, Xiao and Lu \cite{xiaolu}, Xiao and Zhu \cite{xiaozhu2011}-\cite{xiaozhu2018}, Xu \cite{xu2007}, Xu et al.\cite{xujiali2006}, Yuan \cite{yuan1998}-\cite{yuan2022}, and related references therein under the general framework of $p$-vector spaces, in particular, locally $p$-convex spaces for non-self mappings with various boundary conditions  for $0 < p \leq 1$.

The goal of this paper is to establish the general new tools of nonlinear analysis  under the framework of general locally $p$-convex space ($p$-seminorm spaces) for general condensing mappings,  1-set contractive mappings, and semiclosed mappings  (here $0 < p \leq 1$), and we do wish these new results such as best approximation, theorems of Birkhoff-Kellogg type, nonlinear alternative, fixed point theorems for non-self (singl-valued) continuous operators with various boundary conditions, Rothe, Petryshyn type, Altman type, Leray-Schedule types, related others nonlinear problems would play important roles for the nonlinear analysis of $p$-seminorm spaces for $0 < p \leq 1$. In addition, our results also show that fixed point theorem for condensing continuous mappings for closed $p$-convex subsets provide solutions for Schauder's conjecture since 1930'a in the affirmative way under the general setting of $p$-vector spaces (which may not locally convex, see related study given by
Ennassik and Taoudi \cite{ennassik2021}, Kalton \cite{kalton1977}-\cite{kalton1977b}, Kalton et al.\cite{kalton1984}, Jarchow \cite{jarchow}, Roloewicz \cite{rolewicz} on this direction).

The paper has seven sections. Section 1 is the introduction. Section 2 describes general concepts for the $p$-convex subsets of  topological vector spaces ($0 < p \leq 1$). In Section 3, then some basic results of KKM principle related to abstract convex spaces are given. In Section 4, as the application of the KKM principle in abstract convex spaces which including $p$-convex vector spaces as a special class ($0< p \leq 1)$ by combining the embedding lemma for compact $p$-convex subsets from topological vector spaces into locally $p$-convex spaces, we provide general fixed point theorems for condensing continuous mappings for both single-valued version in topological vector spaces; and upper
semi-continuous set-valued version in locally convex spaces defined on closed $p$-convex subsets for $0 < p \leq 1$. The Sections 5, 6 and 7 mainly focus on the study of nonlinear analysis for 1-set contractive (single-valued ) continuous mappings in locally $p$-convex vector spaces  to establish the general existence theorems for solutions of Birkhoff-Kellogg (problem) alternative, general principle of nonlinear alterative, and including Leray-Schauder alternative, Rothe type, Altman type associated with different boundary conditions. The Sections 8, 9 and 10 mainly focus on the study of new results based on semiclosed
1-set contractive (single-valued) continuous mappings related to  nonlinear alternative principles, Birkhoff-Kellogg theorems, Leray-Schauder alternative and non-self operations from general locally $p$-convex spaces to uniformly convex Banach spaces for nonexpansive mappings, or locally convex topological spaces with Opial condition.

For the convenience of our discussion, throughout this paper, we always assume that all
$p$-vector spaces are Hausdorff for $0 < p \leq 1 $ unless specified; and we also denote by $\mathbb{N}$ the set of all positive integers, i.e., $\mathbb{N}:=\{1, 2, \cdots, \}$.

\vskip.1in
\section{Some Basic Results for $p$-Vector Spaces}
\label{sec:1}

For the convenience of self-containing reading dor readers, we recall some notion and definitions for $p$-convex vector spaces below as summarized by Yuan \cite{yuan2022}, see also Balachandran \cite{balachandra}, Bayoumi \cite{bayoumi}, Jarchow \cite{jarchow}, Kalton \cite{kalton1977}, Rolewicz \cite{rolewicz}, Gholizadeh et al.\cite{gholizadeh}, Ennassik and Taoudi \cite{ennassik2021}, Ennassik et al.\cite{ennassik}, Xiao and Lu \cite{xiaolu}, Xiao and Zhu \cite{xiaozhu2011} and references therein for more in details.

\noindent
{\bf Definition 2.1.}
A set $A$ in a vector space $X$ is said to be
$p$-convex for $0 < p \leq 1$ if, for any $x, y\in A$, $0 \leq s, t \leq 1$ with
$s^p + t^p=1$, we have
$s^{1/p}x+t^{1/p}y\in A$; and if $A$ is $1$-convex, it is simply called convex (for $p = 1$) in general vector spaces;
the set $A$ is said to be absolutely $p$-convex if $s^{1/p}x+t^{1/p}y\in A$ for $0 \leq |s|, |t| \leq 1$ with $|s|^p + |t|^p \leq 1$.

\noindent
{\bf Definition 2.2.} If $A$  is a subset of a topological vector space $X$, the closure of $A$ is denoted by $\overline{A}$,
then the $p$-convex hull of $A$ and its closed $p$-convex hull denoted by $C_{p}(A)$, and $\overline{C}_{p}(A)$, respectively, which is the smallest $p$-convex set containing $A$, and the smallest closed $p$-convex set containing $A$, respectively.

\noindent
{\bf Definition 2.3.} Let $A$ be $p$-convex and $x_{1}, \cdots, x_{n}\in A$, and $t_{i}\geq 0, \displaystyle \sum_{1}^{n}t_{i}^{\mathrm{p}}=1$. Then $\displaystyle \sum_{1}^{n}t_{i}x_{i}$ is called a $p$-convex combination of $\{x_i\}$ for $i=1, 2, \cdots, n$.  If $\displaystyle \sum_{1}^{n}|t_{i}|^{\mathrm{p}}\leq 1$,
then $\displaystyle \sum_{1}^{n}t_{i}x_{i}$ is called an absolutely $p$-convex combination.
It is easy to see that $\displaystyle \sum_{1}^{n}t_{i}x_{i}\in A$ for a $p$-convex set $A$.

\noindent
{\bf Definition 2.4.} A subset $A$  of a vector space $X$ is called circled (or balanced) if
$\lambda A \subset A$ holds for all scalars $\lambda$ satisfying $|\lambda| \leq 1$. We say that $A$ is absorbing
if for each $x \in X$, there is a real number $\rho_x >0$ such that $\lambda x \in A$ for all $\lambda > 0$ with $|\lambda |\leq \rho_x$.

By the definition 2.4, it is easy to see that the system of all circled subsets of $X$ is easily seen to be closed under
the formation of linear combinations, arbitrary unions, and arbitrary intersections.
In particular, every set $A \subset X$ determines a smallest circled subset $\hat A$ of $X$ in which it is contained:
$\hat A$ is called the circled hull of $A$. It is clear that $\hat A =\cup_{|\lambda |\leq 1} \lambda A$ holds, so that $A$ is circled if and only if (in short, iff) $\hat A =A$. We use $\overline{\hat A}$ to denote for the closed circled hull of $A\subset X$.

In addition, if $X$ is a topological vector space, we  use the $int(A)$ to denote the interior of set $A \subset X$ and if $0 \in int(A)$, then $int(A)$ is also circled, and using $\partial A$ to denote the boundary of $A$ in $X$ unless specified.

\noindent
{\bf Definition 2.5.} A topological vector space is said to be locally
$p$-convex if the origin has a fundamental set of absolutely $p$-convex $0$-neighborhoods.
This topology can be determined by $p$-seminorms which are defined in the obvious
way (see P.52 of Bayoumi \cite{bayoumi}, Jarchow \cite{jarchow} or Rolewicz \cite{rolewicz}).

\noindent
{\bf Definition 2.6.} Let $X$ is a vector space and $\mathbb{R}^+$ is a non-negative part of a real line $\mathbb{R}$. Then a mapping $P: X\longrightarrow \mathbb{R}^+$ is said to be a $p$-seminorm if it satisfies the requirements for $(0 < p \leq 1)$

(i) $P(x) \geq 0$ for all $x \in X$;

(ii) $P(\lambda x) = |\lambda|^p P(x)$  for all $x\in X$ and $\lambda \in R$;

(iii) $P(x + y) \leq P(x) +  P(y)$  for all $x, y \in X$.

An $p$-seminorm $P$ is called a $p$-norm if $x=0$  whenever $P(x)=0$, so a vector space with a specific $p$-norm is called an $p$-normed space, and of course if $p=1$, $X$ is a normed space as discussed beofe (e.g., see Jarchow \cite{jarchow}).

By Lemma 3.2.5 of Balachandran \cite{balachandra}, the following proposition gives a necessary and sufficient condition for
an $p$-seminorm to be continuous.

\vskip.1in
\noindent
{\bf Proposition 2.1.} Let $X$  be a topological vector space, $P$ is a $p$-seminorm on $X$ and
 $V: =\{ x\in X: P(x) < 1\}$. Then $P$ is continuous if and only if $0 \in int(V)$, where $int(V)$ is the interior of $V$.

Now given an $p$-seminorm $P$, the $p$-seminorm topology determined by $P$ (in short, the $p$-topology) is the class of unions of open balls $B(x, \epsilon): = \{ y \in X: P(y-x) < \epsilon\}$ for $x \in X$ and $\epsilon> 0$.

\noindent
{\bf Definition 2.7.} A topological vector space $X$ is said to be locally $p$-convex if it has a $0$-basis consisting of $p$-convex neighborhoods for $(0 < p \leq 1)$. If $p=1$, $X$ a usual locally convex space.

We also need the following notion for the so-called $p$-gauge (see Balachandran \cite{balachandra}).

\noindent
{\bf Definition 2.8.}  Let $A$ be an absorbing subset of a vector space $X$. For $ x \in X$  and $0 < p \leq 1$, set
$P_A=\inf\{\alpha >0:  x \in \alpha^{\frac{1}{p}}A\}$, then the non-negative real-valued function $P_A$ is called the $p$-gauge (gauge if $p=1$). The $p$-gauge of $A$ is also known as the Minkowski $p$-functional.

By Proposition 4.1.10 of Balachandran \cite{balachandra}, we have the following proposition.

\vskip.1in
\noindent
{\bf Proposition 2.2.} Let $A$ be an absorbing subset of $X$. Then $p$-gauge $P_A$ has  the following properties:

(i) $P_A(0)=0$;

(ii) $P_A(\lambda x) = |\lambda|^p P_A(x)$ if $\lambda \geq 0$;

(iii) $P_A(\lambda x) = |\lambda |^p P_A(x)$ for all $\lambda \in R$ provided $A$ is circled;

(iv) $P_A(x + y) \leq P_A(x) + P_A(y)$  for all $x, y \in A$ provided $A$ is $p$-convex.

In particular, $P_A$ is a $p$-seminorm if $A$  is absolutely $p$-convex (and also absorbing).

As mentioned above, a given $p$-seminorm is said to be an $p$-norm if $x = 0$  whenever $P(x) = 0$. A vector space with
a specific $p$-norm is called a $p$-normed space. The $p$-norm of an element $x \in E$  will
usually be denoted by $\|x\|_p$. If $p = 1$, $X$ is a usual normed space. If $X$ is an $p$-normed
space, then $(X, d_p)$ is a metric linear space with a translation invariant metric $d_p$ such
that $d_p=d_p(x, y)=\|x -y\|_p$ for $x, y \in X$. We point out that $p$-normed spaces are
very important in the theory of topological vector spaces. Specifically, a Hausdorff
topological vector space is locally bounded if and only if it is a $p$-normed space for
some $p$-norm $\| \cdot \|_p$, where $0 < p \leq 1$ (see p.114 of Jarchow \cite{jarchow}).
We also note that examples of $p$-normed spaces include such as $L^p(\mu)$ - spaces and Hardy spaces $H_p$, $0 < p < 1$, endowed with their usual $p$-norms.

\vskip.1in
\noindent
{\bf Remark 2.1.} We like to make the following important two points as follows:

(1) First, by the fact that (e.g., see Kalton et al.\cite{kalton1984}, or Ding \cite{ding}), there is no open convex non-void subset in $L^p[0, 1]$ (for $0< p < 1$) except $L^p[0,1]$ itself, this means that $p$-normed paces with $0< p <1$ are not necessarily locally convex. Moreover, we know that every $p$-normed space is locally $p$-convex; and incorporating Lemma 2.3 below, it seems that $p$-vector spaces (for $0 < p \leq 1$ ) is a nicer space as we can use $p$-vector space to approximate (Hausdorff) topological vector spaces (TVS) in terms of Lemma 2.1 (ii) below for the convex subsets in TVS by using a bigger $p$-convex subsets in $p$-vector spaces for $p\in (0,1)$ by also considering Lemma 2.3 below, in this way, it seems $P$-vector spaces seems having better properties in terms of $p$-convexity than the usually ($1-$) convex subsets used in TVS with $p=1$.

(2) Second, it is worthwhile noting that a $0$-neighborhood in a topological
vector space is always absorbing by Lemma 2.1.16 of Balachandran \cite{balachandra}, or Proposition 2.2.3 of Jarchow \cite{jarchow}.

\vskip.1in
Now by Proposition 4.1.12 of Balachandran \cite{balachandra}, we also have the following Proposition 2.3 and Remark 2.2 (which is the Remark 2.3 of Ennassik and Taoudi \cite{ennassik2021}).

\noindent
{\bf Proposition 2.3.} Let $A$ be a subset of a vector space $X$, which is absolutely $p$-convex $(0 < p \leq 1)$ and absorbing. Then, we have that

(i) The $p$-gauge $P_A$ is a $p$-seminorm such that if
$B_1: =\{x \in X: P_A(x) < 1\}$, and  $\overline{B_1}=\{ x \in X: P_A(x) \leq 1\}$. then
$B_1\subset A \subset \overline{B_1}$; in particular, $ker P_A \subset A$, where $ker P_A: =\{ x \in X: P_A(x) =  0 \}$.

(ii) $ A = B_1$ or $\overline{B_1}$ according as $A$ is open or closed in the $P_A$-topology.

\noindent
{\bf Remark 2.2.} Let $X$ be a topological vector space and let $U$ be an open
absolutely $p$-convex neighborhood of the origin, and let $\epsilon$ be given.
If $y \in \epsilon^{\frac{1}{p}} U$, then $y=\epsilon^{\frac{1}{p}} u$ for some $u \in U$ and $P_U(y)= P_U(\epsilon^{\frac{1}{p}} u)= \epsilon P_U(u) \leq \epsilon$ (as $u \in U$ implies that $P_U(u) \leq 1$). Thus, $P_U$ is continuous at
$zero$, and therefore, $P_U$ is continuous everywhere. Moreover, we have $U=\{ x \in X: P_U(x) < 1\}.$

Indeed, since $U$ is open and the scalar multiplication is continuous, we
have that for any $x \in U$, there exists $0 < t < 1$ such that $x \in t^{\frac{1}{p}} U$ and so $P_U(x) \leq t < 1$.This shows that $U \subset \{ x\in X: P_U(x) < 1\}$. The conclusion follows by Proposition 2.3 above.


The following result is a very important and useful result which allows use to make the approximation for convex subsets in topological vector spaces by $p$-convex subsets in $p$-convex vector spaces. For the reader's self-contained in reading, we provide a sketch of proof below (see also Lemma 2.1 of Ennassik and Taoudi \cite{ennassik}, Remark 2.1 of Qiu and Rolewicz \cite{qiu}).

\noindent
{\bf Lemma 2.1.} Let $A$ be a subset of a vector space $X$, then we have

(i) If $A$ is $p$-convex, with $0 < p < 1$, then $\alpha x \in A$  for any $x \in A$ and any $0 < \alpha \leq 1$.

(ii) If $A$ is convex and $0 \in A$, then $A$ is $p$-convex for any $p \in (0, 1]$.

(iii) If $A$ is $p$-convex for some $p \in (0, 1)$, then $A$ is $s$-convex for any $s \in (0, p]$.

\noindent
{\bf Proof.}  (i) As $r \le 1$, by the fact that
$``$for all $x \in A$ and all $\alpha  \in [2^{(n+1)(1-\frac{1}{p})}, 2^{n(1-\frac{1}{p})}]$, we have $\alpha x \in A$" is true for all integer $n\geq 0$. Taking into account that the fact that $(0, 1]=\cup_{n\geq 0}  [2^{(n+1)(1-\frac{1}{p})}, 2^{n(1-\frac{1}{p})}]$,  thus the result is obtained.

(ii)  Assume that $A$  is a convex subset of $X$ with $ 0 \in A$  and take a real number $s \in (0, 1]$. we show that  $A$ is $s$-convex. Indeed, let $x, y \in A$ and $\alpha, \beta >0$ with $\alpha^p + \beta^p = 1$. Since $A$ is convex, then
$\frac{\alpha}{\alpha + \beta} x + \frac{\beta}{\alpha + \beta}y \in A$. Keeping in mind that $0 < \alpha + \beta < \alpha^p + \beta^p=1$, it follows that $\alpha x + \beta y=(\alpha + \beta )(\frac{\alpha}{\alpha + \beta}x + \frac{\beta}{\alpha +\beta}y ) + (1-\alpha -\beta) 0 \in A$.

(iii) Now, assume that $A$ is $r$-convex for some $p \in (0,1)$ and pick up any real $s \in (0, p]$.
We show that $A$ is $s$-convex. To see this, let $x, y \in A$  and $\alpha, \beta > 0$  such that $\alpha^s + \beta^s=1$. First notice that $ 0 < \alpha^{\frac{p - s}{p}} \leq 1$ and $0 < \beta^{\frac{p - s}{p}} \leq 1$, which imply that $\alpha^{\frac{p - s}{p}} x \in A$ and
$\beta^{\frac{p - s}{p}} y \in A$. By the $p$-convexity of $A$ and the equality $(\alpha^{\frac{s}{p}})^p + (\beta^{\frac{s}{p}})^p =1$, it follows that $\alpha x + \beta y = \alpha^{\frac{s}{p}}(\alpha^{\frac{p-s}{p}}x) + \beta^{\frac{s}{p}}(\beta^{\frac{p-s}{p}} y) \in A$.
This competes the sketch of the proof.  $\square$

\noindent
{\bf Remark 2.3.} We like to point out that the results (i) and (iii) of Lemma 2.1 do not hold for $p = 1$. Indeed, any singleton
$\{x\} \subset X$ is convex in topological vector spaces; but if $x \neq 0$, then it is not $p$-convex for any $p \in (0, 1)$.

We also need the following Proposition which is proposition 6.7.2 of Jarchow \cite{jarchow}.

\noindent
{\bf Proposition 2.4.} Let $K$ be compact in a topological vector $X$ and $(1< p \leq 1)$. Then the closure $\overline{C}_p(K)$ of the $p$-convex hull, and the closure $\overline{AC}_p(K)$ of absolutely $p$-convex hull of $K$ are compact if and only if $\overline{C}_p(K)$ and $\overline{AC}_p(K)$ are complete, respectively.

\vskip.1in
We also need following fact, which is a special case of Lemma 2.4 of Xiao and Zhu \cite{xiaozhu2011}.

\noindent
{\bf Lemma 2.2.} Let $C$ be a bounded closed $p$-convex subset of $p$-seminorm $X$ with $0 \in int C$, where $(0< p\leq 1)$. For every $x\in X$ define an operator by $r(x):=\frac{x}{\max\{1, (P_C(x))^{\frac{1}{p}}\}}$, where $P_C$ is the Minkowski $p$-functional of $C$. Then $C$ is a retract of $X$ and $r: X \rightarrow C$ is a continuous such that

(1) if $x \in C$, then $r(x)=x$;

(2) if $x \notin C$, then $r(x) \in \partial C$;

(3) if $x \notin C$, then the  Minkowski $p$-functional $P_C(x) >1 $.

\noindent
{\bf Proof.} Taking $s = p$ in Lemma 2.4 of Xiao and Zhu \cite{xiaozhu2011}, Proposition 2.3 and Remark 2.2, thus the proof is compete. $\square$

\vskip.1in
\noindent
{\bf Remark 2.4.} As discussed by Remark 2.2, Lemma 2.2 still holds if $``$the bounded closed $p$-convex subset $C$ of the $p$-normed space $(X, \|\cdot \|_p)$"  is replaced by $``$$X$ is a $p$-seminorm vector space and $C$ is a bounded closed absorbing $p$-convex subset with $0 \in int C$ of $X$".

\vskip.1in
Before we close this section, we like to point out that the structure of $p$-convexity when $p \in (0, 1)$ is really different from what we normally have for the concept of $``$convexity" used in topological vector spaces (TVS), in particular, maybe the following fact is one of reasons for us to use better ($p$-convex) structures in $p$-vector spaces to approximate the corresponding structure of the convexity used in TVS (i.e., the $p$-vector space when $p=1$). Based on the discussion in P.1740 of Xiao and Zhu \cite{xiaozhu2011}(see also Bernu\'{e}s and Pena \cite{bernues1997} and Sezer et al.\cite{zezer2021}), we have the following fact which indicates that each $p$-convex subset is $``$bigger" than the convex subset in topological vector spaces for $0 < p < 1$.

\vskip.1in
\noindent
{\bf Lemma 2.3.} Let $x$ be a point of $p$-vector space $E$, where assume $0 < p < 1$, then the $p$-convex hull and the closure of $\{x\}$ is given by
\begin{equation}\label{eq10}
C_p(\{x\})=\left\{
\begin{aligned}
&\{tx: t \in (0, 1]\}, &\mbox{ if }& x \neq 0,\\
&\{0\},                &\mbox{ if }& x = 0;  \\
\end{aligned}
\right.
\end{equation}
and
\begin{equation}\label{eq10}
\overline{C_p(\{x\})}=\left\{
\begin{aligned}
&\{tx: t \in [0, 1]\}, &\mbox{ if }& x \neq 0,\\
&\{0\},                &\mbox{ if }& x= 0.  \\
\end{aligned}
\right.
\end{equation}
But note that if $x$ is a given one point in $p$-vector space $E$, when $p=1$, we have that $\overline{C_1(\{x\})} =C_1(\{x\})=\{ x\}$, This shows significantly different for the structure of $p$-convexity between $p=1$ and $p\neq 1$!

\vskip.1in
As an application of Lemma 2.3, we have the following fact for (set-valued) mappings with non-empty closed $p$-convex values in $p$-vector spaces for $p \in (0, 1)$, which are truly different from  any (set-valued) mappings defined in topological vector spaces (i.e., for a $p$-vector space with $p =1)$.

\vskip.1in
\noindent
{\bf Lemma 2.4.} Let $U$ be a non-empty subset of a $p$-vector space $E$ (where $0 < p < 1)$,  with zero $0 \in U$, and assume a (set-valued) mapping $T: U \rightarrow 2^E$ is  with non-empty closed $p$-convex values. Then  $T$ has at least one fixed point in $U$, which is the element zero, i.e., $0 \in \cap_{x \in U} T(x) \ne \emptyset$.

\noindent
{\bf Proof.} For each $x \in U$, as $T(x)$ is non-empty closed $p$-convex, by Lemma 2.3, we have at leat $0 \in T(x)$. It implies that $0 \in \cap_{x \in U} T(x)$ and thus zero of $E$ is a fixed point of $T$. This completes the proof. $\square$

\vskip.1in
\noindent
{\bf Remark 2.5}. We like to point out that Lemma 2.4 shows that any set-valued mapping with closed $p$-convex values in $p$-spaces for $0< p < 1$ has the zero point as its trivial fixed point, thus it is very important to study fixed point and related principle of nonlinear analysis for single-valued (instead of set-valued) mappings for $p$-vector spaces (for $0 < p  < 1$), as pointed out and the discussion given in P.40-41 by Yuan \cite{yuan2022}, thus the most new results established in this paper is for singl-valued mappings for the three classes of (single-valued) continuous mappings which are: 1) condensing; 2) 1-set contractive; and 3) semiclosed 1-set contractive mappings. This is key different from those results obtained by Yuan \cite{yuan2022} recently for the study of set-valued mappings in $p$-vector spaces for $0 < p \leq 1$.

\vskip.1in
 By following Definitions 2.5 and 2.6, the discussion given by Proposition 2.3 and remarks thereafter, each given (open) $p$-convex subset $U$ in a $p$-vector space $E$ with the zero $0 \in int(U)$ always corresponds to a $p$-seminorm $P_U$, which is indeed the Minkowski $p$-functional of $U$ in $E$, and $P_U$ is continuous in $E$. In particular,  a topological vector space is said to be locally $p$-convex if the origin $0$ of $E$ has a fundamental set (denoted by) $\mathfrak{U}$, which is a family of absolutely $p$-convex $0$-neighborhoods (each denoted by $U$).
 This topology can be determined by $p$-seminorm $P_U$, which are indeed the family $\{P_U\}_{U \in \mathfrak{U}}$,
where $P_U$ is just the Minkowski $p$-functional for each $U \in \mathfrak{U}$ in $E$ (see also P.52 of Bayoumi [7], Jarchow [50] or Rolewicz [99]).

\vskip.1in
Throughout this paper, by following Remark 2.5, without loss of generality unless specified,
for a given $p$-vector space $E$, where $p \in (0, 1]$,
we always denote by $\mathfrak{U}$ the base of the $p$-vector space $E$'s topology structure, which is the family of its $0$-neighborhoods. For each $U \in \mathfrak{U}$, its corresponding $P$-seminorm $P_U$ is the Minkowski $p$-functional of $U$ in $E$. For a given point $x \in E$ and a subset $C \subset E$, we denote by $d_{P_U}(x, C): =\inf \{P_U(x-y): y \in C\}$ for the distance of $x$ and $C$ by the seminorm $P_U$, where $P_U$ is the Minkowski $p$-functional for each $U \in \mathfrak{U}$ in $E$.

\vskip.1in
\section{The KKM Principle in Convex Vector Spaces}

Since Knaster, Kuratowski and Mazurkiewicz (in short, KKM)\cite{kkm} in 1929 obtained
the so-called KKM principle (theorem) to give a new proof for the Brouwer fixed point theorem in finite dimensional spaces, and  later in 1961, Fan \cite{fan1960} (see also Fan \cite{fan1972}) extended the KKM principle (theorem) to any topological vector spaces and applied it to various results including the Schauder fixed point theorem, then there have appeared
a large number of works devoting applications of the KKM principle (theorem). In 1992, such research field was
called the KKM theory first time by Park \cite{park1991}, then the KKM theory has been extended to general abstract convex
spaces by Park \cite{park2010}(see also Park \cite{park2016} and \cite{park100}) which actually include locally $p$-convex spaces ($0 < p \leq 1)$ as a special class. The same as last section, for the convenience of self-reading, we recall some notion and definitions for $KKM$ principle in convex vector spaces which include $p$-vector space as a special class as summarized by Yuan \cite{yuan2022} below.

The same as last section, for the convenience of readers' self-containing reading, we recall again give some notion and definition on the abstract convex spaces which play important role for the development of KKM principle and related applications. Once again, the corresponding comprehensive discussion on KKM theory and its various applications to nonlinear analysis and related topics, we refer to Mauldin \cite{Mauldin}, Granas and Dugundji \cite{granas}, Park \cite{park100} and \cite{park2022}, Yuan \cite{yuan1999}-\cite{yuan2022} and related comprehensive reference there.

Let $\langle D\rangle$ denote the set of all nonempty finite subsets of a given non-empty set $D$, and $2^D$ denotes the family of all subsets of $D$.  We have the following definition for abstract convex spaces essentially by Park \cite{park2010}.

\noindent
{\bf Definition 3.1.} An abstract convex space $(E, D; \Gamma)$ consists of a
topological space $E$, a nonempty set $D$, and a set-valued mapping $\Gamma: \langle D\rangle\rightarrow 2^E$ with nonempty values $\Gamma_{A}: = \Gamma(A)$ for each $A \in \langle D\rangle$, such that the $\Gamma$-convex hull of any $D' \subset D$  is denoted and defined by $\mathrm{c}\mathrm{o}_{\Gamma}D': = \cup\{\Gamma_{A}| A \in \langle D'\rangle\}\subset E.$

A subset $X$ of $E$ is said to be a $\Gamma$-convex subset of $(E, D; \Gamma)$ relative to $D' $ if for any $N \in \langle D' \rangle$, we have $\Gamma_{N} \subseteq X$, that is, $\mathrm{c}\mathrm{o}_{\Gamma}D'\subset X$.
For the convenience of our discussion, in the case  $E=D$, the space $(E, E; \Gamma)$ is simply denoted by $(E; \Gamma)$ unless specified.

\noindent
{\bf Definition 3.2.} Let $(E, D; \Gamma)$ be an abstract convex space and $Z$ a topological space.
For a set-valued mapping (or say, multimap) $F:  E \rightarrow 2^Z$ with nonempty values, if a set-value mapping $G: D\rightarrow 2^Z$ satisfies
$F(\Gamma_{A})\displaystyle \subset G(A):=\bigcup_{y\in A}G(y)$ for all $A\in \langle D \rangle$, then $G$ is called a KKM mapping with respect to $F$. A KKM mapping $G: D\rightarrow 2^E$  is a KKM mapping with respect to the identity map $1_{E}$.

\noindent
{\bf Definition 3.3.} The partial KKM principle for an abstract convex space $(E, D; \Gamma)$ is
that, for any closed-valued KKM mapping $G: D\rightarrow 2^E$, the family $\{G(y)\}_{y\in D}$ has the finite intersection property. The KKM principle is that, the same property also holds for any open-valued KKM mapping.

An abstract convex space is called a (partial) KKM space if it satisfies the (partial) KKM principle (resp.). We now gave some known examples of (partial) KKM spaces (see Park \cite{park2010}, and also \cite{park2016}) as follows.

\noindent
{\bf Definition 3.4.} A $\phi_{A}$-space $(X, D;\{\phi_{A}\}_{A\in \langle D\rangle})$ consists of a
topological space $X$, a nonempty set $D$, and a family of continuous functions $\phi_{A}: \Delta_{n}\rightarrow 2^X$ (that is,
singular $n$-simplices) for $A \in \{D\}$ with $|A|=n+1$. By putting $\Gamma_{A}: = \phi_{A}(\Delta_{n})$ for
each $A\in \langle D \rangle$, the triple $(X, D; \Gamma)$ becomes an abstract convex space.

\noindent
{\bf Remark 3.1.} For a $\phi_{A}$-space $(X, D;\{\phi_{A}\})$ , we see easily that any set-valued mapping $G: D\rightarrow 2^X$ satisfying
$\phi_{A}(\Delta_{J})\subset G(J)$ for each $ A \in \langle D \rangle $  and  $J \in \langle A \rangle$ is a KKM mapping.

By the definition, it is clear that every $\phi_{A}$-space is a KKM space, thus we have the following fact (see Lemma 1 of Park \cite{park2016}).

\noindent
{\bf Lemma 3.1.}  Let $(X, D; \Gamma)$  be a $\phi_{A}$-space and $G:  D \rightarrow 2^X$ a set-valued (multimap) with
nonempty closed [resp. open] values. Suppose that $G$ is a KKM mapping, then $\{G(a)\}_{a\in D}$ has the finite intersection property.

\vskip.1in
By following Definition 2.7, we recall that a topological vector space is said to be locally
$p$-convex if the origin has a fundamental set of absolutely $p$-convex $0$-neighborhoods.
This topology can be determined by $p$-seminorms which are defined in the obvious
way (see Jarchow \cite{jarchow}, or P.52 of Bayoumi \cite{bayoumi}).

Now we have a new KKM space as follows inducted by the concept of $p$-convexity (see Lemma 2 of Park \cite{park2016}).

\noindent
{\bf Lemma 3.2.} Suppose that $X$  is a subset of  topological vector space $E$ and $p \in (0,1]$, and $D$ is a
nonempty subset of $X$ such that $C_{p}(D)\subset X$. Let $\Gamma_{N}: =C_{p}(N)$ for each $N\in \langle D\rangle$.
Then $(X, D; \Gamma)$ is a $\phi_{A}$-space.

\vskip.1in
\noindent
{\bf Proof.} Since $C_{p}(D)\subset X$, $\Gamma_{N}$ is well-defined. For each  $N=\{x_{0}, x_{1}, \cdots, x_{n}\}\subset D$, we define $\phi_{N}: \Delta_{n}\rightarrow \Gamma_{N}$ by
$\sum_{i=0}^{n}t_{i}e_{i}\mapsto\sum_{i=0}^{n}(t_{i})^{\frac{1}{\mathrm{p}}}x_{i}$.  Then clearly $(X, D; \Gamma)$  is a $\phi_{A}$-space. This completes the proof.  $\square$

\vskip.1in
\section{Fixed Point Theorems for Condensing Mappings in $p$-Vector Spaces}
\label{sec:1}

In this section, we will establish fixed point theorems for upper semi-continuous and condensing mappings for $p$-convex subsets under the general framework of $p$-vector spaces, which will be a tool used in section 5 and section 6 to establish the best approximation, fixed points, the principle of nonlinear alternative, Birkhoff-Kellogg problems, Leray-Schauder alternative which would be useful tools in nonlinear analysis for the study of  nonlinear problems arising from theory to the practice.
Here, we first gather together necessary definitions, notations, and known facts needed in this section.

\noindent
{\bf Definition 4.1.} Let $X$ and $Y$ be two topological spaces.
A set-valued mapping (also saying, multifunction) $T: X \longrightarrow 2^Y$ is a point to set function such
that for each $x \in X$, $T(x)$ is a subset of $Y$. The mapping $T$ is said to be upper semi-continuous (USC)
if the subset $T^{-1}(B): = \{ x\in X: T(x) \cap B \neq \emptyset\}$ (resp., the set $\{x \in X: T(x) \subset B\}$)  is closed (resp., open)  for any closed (resp., open) subset $B$ in $Y$.
The function  $T: X \rightarrow 2^Y$  is said to be lower semi-continuous (LSC) if the set
 $T^{-1}(A)$ is open for any open subset $A$ in $Y$.

As an application of KKM principle for general abstract convex spaces with the help of embedding lemma for Hausdorff compact $p$-convex subsets from topological vector spaces (TVS) into locally $p$-convex vector spaces, we have the following general existence result for the $``$approximation" of fixed points for upper and lower semi-continuous set-valued mappings in $p$-convex vector spaces for $0 < p  \leq 1$ (see  the corresponding related results given by Theorem 2.7 of Gholizadeh et al. \cite{gholizadeh}, Theorem 5 of Park \cite{park2016} and related discussion therein).

The following result is originally given by given by Yuan \cite{yuan2022}, here we provide the sketch of its proof for the purpose of reading's self-containing.

\vskip.1in
\noindent
{\bf Theorem 4.1.} Let $A$ be a $p$-convex compact subset of a locally $p$-convex vector space $X$, where $0 < p \leq 1$.
Suppose that $T: A \rightarrow 2^A$ is lower (resp. upper) semi-continuous with non-empty $p$-convex values.
Then for any given $U$ which is a $p$-convex neighborhood of zero in $X$, there exists $x_U \in A$ such that
$T(x_U) \cap (x_U + U) \neq \emptyset$.

\noindent
{\bf Proof.} Suppose $U$ is any given element of $\mathfrak{U}$, there is a symmetric open neighborhood $V$ of zero for which
$\overline{V} + \overline{V} \subset U$ in locally $p$-convex neighborhood of zero, we prove the results by two cases for $T$ is lower semicontinuous (LSC) and upper semicontinuous (USC).

Case 1, by assuming $T$ is lower semi-continuous: As $X$ is locally $p$-convex vector space, suppose that $\mathfrak{U}$ is the family of neighborhoods of $0$ in $X$. For any element $U$ of $\mathfrak{U}$, there is a symmetric open neighborhood $V$ of zero  for which
$\overline{V} + \overline{V} \subset U$.
Since $A$ is compact, so there exist $x_0, x_1, \cdots, x_n$ in $A$ such that $A \subset \cup_{i=0}^n (x_i + V)$. By using the fact that $A$ is $p$-convex, we find $D: =\{b_0, b_2, \cdots, b_n\} \subset A$ for which
$b_i - x_i \in V$ for all $i \in \{0, 1, \cdots, n\}$ and we define $C$ by $C: = C_p(D) \subset A$. By the fact that $T$ is LSC, it follows that the subset $F(b_i): = \{c \in C:  T(c) \cap (x_i +V) = \emptyset\}$ is closed in $C$ (as the set $x_i +V$ is open)  for each  $i \in \{0, 1, \cdots, n\}$. For any $c \in C$, we have $\emptyset \neq T(c)\cap A \subset T(c)\cap \cup_{i=0}^n(x_i+ V)$, it follows that $\cap_{i=0}^n F(b_i)=\emptyset$. Now applying Lemma 3.1 and Lemma 3.2, which implies that
that there is $N:= \{b_{i_0},  b_{i_1}, \cdots, b_{i_k}\} \in \langle D \rangle$
and $x_U \in C_p(N) \subset A$ for which $x_U \notin F(N)$, and so  $T(x_u) \cap (x_{i_j} + V) \neq \emptyset$ for all
$ j \in \{0, 1, \cdots, k\}$.
 As $ b_i  -  x_i \in V$ and $\overline{V} + \overline{V} \subset U$, which imply that $x_{i_j} + \overline{V} \subset b_{i_j} + U$, which means that $T(x_U) \cap ((b_{i_j} + U) \neq \emptyset$, it follows that
 $N \subset \{c \in C: T(x_U) \cap (c + U)\neq \emptyset \}$.
By the fact that the subsets $C, T(x_U)$ and $U$ are $p$-convex, we have that $x_U \in \{c \in C: T(x_U) \cap (c+U)\neq \emptyset\}$, which means that $T( x_U) \cap (x_U  + U ) \neq \emptyset$.

Case 2, by assuming $T$ is upper semi-continuous:  We define $F(b_i): = \{c \in C: T(c) \cap (x_i + \overline{V}) = \emptyset\}$, which is then open in $C$ (as the subset $x_i + \overline{V}$  is closed) for each $i=0, 1, \cdots, n$. Then the argument is similar to the proof for the case T is USC, and by applying Lemma 3.1 and Lemma 3.2 again, it follows that there exists $x_U \in A$ such that $T(x_U) \cap (x_U + U) \neq \emptyset$. This completes the proof.  $\square$

\vskip.1in
\noindent
By Theorem 4.1, we have the following Fan-Glicksberg fixed point theorems (Fan \cite{fan1952}) in locally $p$-convex vector spaces for $(0 < p \leq 1)$, which also improve or generalize the corresponding results given by Yuan \cite{yuan1999}, Xiao and Lu \cite{xiaolu}, Xiao and Zhu \cite{xiaozhu2011}-\cite{xiaozhu2018} into locally $p$-convex vector spaces.

\vskip.1in
\noindent
{\bf Theorem 4.2.} Let $A$ be a $p$-convex compact subset of a locally $p$-convex vector space $X$, where $0 < p \leq 1$.
Suppose that $T: A \rightarrow 2^A$ is upper semi-continuous with non-empty $p$-convex closed values.
Then $T$ has at least one fixed point.

\vskip.1in
\noindent
{\bf Proof.} Assume $\mathfrak{U}$ is the family of open $p$-convex neighborhoods of $0$ in $X$, and $U \in \mathfrak{U}$, by Theorem 4.1, there exists $x_U \in A$ such that $T(x_U) \cap (x_U + U) \neq \emptyset$.
Then there exists $a_U, b_U \in A$ for which $b_U \in T(a_U)$ and $b_U \in a_U + U$.
Now, two nets $\{a_U\}$ and $\{b_U\}$ in Graph$(T)$, which is a compact graph of mapping $T$ as $A$ is compact and $T$ is semi-continuous, we may assume that $a_U$ has a subnet converging to $a$, and $\{b_U\}$ has a subnet converging to $b$. As $\mathfrak{U}$ is the family of neighborhoods for $0$, we should have $a=b$ (e.g., by the Hausdorff separation property), and $a=b \in T(b)$ due to the fact that Graph(T) is close (e.g., see Lemma 3.1.1 in P.40 of Yuan \cite{yuan1998}), thus the proof is compete. $\square$

For a given set $A$ in vector space $X$, we denote by $``$$lin(A)$" the $``$linear hull" of $A$ in $X$.

\vskip.1in
\noindent
{\bf Definition 4.2.} Let $A$ be a subset of a topological vector space $X$ and let $Y$ be another topological vector space. We shall say that $A$ can be linearly embedded in $Y$ if there is a linear map $L: lin(A) \rightarrow Y$ (not necessarily
continuous) whose restriction to $A$ is a homeomorphism.

\vskip.1in
The following embedded Lemma 4.1 is a significant result due to Theorem 1 of Kalton \cite{kalton1977}, which says though not every compact convex set can be linearly imbedded in a locally convex space (e.g., see Roberts \cite{roberts1977} and
 Kalton et al.\cite{kalton1984}), but for $p$-convex sets when $0 < p <1$, every compact $p$-convex set in topological vector spaces can be considered as a subset of a locally $p$-convex vector space, hence every such set has sufficiently many $p$-extreme points.

Secondly, by the property (ii) of Lemma 2.1 above, each convex subset of a topological vector space containing zero is always $p$-convex for $0 < p \leq 1$, thus it is possible for us to transfer the problem involved $p$-convex subsets from topological vector spaces into the locally $p$-convex vector spaces, which indeed allows us to establish the existence of fixed points for upper semi-continuous set-valued mappings for compact  $p$-convex subsets in locally convex spaces for $ 0 < p \leq 1$, but we note that by Lemma 2.4, any set-valued mapping with closed $p$-convex values in $p$-spaces for $0< p < 1$ has the zero point as its trivial fixed point, thus it is essentially to study fixed point and related principle of nonlinear analysis for single-valued (instead of set-valued) mappings in $p$-vector spaces as pointed out by Remark 2.5 (see also the discussion in P.40-41 given by Yuan \cite{yuan2022}).

Indeed, a fixed point theorem a topological vector space for (single-valued) continuous and condensing mappings given by Theorem 4.5 below which will be proved below (also see Theorem 4.3 below essentially due to Ennassik and Taoudi \cite{ennassik2021}) provides the answer for Schauder's  conjecture in the affirmative.

\vskip.1in
\noindent
{\bf Lemma 4.1.} Let $K$ be a compact $p$-convex subset $(0 <
p < 1$) of a topological vector space $X$. Then, $K$ can be linearly
embedded in a locally $p$-convex topological vector space.

\noindent
{\bf Proof.} It is Theorem 1 of Kalton \cite{kalton1977}, which completes the proof. $\square$

\vskip.1in
\noindent
{\bf Remark 4.1.} At this point, it is important to note that Lemma 4.1 does
not hold for $p = 1$. By Theorem 9.6 of Kalton et al.\cite{kalton1984}, it was shown that the spaces
$L_p = L_p(0, 1)$, where $0 < p < 1$, contain compact convex sets with no extreme points, which thus  cannot be linearly embedded in a locally convex space, see also Roberts \cite{roberts1977}.

Now we give the the following fixed point theorem for single-valued continuity mappings, which are essentially Theorem 3.1 and Theorem 3.3 given first by Ennassik and Taoudi \cite{ennassik2021}. Here we include the argument for the second part of the conclusions below only.

\vskip.1in
\noindent
{\bf Theorem 4.3.} If $K$  is a nonempty compact $p$-convex subset of a locally $p$-convex space $E$ for $ 0 < p \leq 1$, then the (single-valued) continuous mapping $T: K \rightarrow K$ has at least a fixed point. Secondly, if $K$  is a nonempty compact $p$-convex subset of a Hausdorff TVS $E$, then the (single-valued) continuous mapping $T: K \rightarrow K$ has at least a fixed point.

\noindent
{\bf Proof.} The first part is Theorem 3.1 of Ennassik and Taoudi \cite{ennassik2021}, and the second part is indeed Theorem 3.3 of Ennassik and Taoudi \cite{ennassik2021}, but here we include their very smart proof as below.

Case 1: For $0< p < 1$, $K$ is a nonempty compact $p$-convex subset of a
topological vector space $X$ for $ 0 < p < 1$, by Lemma 4.1, it follows that $K$ can be linearly embedded in a locally $p$-convex
space $E$, which means that there exists a linear map
$L: \operatorname{lin}(K) \rightarrow E$ whose restriction to $K$ is a homeomorphism.
Define the mapping $S: L(K) \rightarrow L(K)$ by $(Sx): = L(Tx)$ for
$x \in X$. This mapping is easily checked to be well defined. The mapping $S$ is continuous since $L$ is a (continuous) homeomorphism and $T$ is continuous on $K$. Furthermore,
the set $L(K)$ is compact, being the image of a compact set under a continuous
mapping $L$. It is also $p$-convex since it is the image of a $p$-convex set
under a linear mapping. Then, by the conclusion in the first part (see also Theorem 3.1 in \cite{ennassik2021}), there exists $x \in K$ such
that $Lx =S(Lx) = L(Tx)$, thus it implies that $x =T(x)$ since $L$ is a homeomorphism, which is the fixed point of $T$.

Case 2: For $p=1$, taking any point $x_0 \in K$, and let $K_0: =K-\{x_0 \}$. Now define a new mapping $T_0: K_0 \rightarrow K_0$ by
$T_0(x)= T(x)-x_0$ for each $x \in K_0$. By the fact that now $K_0$ is $p$-convex for any $0< p < 1$ by Lemma 2.1(ii), then the $T_0$ has a fixed point in $K_0$ by the proof in Case 1 above, so $T$ has a fixed point in $K$. The proof is complete.
$\square$

\vskip.1in
{\bf Remark 4.2}
Theorem 4.3 is indeed the results of Theorem 3.1 and Theorem 3.3 (of \cite{ennassik2021}) for $0 < p \leq 1$ which provides an answer to Schauder's
conjecture under the TVS. Here we also mention a number of related works and discussion by authors in this drection,
see Mauldin \cite{Mauldin}, Granas and Dugundji \cite{granas}, 
Park \cite{park100, park2022} and the references therein.
\vskip.1in
We recall that for two given topological spaces $X$ and $Y$, and a set-valued mapping $T: X \rightarrow 2^Y$ is said to be compact if there is compact subset set $C$ in $Y$ such that $F(X) (=\{y \in F(X), x \in X\})$ is contained in $C$, i.e., $F(X) \subset C$. Now we  have the following non-compact version of fixed point theorems for compact set-valued mappings defined on a general $p$-convex subset in  $p$-vector spaces for $ 0 < p \leq 1$.

\vskip.in
As an immediate consequence of Theorem 4.2 for $p=1$, we have following result for upper semi-continuous version in locally convex spaces (LCS).
\vskip.1in
\noindent
{\bf Theorem 4.4}  If $K$ is a nonempty compact convex subset of a locally convex space $X$, then any upper semi-continuous set-valued mappings $T: K \rightarrow 2^K$ with non-empty closed convex values has at least a fixed point.

\noindent
{\bf Proof.} Apply Theorem 4.2 with $p = 1$, this completes the proof. $\square$

\vskip.1in
Theorem 4.4 also improves or unifies corresponding results given by Askoura and Godet-Thobie \cite{Askoura}, Cauty \cite{cauty2005}, Cauty \cite{cauty2007}, Chen \cite{chen}, Isac \cite{isac}, Li \cite{li}, Nhu \cite{nhu1996}, Okon \cite{okon}, Park \cite{park2022}, Reich \cite{reich}, Smart \cite{smart}, Yuan \cite{yuan1999}, Theorem 3.14 of Gholizadeh et al.\cite{gholizadeh}, Xiao and Lu \cite{xiaolu}, Xiao and Zhu \cite{xiaozhu2011}-\cite{xiaozhu2018}  under the framework of LCS for set-valued mappings instead of single-valued functions.

\vskip.1in
In order to establish fixed point theorems for the classes of 1-set contractive and condensing mappings in $p$-vector spaces by using the concept of the measure of noncompactness (or saying, the noncompactness measures) which were introduced and widely accepted in mathematical community by Kuratowski \cite{kuratowski}, Darbo \cite{darbo} and related references therein, we first need to have a brief introduction for the concept of non-compactness measures for the so-called Kuratowski or Hausdorff measures of noncompactness in normed spaces (see Alghamdi et al.\cite{alghamdiporegan}, Machrafi and Oubbi \cite{machrafioubbi}, Nussbaum \cite{nussbaum},  Sadovskii \cite{sadovskii}, Silva et al.\cite{silva}, Xiao and Lu \cite{xiaolu} for the general concepts under the framework of $p$-seminorm or, just for locally convex $p$-convex settings for $0< p \leq 1$ for which will be discussed below, too).

For a given metric space $(X, d)$ (or a $p$-normed space $(X, \|\cdot\|_p)$), we recall notions of completeness, boundedness, relatively compactness and compactness as follows. Let $(X, d)$ and $(Y, d)$ be two metric spaces and $T: X \rightarrow Y$ is a mapping (or saying, operator). Then:  1) $T$ is said to be bounded if for each bounded set $A\subset X$, $T(A)$ is  bounded
set of $Y$; 2) $T$  is said to be continuous if for every $x \in X$, the $\lim_{n \rightarrow \infty} x_n = x$  implies
that $\lim_{n\rightarrow \infty} T(x_n)= T$; and 3) $T$ is said to be completely continuous if $T$ is continuous and $T(A)$ is relatively compact for each bounded subset $A$ of $X$.

Let $A_1$, $A_2 \subset X$ be bounded of a metric space $(X, d)$, we also recall that the Hausdorff metric $d_H(A_1, A_2)$ between $A_1$ and $A_2$ is defined by
$$d_H(A_1, A_2): =\max\{\sup_{x\in A_1}\inf_{y \in A_2} d(x, y), \sup_{y \in A_2} \inf_{x \in A_1} d(x, y)\}$$
The Hausdorff and Kurotowskii measures of noncompactness (denoted by $\beta_H$ and $\beta_K$, respectively) for nonempty bounded subset $D$ in $X$ are are the nonnegative real numbers $\beta_H(D)$ and $\beta_K(D)$ defined by
$$\beta_H(D): = \inf\{\epsilon > 0: D \mbox{ has a finite} \epsilon\mbox{-net}\},$$
and
$$\beta_K(D): = \inf\{ \epsilon > 0: D\subset \cup_{i=1}^n D_i, \mbox{ where $D_i$ is bounded and $diam D_i \leq  \epsilon$, $n$ is an integer } \},$$
here $diam D_i$ means the diameter of the set $D_i$, and it is well known that $\beta_H \leq \beta_K \leq 2 \beta_H$.
We also point out that the notions above can be well defined under the framework of $p$-seminorm spaces $(E, \|\cdot\|_p)_{p \in \mathfrak{P}}$ by following the similar idea and method used by Chen and Singh \cite{chensingh}, Ko and Tasi \cite{kotsai}, Kozlov et al.\cite{kozlovetal} and references therein  for more in details.

Let $T$ is a mapping from $D\subset X$  to $X$. Then  we have that: 1) $T$ is said to be a $k$-set contraction with respect to $\beta_K$ (or $\beta_H$) if there is a number $k \in (0, 1]$  such that $\beta_K(T(A)) \leq k \beta_K(A)$ (or $\beta_H(T(A)) \leq k\beta_H(A))$ for all bounded sets $A$ in $D$; and 2) $T$ is said  said to be $\beta_K$-condensing (or $\beta_H$-condensing) if $(\beta_K(T(A)) < \beta_K(A))$ (or $\beta_H (T(A)) < \beta_H(A))$  for all bounded sets $A$ in $D$ with $\beta_K(A)> 0$ (or $\beta_H(A)> 0$).

For the convenience of our discussion, throughout the rest part of this paper, if a mapping $``$is $\beta_K$-condensing (or $\beta_H$-condensing)", we simply say it is $``$a condensing mapping" unless specified.

Moreover, it is easy to see that: 1) if $T$ is a compact operator, then $T$ is a $k$-set contraction; and 2) if $T $ is a $k$-set contraction for $k \in (0, 1)$, then $T$ is condensing.

In order to establish the fixed points of set-valued condensing mappings in $p$-vector spaces for $0 < p \leq 1$, we need to
recall some notions introduced by Machrafi and Oubbi \cite{machrafioubbi}  for the measure of noncompactness in locally $p$-convex vector spaces, which also satisfies some necessary (common) properties of the classical measures of noncompactness such as $\beta_K$ and $\beta_H$ mentioned above introduced by Kuratowski \cite{kuratowski}, Sadovskii \cite{sadovskii}(see, also related discussion by Alghamdi et al.\cite{alghamdiporegan}, Nussbaum \cite{nussbaum}, Silva et al.\cite{silva}, Xiao and Lu \cite{xiaolu} and references therein). In particular, the measures of noncompactness in locally $p$-vector spaces (for $0 < p \leq 1$) should have the stable property which means the measure of noncompactness $A$ is the same by transition to the (closure) for the $p$-convex hull of subset $A$.

For the convenience of discussion, we follow up to use $\alpha$ and $\beta$  to denote the Kuratowski and the Hausdorf measures of noncompactness in topological vector spaces, respectively (see the same way used by  Machrafi and Oubbi \cite{machrafioubbi}),  unless otherwise stated. The $E$ is used to denote a Hausdorf topological vector space over the field $\mathbb{K} \in \{\mathbb{R}, \mathbb{Q}\}$, here $\mathbb{R}$ denotes for all real numbers, and $\mathbb{Q}$ for all complex numbers, and $p \in (0, 1]$. Here, the base set of family of all balanced zero neighborhoods in $E$ is denoted by $\mathfrak{V}_0$.

We recall that $U \in \mathfrak{V}_0$ is said to be shrinkable, if it is absorbing, balanced, and $r U \subset U$ for all $r \in (0, 1)$, and we know that any topological vector space admits a local base at zero consisting of shrinkable sets (see Klee \cite{klee1960}, or Jarchow \cite{jarchow} for details).

Recalling again that a topological vector space $E$ is said to be a locally $p$-convex space, if $E$ has a
local base at zero consisting of $p$-convex sets. The topology of a locally $p$-convex space is always given by an upward directed family $P$ of $p$-semi-norms, where a $p$-semi-norm on $E$ is any non negative real-valued and subadditive functional $\|\cdot\|_p$ on on $E$ such that  $\| \lambda x\|_p=|\lambda|^p\|x\|_p$ for each $x \in E$ and $\lambda \in \mathbb{R}$ (i.e., the real number line). When $E$ is Hausdorff, then for every $x \neq 0$,  there is some $p \in P$ such that $P(x) \neq 0$.
Whenever the family $P$ is reduced to a singleton, one says that  $(E, \| \cdot \|)$ is a $p$-semi-normed space.
A $p$-normed space is a Hausdorff $p$-seminormed space, and when $p=1$ which is the usual locally convex case. Furthermore, a $p$-normed space is a metric vector space with the translation invariant metric $d_p(x, y): = \| x- y\|_p$ for all $x, y \in E$, which is the same notation used above.

By Remark 2.5, if $P$ is a continuous $p$-seminorm on $E$, then the ball $B_p(0, s): = \{x \in E: P(x) < s  \}$ is
shrinkable for each $s > 0$. Indeed, if $r \in (0, 1)$ and $x \in \overline{r B_p(0, s)}$,  then there exists
 a net $(x_i)_{i \in I}  \subset B_p(0, s)$  such that $r x_i$  converges to $x$.  By continuity of $P$,  we get
$P(x) \leq r^p s < s$, which means that $r \overline{B_p(0,s)} \subset B_P(0,s)$. In generally, it can be shown
that every $p$-convex $U \in \mathfrak{V}_0$ is shrinkable.

We recall that a given such neighborhood $U$, a subset $A \subset E$  is said to be
$U$-small if  $A - A \subset U$ (or saying, small of order $U$ by Robertson \cite{robertson}). Now by following the idea of
Kaniok \cite{kaniok} in the setting of a topological vector space $E$ to use zero neighborhoods in $E$  instead of semi-norms to to define the measure of noncompactness in (local convex) $p$-vector spaces ($0< p \leq 1$) as follows: For each $A \subset E$, the $U$-measures of noncompactness $\alpha_U(A)$ and $\beta_U(A)$ for $A$ are defined by:
$$\alpha_U(A): =\inf\{: r > 0: \mbox{$A$ is covered by a finite number of $rU$-small sets $A_i$ for $i=1, 2, \cdots, n$} \},$$
and
$$\beta_U(A): = \inf\{ r > 0: \mbox{there exists $x_1, \cdots, x_n \in E$ such that $A \subset \cup_{i=1}^n (x_i+ rU)$ }\},$$
here we set $\inf \emptyset: = \infty$.

By the definition above, it is clear that when $E$ is a normed space and $U$ is the closed unit ball of $E$,
$\alpha_U$ and $\beta_U$ are nothing else but the Kuratowski measure $\beta_K$ and Hausdorf measure $\beta_H$ of noncompactness, respectively. Thus, if $\mathfrak{U}$ denotes a fundamental system of balanced and closed zero neighborhoods in $E$ and
$\mathfrak{F}_{\mathfrak{U}}$ is the space of all functions $\phi: \mathfrak{U} \rightarrow R$, endowed with the pointwise ordering, then, the $\alpha_{U}$ (resp., $\beta_{U}$) measures for noncompactness introduced by Kaniok \cite{kaniok} can be expressed by the Kuratowski (resp., the Hausdorf) measure of  noncompact $\alpha(A)$(resp., $\beta(A)$) for a subset $A$ of $E$ as the function defined from $\mathfrak{U}$ into
$[0, \infty)$ by
$$\alpha(A)(U): = \alpha_{U}(A)  \mbox{ (resp.,  } \beta(A)(U):=\beta_{U}(A)).$$

By following Machrafi and Oubbi \cite{machrafioubbi}, in order to define the measure of noncompactness in (locally convex) $p$-vector space $E$, we need the following notions of basic and sufficient collections for zero neighborhoods in a topological vector space. To do this, let us introduce an equivalence relation on $V_0$ by saying that $U$ is related to $V$, written
$U\mathfrak{R}V$, if and only if there exist $r, s > 0$  such that $r U \subset V \subset s U$. We now have the following definition.

\vskip.1in
\noindent
{\bf Definition 4.2 (BCZN)}. We say that $\mathfrak{B} \subset \mathfrak{V}_0$ is a basic collection of zero neighborhoods (in short, BCZN) if it contains at most one representative member from each equivalence class with respect to $\mathfrak{R}$. It will be said to be sufficient (in short, SCZN) if it is basic
and, for every $V \in \mathfrak{V}_0$, there exists some $U \in \mathfrak{B}$ and some $r > 0 $ such that  $r U \subset V$.

\vskip.1in
\noindent
{\bf Remark 4.3}. By Remark 2.5, it follows that for a locally $p$-convex space $E$, its base set $\mathfrak{U}$, the family of all open $p$-convex subsets for $0$ is BCZB. We also note that: 1) In the case if $E$ is a normed space, if $f$ is a continuous functional on $E$, $U: =\{x \in E: |f(x)| < 1\}$, and $V$ is the open unit ball of $E$, then  $\{U\}$ is basic but not sufficient, but $\{V\}$ is sufficient;
2) Secondly, if $(E, \tau)$ is a locally convex space, whose topology is given by an upward directed
family $P$ of seminorms, so that no two of them are equivalent, the collection $(B_p)_{p \in \mathbb{P}}$ is a SCZN, where $B_p$ is the open unit ball of $p$. Further, if $\mathfrak{W}$ is a fundamental system of zero neighborhoods in a topological
vector space $E$, then there exists an SCZN consisting of $\mathfrak{W}$ members; and 3) By following Oubbi \cite{oubbi}, we recall that a subset $A$ of $E$  is called uniformly bounded with respect to a sufficient collection $\mathfrak{B}$ of zero neighborhoods, if there exists $r > 0 $ such that $A \subset r V$ for all $V \in \mathfrak{B}$. Note that in the locally convex space $C_c(X): = C_c(X, \mathbb{K})$, the
set $B_{\infty}:=\{ f\in C(X): \|f\|_{\infty} \leq 1\}$ is uniformly bounded with respect to the SCZN
$\{B_k, k \in \mathbb{K}\}$, where $B_k$ is the (closed or) open unit ball of the seminorm $P_k$, where $k \in \mathbb{K}$.

\vskip.1in
Now we are ready to give the definition for the measure of non-compactness  in (locally $p$-convex) topological
vector space $E$ as follows.

\vskip.1in
\noindent
{\bf Definition 4.3.} Let $\mathfrak{B}$ be a SCZN in $E$. For each $A \subset E$, we define the measure of noncompactness of $A$ with respect to $\mathfrak{B}$ by $\alpha_{\mathfrak{B}}(A):=\sup_{U\in \mathfrak{B}}\alpha_{U}(A)$.

\vskip.1in
By the definition above, it is clear that:  1) The measure of noncompactness $\mathfrak{B}$ holding the semi-additivity, i.e.,
$\alpha_{\mathfrak{B}}(A \cup B) = \max\{\alpha_{\mathfrak{B}}(A), \alpha_{\mathfrak{B}}(B)\}$; and 2) $\alpha_{\mathfrak{B}}(A) = 0 $ if and only if $A$ is a precompact subset of $E$ (for more properties in details, see Proposition 1 and related discussion by Machraf and Oubbi \cite{oubbi}).

\vskip.1in
As we know, under the normed spaces (and even semi-normed spaces), Kuratowski \cite{kuratowski}, Darbo \cite{darbo} and Sadovskii \cite{sadovskii} introduced the notions of $k$-set-contractions for $k \in (0, 1)$, and the condensing mappings to establish fixed point theorems in the setting of Banach spaces, normed or semi-norm spaces. By following the same idea, if $E$ is a Hausdorf locally $p$-convex space, we have the following definition for general (nonlinear) mappings.

\vskip.1in
\noindent
{\bf Definition 4.4.} A mapping $T: C \rightarrow 2^C$ is said to be a $k$-set contraction (resp., condensing),
if there is some SCZN $\mathfrak{B}$ in $E$ consisting of $p$-convex sets, such that (resp., condensing) for any $U \in \mathfrak{B}$, there exists $k \in (0,1)$ (resp., condensing) such that $\alpha_U(T(A)) \leq k \alpha_U(A)$ for $A \subset C$
(resp., $\alpha_U(T(A)) < \alpha_U(A)$ for each $A \subset C$ with $\alpha_U(A) > 0$).

\vskip.1in
It is clear that a contraction mapping on $C$  is a $k$-set contraction mapping (where we always mean $k \in (0, 1)$), and a
$k$-set contraction mapping on $C$ is condensing; and they all reduce to the usually cases by the definitions for the $\beta_K$ and $\beta_H$ which are the Kuratowski measures and Hausdorff measure of noncompactness, respectively in normed spaces (see Kuratowski \cite{kuratowski}).

From now on, we denote by $\mathfrak{V}_0$ the set of all shrinkable zero neighborhoods in $E$,
we have the following result which is Theorem 1 of Machrafi and Oubbi \cite{machrafioubbi}, saying that in the general setting of locally $p$-convex spaces, the measure of noncompactness $\alpha$ for $U$ given by Definition 4.3 above is stable from $U$ to its $p$-convex hull $C_p(A)$ of the subset $A$ in $E$, which is key for us to establish the fixed points for condensing mappings in locally $p$-convex spaces for $0< p \leq 1$. This also means that the key property for the measures due to the Kurotowski and Hausdorff measures of noncompactness in normed (or $p$-semi-norm) spaces, which also holds for the measure of noncompctness by Definition 4.3 in the setting of locally $p$-convex spaces with $(0 < p \leq 1)$ (see  more similar and related discussion  in details by Alghamdi et al.\cite{alghamdiporegan} and Silva et al.\cite{silva}).

\vskip.1in
\noindent
{\bf Lemma 4.2.} If $U \in \mathfrak{V}_0$ is $p$-convex for some $ 0 < p  \leq 1$, then $\alpha(C_p(A)) = \alpha(A)$  for
every $A \subset E$.

\noindent
{\bf Proof.} It is Theorem 1 of  Machrafi and Oubbi \cite{machrafioubbi}. The proof is complete. $\square$

\vskip.1in
Now based on the definition for the measure of noncompactness given by Defintion 4.3 (originally from Machrafi and Oubbi \cite{machrafioubbi}), we  have the following general extension version of Schauder, Darbo and Sadovskii type fixed point theorems in the context of locally $p$-convex vector spaces for condensing mappings.

\vskip.1in
\noindent
{\bf Theorem 4.5 (Schauder Fixed Point Theorem for single-valued condensing mappings).}
Let $C \subset E$  be a complete $p$-convex subset of a Hausdorf locally $p$-convex, or Hausdorf topological vector space space $E$, with $0 < p \leq 1$. If $T: C \rightarrow C$ is continuous and $(\alpha$) condensing, then $T$ has a fixed point in $C$ and the
set of fixed points of $T$ is compact.

\noindent
{\bf Proof.} We first prove the conclusion by assuming $E$ is a locally $p$-convex space, then prove the conclusion when $E$ is a topological vector space.

Case A: Assuming $E$ is locally $p$-convex. In this case, let $\mathfrak{B}$ be a sufficient collection of $p$-convex zero neighborhoods in $E$ with
respect to which $T$ is condensing and for any given $U \in \mathfrak{B}$. We choose some $x_0 \in C$ and let
$\mathfrak{F}$ be the family of all closed $p$-convex subsets $A$ of $C$ with $x_0 \in A$ and $T(A) \subset A$.
 Note that $\mathfrak{F}$ is not empty since $C \in \mathfrak{F}$. Let $A_0=\cap_{A \in \mathfrak{F}} A$.
 Then $A_0$ is a non empty closed $p$-convex subset of $C$, such that $T(A_0) \subset A_0$, and then the conclusion follows by Theorem 4.3 for the continuous mapping $T$ from $A_0$ to $A_0$.\, subject to show that $A_0$ is  compact. Now we prove $A_0$ is compact. Indeed, let
  $A_1=\overline{C_p(T(A_0) \cup \{x_0\})}$. Since $T(A_0)\subset A_0$ and $A_0$ is closed and p-convex, $A_1\subset A_0$.
 Hence, $T(A_1)\subset T(A_0)\subset A_1$. It follows that $A_1 \in \mathfrak{F}$ and therefore $A_1=A_0$.
 Now by Proposition 1 of Machrafi and Oubbi \cite{machrafioubbi} and Lemma 4.2 above (i.e., Theorem 1 and Theorem 2 in \cite{machrafioubbi}), we get $\alpha_U(T(A_0)) = \alpha_U(A)$. Our assumption on $T$ shows
that $\alpha_U(A_0)=0$ since $T$ is condensing. As $U$ is arbitrary from the family $\mathfrak{B}$, thus $A_0$ is $p$-convex and compact (see Proposition 4 in \cite{machrafioubbi}). Now, the conclusion follows by Theorem 4.3 above.
Secondly, let $C_0$ be the set of fixed points of $T$ in $C$. Then it follows that $C_0 \subset T(C_0)$ and the
the upper semi-continuity of $T$ implies that its graph is closed, so is the set $C_0$. As $T$ is condensing, we have $\alpha_U(T(C_0)) \leq \alpha_U(C_0)$, which implies that $\alpha_U(C_0)=0$.  As $U$ is arbitrary from the family $\mathfrak{B}$, which implies that $C_0$ is compact (by the Proposition 4 in \cite{machrafioubbi} again).

Case B: We now prove the conclusion by assuming $E$ is a topological vector space. Based on the argument in the Case A's proof above, when $T$ is condensing, there exists a non-empty compact $p$-convex subset $A_0$ such that $T: A_0 \rightarrow A_0$.
We prove the conclusion by considering two situations:  (1)  $0< p < 1$; and (2) $p=1$.

Now for  the case (1)  $0 <p< 1$: By the proof above, $A_0$ is a nonempty compact $p$-convex subset of a
topological vector space $E$, by Lemma 4.1, it follows that $A_0$ can be linearly embedded in a locally $p$-convex
space $X$, which means that there exists a linear mapping
$L: \operatorname{lin}(A_0) \rightarrow X$ whose restriction to $A_0$ is a homeomorphism.
Define the mapping $S: L(A_0) \rightarrow L(A_0)$ by $(Sx): = L(Tx)$ for
$x \in A_0$. This mapping is easily checked to be well defined. The mapping $S$ is continuous (and condensing) since $L$ is a (continuous) homeomorphism and $T$ is continuous (and condensing) on $A_0$. Furthermore,
the set $L(A_0)$ is compact, being the image of a compact set under a continuous
mapping $L$. It is also $p$-convex as it is the image of a $p$-convex set
under a linear mapping. Then, by the conclusion in the first part above for $S$ on $A_0$, there exists $x \in A_0$ such
that $Lx =S(Lx) = L(Tx)$, thus it implies that $x =T(x)$ since $L$ is a homeomorphism, which means $x$ is the fixed point of $T$.

Now for the case (2) $p=1$: taking any point $x_0 \in A_0$, and let $K_0: = A_0 -\{x_0 \}$. Now define a new mapping $T_0: K_0 \rightarrow K_0$ by $T_0(x)= T(x)-x_0$ for each $x \in A_0$. By the fact that now $K_0$ is $p$-convex for any $0< p < 1$ by Lemma 2.1(ii), then the $T_0$ has a fixed point in $K_0$ by the proof above for the case (1) when $0 < p < 1$,  so $T_0$ has a fixed point in $K_0$ implies that $T$ has a fixed point in $A_0$.

This completes the proof.
$\square$

\vskip.1in
\noindent
{\bf Remark 4.4} We first note that Theorem 4.5 improves Theorem 4.5 of Yuan \cite{yuan2022}. Secondly, as pointed by Remark 2.2 (for Theorem 3.1 and Theorem 3.3 given by Ennassik and Taoudi \cite{ennassik2021}),
Theorem 4.5 above provides an answer to Schauder's  conjecture in the affirmative way under the general framework of closed $p$-convex subsets in topological vector spaces for $ 0 < p \leq 1$ of (single-valued) continuous condensing mappings.
Here we also mention a number of related works and discussion by authors in this drection, see Mauldin [74], Granas and Dugundji [46], Park [90, 91] and the references therein.

\vskip.1in
By following the same argument used by Theorem 4.5, we have the following results for upper semicontinuous set-valued mappings in locally convex spaces as an application of Theorem 4.2.

\vskip.1in
\noindent
{\bf Theorem 4.6 (Schauder Fixed Point Theorem for upper semicontinous condensing mappings).}
Let $C$ be a convex subset of a locally convex space $E$. If $T: C \rightarrow 2^C$ is upper semicontinuous, $(\alpha$) condensing with closed convex values, then $T$ has a fixed point in $C$ and the set of fixed points of $T$ is compact.

\noindent
{\bf Proof.} By the same argument for Theorem 4.5 above by applying Theorem 4.4. $\square$

\vskip.1in
As applications of Theorem 4.5, we have the following fixed points for condensing mappings in locally $p$-convex, or topological vector spaces for $0 < p \leq 1$.

\noindent
{\bf Corollary 4.2 (Darbo type fixed point theorem).}
Let $C$ be a complete $p$-convex subset of a Hausdorf locally $p$-convex space or topological vector space $E$ with $0 < p \leq 1$.
If $T: C \rightarrow C$ is a (k)-set-contraction (where $k \in (0, 1))$, then $T$ has a fixed point.

\vskip.1in
\noindent
{\bf Corollary 4.3 (Sadovskii type fixed point theorem).}  Let $(E, \| \cdot \|)$ be a complete $p$-normed space and $C$  be a bounded, closed and p-convex subset of $E$, where $0 < p \leq 1$. Then, every continuous and condensing mapping $T: C \rightarrow C$ has a fixed point.

\noindent
{\bf Proof.} In Theorem 4.5, let $\mathfrak{B}: =\{B_p(0, 1) \}$, where $B_p(0,1)$ stands for the closed unit ball
of $E$, and by the fact that it is clear that $\alpha(A)=(\alpha_{\mathfrak{B}}(A))^p$ for each $A \subset E$.  Then  $T$ satisfies all conditions of Theorem 4.5. This completes the proof. $\square$

\vskip.1in
\noindent
{\bf Corollary 4.4 (Darbo type).} Let $(E, \| \cdot \|)$  be a complete $p$-normed
space and $C$ be a bounded, closed and $p$-convex subset of $E$, where $0 < p \leq 1$. Then each single-valued mapping $T: C \rightarrow C$  has a fixed point.
\vskip.1in
Theorem 4.5 and Theorem 4.6 improve Theorem 5 of  Machrafi and Oubbi \cite{machrafioubbi} for general condensing mappings, and also unify corresponding the results in the existing literature, e.g., see Alghamdi et al.\cite{alghamdiporegan}, G\'{o}rniewicz \cite{gorniewicz}, G\'{o}rniewicz et al.\cite{gorniewiczetal},
Nussbaum \cite{nussbaum}, Silva et al.\cite{silva}, Xiao and Lu \cite{xiaolu}, Xiao and Zhu \cite{xiaozhu2011}-\cite{xiaozhu2018} and references therein.

Before the ending of this section, we also like to remark that by comparing with topological method or related arguments used by Askoura et al.\cite{Askoura}, Cauty \cite{cauty2005}-\cite{cauty2007}, 
Nhu \cite{nhu1996}, Reich\cite{reich}, the fixed points given in this section improve or unify the corresponding ones  given by Alghamdi et al.\cite{alghamdiporegan}, Darbo \cite{darbo}, Liu\cite{liu2001}, Machrafi and Oubbi \cite{machrafioubbi}, Sadovskii \cite{sadovskii}, Silva et al.\cite{silva}, Xiao and Lu \cite{xiaolu} and those from references therein.

\vskip.1in
\section{Best Approximation for the Class of 1-Set Contractive Mappings in Locally $p$-Convex Spaces}

The goal of this section is first to establish one general best approximation results for the classes of  1-set continuous and hemicompact (see its definition below) non-self mappings, which in turn are used as a tool to  derive the general principle for the existence of solutions for Birkhoff-Kellogg Problems (see Birkhoff and Kellogg \cite{bk}),  fixed points for non-self 1-set contractive mappings.

Here, we recall that since the Birkhoff-Kellogg theorem was first introduced and proved by Birkhoff and Kellogg \cite{bk} in 1922 in discussing the existence of solutions for the equation $ x = \lambda F(x)$, where $\lambda$ is a real parameter, and $F$ is a general nonlinear non-self mapping defined on an open convex subset $U$ of a topological vector space $E$, now the general form of the Birkhoff-Kellogg  problem is to find the so-called an invariant direction for the nonlinear single-valued or set-valued mappings $F$, i.e., to find $x_0 \in \overline{U}$ (or $x_0 \in \partial \overline{U}$)  and $\lambda > 0$ such that $\lambda x_0 = F(x_0)$, or $\lambda x_0 \in F(x_0)$. But current paper focuses on the study for single-valued mappings for $p$-vector spaces for $0 < 1 \leq 1$.

Since Birkhoff and Kellogg theorem given by Birkhoff and Kellogg in 1920's, the study on Birkhoff-Kellogg problem has been received a lot of attention by scholars since then, for example, one of the fundamental results in nonlinear functional analysis, called the Leray-Schauder alternative by Leray and Schauder \cite{lerayschauder} in 1934, was established via topological degree. Thereafter, certain other types of Leray-Schauder alternatives were proved using different techniques other than topological degree, see work given by Granas and Dugundji \cite{granas}, Furi and Pera \cite{furipera} in the Banach space setting and applications to the boundary value problems for ordinary differential equations, and a general class of mappings for nonlinear alternative of Leray-Schauder type in normal topological spaces, and also Birkhoff-Kellogg type theorems for general class mappings in TVS by Agarwal et al.\cite{Agarwal}, Agarwal and O'Regan \cite{agarwalorgan2003}-\cite{agarwalorgan2004}, Park \cite{park1997}; in particular, recently O'Regan \cite{oregan2019} using the Leray-Schauder type coincidence theory to establish some Birkhoff-Kellogg problem, Furi-Pera type results for a general class of single-valued or set-valued mappings, too.  

In this section, one best approximation result for 1-set contractive mappings in locally $p$-convex spaces is first established, which is then used to establish the solution principle for Birkhoff-Kellogg problems and related nonlinear alternatives; these new results allow us to give general principle for Leray - Schaduer type, and related fixed point theorems of non-self mappings in locally $p$-convex spaces for $(0 < p \leq 1)$. The new results given in this part not only include the corresponding results in the existing literature as special cases, but also would expected to play the fundamental role for the development of nonlinear problems arising from theory to practice for 1-set contractive mappings under the framework of $p$-vector spaces, which include the general topological vector spaces as a special class.

We also note that the general principles for nonlinear alternative related to Leray - Schauder alternative and other types under the framework of locally  $p$-convex spaces for $(0< p \leq 1)$ given in this section would be useful tools for the study of nonlinear problems. In addition, we also note that corresponding results in the existing literature for Birkhoff-Kellogg problems and the Leray - Schauder alternatives have been studied comprehensively by Granas and Dugundji \cite{granas}, Isac \cite{isac}, Kim et al.\cite{kim2002}, Park \cite{park2010}-\cite{park100}, Carbone and Conti \cite{carbone1991}, Chang et al.\cite{chang1999}-\cite{chang2001}, Chang and Yen \cite{chang1996}, Shahzad \cite{shahzad2004}-\cite{shahzad2006}, Singh \cite{singh1997}; and in particular, many general forms recently obtained by O'Regan \cite{oregan2021}, Yuan \cite{yuan2022} and references therein.

In order to study the general existence of fixed points for non-self mappings in locally $p$-convex spaces, we need some definitions and notations given below.
\vskip.1in
\noindent
{\bf Definition 5.1 (Inward and Outward sets in $p$-vector spaces).}
Let $C$ be a subset of a $p$-vector space $E$ and $x \in E$ for $0 < p \leq 1$. Then the $p$-Inward set $I^p_C(x)$ and
$p$-Outward set $O^p_C(x)$ are defined by

$I^p_C(x): =\{ x + r(y-x): y \in C, \mbox { for any } r \geq 0 \mbox { (1) if  } 0 \leq r \leq 1, \mbox { with } (1-r)^p + r^p =1;
\mbox{ or (2) if } r \geq 1, \mbox { with  } (\frac{1}{r})^p + (1- \frac{1}{r})^p = 1 \};$  and

$O^p_C(x): =\{x + r(y-x): y \in C, \mbox{ for any }  r \leq 0 \mbox{ (1) if } 0 \leq |r| \leq 1, \mbox{  with  }
(1-|r|)^p + |r|^p = 1;  \mbox{ or (2) if  } |r| \geq 1,  \mbox { with } (\frac{1}{|r|})^p + (1-\frac{1}{|r|})^p =1 \}.$

From the definition, it is obviously that when $p=1$, the both inward and outward sets $I^p_C(x)$, $O^p_C(x)$ are reduced to the definition for the inward set $I_C(x)$ and the outward set $O_C(x)$, respectively in topological vector spaces introduced by Halpern and Bergman \cite{hergman1968} and used for the study of non-self mappings related to nonlinear functional analysis in the literature. In this paper, we will mainly focus on the study of the $p$-inward set $I_U^p(x)$ for the best approximation and related to the boundary condition for the existence of the fixed points in locally $p$-convex spaces. By the special property of $p$-convex concept for $p \in (0, 1)$ and $p=1$, we have the following fact.

\vskip.1in
\noindent
{\bf Lemma 5.1.}
Let $C$ be a subset of a $p$-vector space $E$ and $x \in E$, where for $0 < p \leq 1$. Then for both $p$-Inward and Outward sets  $I^p_C(x)$ and $O^p_C(x)$ defined above, we have

(I) when $p \in (0, 1)$, $I^p_C(x)= [\{x\}\cup C]$, and $O^p_C(x)=[\{x \} \cup \{2x\} \cup - C ]$,

(II) when $p=1$, in general $ [\{x \}\cup C] \subset I^p_C(x)$, and $[\{ x \} \cup \{2x\} \cup -C] \subset O^p_C(x)$.

\noindent
{\bf Proof.} First, when $ p\in (0, 1)$, by the definitions of $I^p_C(x)$, the only real number $r \geq 0$ satisfying the equation $(1-r)^p + r^p =1$ for $r\in [0,1]$ is $r=0$ or $r=1$, and when $r \geq 1$, the equation
 $(\frac{1}{r})^p + (1- \frac{1}{r})^p = 1$ implies that $r=1$. The same reason for $O^p_C(x)$, it follows that $r=0$ and $r= -1$.

Secondly when $p=1$, all $r\geq 0$, and all $r\leq 0$ satisfy the requirement of definition for $I^p_C(x)$ and $O^p_C(x)$, respectively, thus the proof is compete. $\square$

\vskip.1in
By following the original idea by Tan and Yuan \cite{tanyuan1994} for hemicompact mappings in metric spaces, we introduce the following definition for a mapping being hemicompact in $p$-seminorm spaces for $p \in (0,1]$, which is indeed the $``${\bf (H) condition"} below used in Theorem 5.1 to prove the existence of best approximation results for 1-set contractive mappings in locally $p$-convex spaces for $p \in (0, 1]$.

\vskip.1in
\noindent
{\bf Definition 5.2 (Hemicompact mapping)}.
Let $E$ be a locally $p$-convex space for $1 < p \leq 1$.
For a given bonded (closed) subset $D$ in $E$, a mapping $F: D \rightarrow 2^E$ is said to be hemicompact
if each sequence $\{x_n\}_{n\in N}$ in $D$ has a convergent subsequence with limit $x_0$ such that $x_0 \in F(x_0)$, whenever $\lim_{n \rightarrow \infty} d_{P_U}P(x_n, F(x_n)) =0$ for each $U \in \mathfrak{U}$, where $d_{P_U}P(x, C):= \inf\{P_U(x-y): y \in C\}$ is the distance of a single point $x$ with the subset $C$ in $E$ based on $P_U$, $P_U$ is the Minkowski $p$-functional in $E$ for $U \in \mathfrak{U}$, which is the base of the family consisted by all open $p$-convex subsets for $0$-neighborhoods in $E$.

\vskip.1in
\noindent
{\bf Remark 5.1}. We like to point that the Definition 5.2 is indeed an extension for a $``$hemicompact mapping" defined from a metric space to a (locally) $p$-convex space with the $p$-seminorm, where $p \in (0, 1]$ (see Tan and Yuan \cite{tanyuan1994}).
By the monotonicity of Minkowski $p$-functionals, i.e., the biggger $0$-neighborhoods, the smaller Minkowski $p$-functionals' values (see also p.178 of Balachandran \cite{balachandra}), the Definition 5.2 describes the converge for the distance between $x_n$ and $F(x_n)$ by using the language of seminorms in terms of Minkowski $p$-functionals for each $0$-neighborhood in $\mathfrak{U}$ (the base), which is the family consisted of its open $p$-convex $0$-neighborhoods in $p$-vector space $E$.

\vskip.1in
Now we have the following Schauder fixed point theorem for 1-set contractive mappings in locally $p$-convex spaces for $p \in(0, 1]$.

\vskip.1in
\noindent
{\bf Theorem 5.1 (Schauder fixed point theorem for single-valued 1-set contractive mappings).}
Let $U$ be  a non-empty bounded open subset of a (Hausdorf) locally $p$-convex space $E$ and its zero $0 \in U$, and $C \subset E$  be a closed $p$-convex subset of $E$ such that $0 \in C$, with $0 < p \leq 1$.
 If $F: C \cap \overline{U} \rightarrow C \cap \overline{U}$ is continuous and $1$-set contractive single-valued mapping and satisfying the following (H) or (H1) condition:

{\bf (H) Condition}: The sequence $\{x_n\}_{n\in \mathbb{N}}$ in $\overline{U}$ has a convergent subsequence with limit $x_0 \in \overline{U}$  such that $x_0 \in F(x_0)$,  whenever  $\lim_{n \rightarrow \infty} d_{P_U}(x_n, F(x_n)) =0$, where, $d_{P_U}(x_n, F(x_n)):=P_U(x_n -  F(x_n)\}$, where $P_U$ is the Minkowski $p$-functional for any $U \in \mathfrak{U}$, which is the family of all non-empty open $p$-convex subsets of zero in $E$.

{\bf (H1) Condition}: There exists $x_0$ in $\overline{U}$ with $x_0 = F(x_0)$ if there exists $\{x_n\}_{n\in \mathbb{N}}$ in $\overline{U}$ such that $\lim_{n \rightarrow \infty} d_{P_U}(x_n, F(x_n)) =0$, where, $P_U$ is the Minkowski $p$-functional for any $U \in \mathfrak{U}$, which is the family of all non-empty open $p$-convex subsets of zero in $E$.

Then $F$ has at least one fixed point in $C \cap \overline{U}$.

\noindent
{\bf Proof.} Let $U$ be any element in $\mathfrak{U}$, which is the family of all non-empty open $p$-convex subset for zero in $E$.
As the mapping $T$ is 1-set contractive, taking an increasing sequence $\{\lambda_n\}$ such that $0 < \lambda_n < 1$ and  $\lim_{n \rightarrow \infty} \lambda_n =1$, where $n \in \mathbb{N}$.
Now we define a mapping $F_n: C \rightarrow C$ by
$F_n(x): = \lambda_n F(x)$ for each $x \in  C$ and $n\in \mathbb{N}$. Then it follows that $F_n$ is a $\lambda_n$-set-contractive mapping with $ 0 < \lambda_n < 1$ . By Theorem 4.5 on the condensing mapping $F_n$ in $p$-vector space with $p$-seminorm $P_U$ for each $n \in \mathbb{N}$, there exists $x_n \in C $ such that $x_n \in F_n(x_n)=\lambda_n F(x_n)$.
As $P_U$ is the Minkowski $p$-functional of $U$ in $E$, it follows that $P_U$ is continuous as $0 \in int(U)=U$.
Note that for each $n \in \mathbb{N}$, $\lambda_n x_n \in \overline{U} \cap C$, which imply that
$x_n = r(\lambda_n F(x_n)) = \lambda_n F(x_n)$, thus $P_U(\lambda_n F(x_n)) \leq 1$ by Lemma 2.2.
Note that
$$P_U(F(x_n) - x_n)=P_U(F(x_n) - \lambda_n F(x_n)) =P_U(\frac{(1-\lambda_n) \lambda_n F(x_n)}{\lambda_n})
\leq (\frac{1-\lambda_n}{\lambda_n})^p P_U(\lambda_n F(x_n)) \leq (\frac{1-\lambda_n}{\lambda_n})^p,$$
which implies that
$\lim_{n\rightarrow \infty} P_U(F(x_n) - x_n)=0$ for all $U \in \mathfrak{U}$.

Now (1) if $F$ satisfies the (H) condition,  it implies that the consequence
$\{x_n\}_{n \in \mathbb{N}}$ has a convergent subsequence which converges to $x_0$ such that $x_0 = F(x_0)$.
Without loss of the generality, we assume that $\lim_{n \rightarrow \infty} x_n=x_0$, is with $x_n=\lambda_n F(x_n)$, and $\lim_{n \rightarrow \infty} \lambda_n=1$, it implies that
$x_0=\lim_{n \rightarrow \infty} (\lambda_n F(x_n))$, which means  $\lim_{n\rightarrow \infty} F(x_n)= x_0$.

(ii) if $F$ satisfies the (H1) condition, then by the (H1) condition, it follows that there exists $x_0$ in $\overline{U}$ such that
$x_0 = F(x_0)$, which is a fixed point of $F$. We complete the proof. $\square$

\vskip.1in
\noindent
{\bf Theorem 5.2 (Best approximation for single-valued 1-set-contractive mappings).}
Let $U$ be a bounded open $p$-convex subset of a locally $p$-convex space $E$ ($0 \leq p \leq 1)$  the zero $0 \in U$, and $C$ a (bounded) closed $p$-convex subset of $E$ with also zero $0\in C$. Assume $F: \overline{U}\cap C \rightarrow C$ is a
is 1-set contractive, and for each $x \in \partial_C U$ with $F(x) \in C \diagdown \overline{U}$, $(P^{\frac{1}{p}}_U(F(x))- 1)^p \leq P_U (F(x)-x)$ for $0< p \leq 1$ (this is trivial when $p=1$). In addition, if $F$ satisfies the following (H) or (H1) condition:

{\bf (H) Condition}: The sequence $\{x_n\}_{n\in \mathbb{N}}$ in $\overline{U}$ has a convergent subsequence with limit $x_0 \in \overline{U}$   such that $x_0 = F(x_0)$,  whenever  $\lim_{n \rightarrow \infty} d_{P_U}(x_n, F(x_n)) =0$, where, $d_{P_U}(x_n, F(x_n)):=\inf\{P_U(x_n- F(x_n)\}$, where $P_U$ is the Minkowski $p$-functional for any $U \in \mathfrak{U}$, which is the family of all non-empty open $p$-convex subset containing the zero in $E$.

{\bf (H1) Condition}: There exists $x_0$ in $\overline{U}$ with $x_0 = F(x_0)$ if there exists $\{x_n\}_{n\in \mathbb{N}}$ in $\overline{U}$ such that $\lim_{n \rightarrow \infty} d_{P_U}(x_n, F(x_n)) =0$, where, $P_U$ is the Minkowski $p$-functional for any $U \in \mathfrak{U}$, which is the family of all non-empty open $p$-convex subset containing the zero in $E$.

Then we have that there exist $x_0 \in C \cap \overline{U}$ such that
$$ P_U (F(x_0) - x_0) = d_P(y_0, \overline{U}\cap C) =  d_p(F(x_0), \overline{I^p_{\overline{U}}( x_0)} \cap C),$$
where $P_U$ is the Minkowski $p$-functional of $U$. More precisely, we have the following either (I) or (II) holding:

(I)  $F$ has a fixed point $x_0 \in \overline{U} \cap C$, i.e., $0=P_U (F(x_0) - x_0) = d_P(F(x_0), \overline{U}\cap C) =  d_p(F(x_0), \overline{I^p_{\overline{U}}( x_0)} \cap C)$,

(II)  there exists $x_0 \in \partial_C(U)$ and $F(x_0) \notin \overline{U}$ with
$$ P_U (F(x_0) - x_0) = d_P(F(x_0), \overline{U}\cap C) = d_p(F(x_0), \overline{I^p_{\overline{U}}( x_0)} \cap C)=(P^{\frac{1}{p}}_U(F(x_0))-1)^p  > 0.$$

\noindent
{\bf Proof.} As $E$ is a locally $p$-convex space $E$, it suffices to prove that for each open $p$-convex subset $U$ in
$\mathfrak{U}$ (which is the family of all non-empty open $p$-convex subset containing the zero in $E$),
there exists a sequence $(x_n)_{n \in \mathbb{N}}$ in $\overline{U}$ such that
$\lim_{n\rightarrow \infty} P_U(F(x_n)-x_n)=0$, and the conclusion follows by applying the (H) condition.

Let $r: E \rightarrow U$ be a retraction mapping defined by $r(x): = \frac{x}{\max\{1, (P_U(x))^{\frac{1}{p}}\}}$
for each $x \in E$, where $P_U$ is the Minkowski $p$-functional of $U$.
Since the space $E$'s zero $0 \in U(=intU$ as $U$ is open), it follows that $r$ is continuous by Lemma 2.2.
As the mapping $F$ is 1-set contractive, taking an increasing sequence $\{\lambda_n\}$ such that $0 < \lambda_n < 1$ and  $\lim_{n \rightarrow \infty} \lambda_n =1$, where $n \in \mathbb{N}$.
Now for each $n\in \mathbb{N}$, we define a mapping $F_n: C \cap \overline{U} \rightarrow C$ by
$F_n(x): = \lambda_n F\circ r(x)$ for each $x \in C \cap \overline{U}$.
By the fact that $C$ and $\overline{U}$ are $p$-convex, it follows that $r(C) \subset C$ and $r(\overline{U}) \subset \overline{U}$, thus $r( C \cap \overline{U}) \subset C \cap \overline{U}$. Therefore $F_n$ is a mapping from $\overline{U}\cap C$ to itself. Nor each $n \in \mathbb{N}$, by the fact that $F_n$ is a $\lambda_n$-set-contractive mapping with $ 0 < \lambda_n < 1$  it follows by Theorem 4.5 for the condensing mapping that there exists $z_n \in C \cap \overline{U}$ such that $F_n(z_n)=\lambda_n F \circ r(z_n)$.
As $r( C \cap \overline{U}) \subset C \cap \overline{U}$, let $x_n= r(z_n)$.  Then we have that $x_n \in  C\cap \overline{U}$ and
with $x_n = r(\lambda_n F_n(x_n))$ such that the following (1) or (2) holding for each $n \in \mathbb{N}$:  (1) $\lambda_n F_n(x_n)\in C\cap \overline{U}$;  or (2) $\lambda_n F_n(x_n) \in C \diagdown \overline{U}$.

Now we prove the conclusion by considering the following two cases under the (H) condition and (H1) condition.

Case (I) For each $n \in N$, $\lambda_n F(x_n) \in C \cap \overline{U}$; or

Case (II) There there exists a positive integer $n$ such that $\lambda_n F(x_n) \in C \diagdown \overline{U}$.

First, by the case (I), for each $n \in \mathbb{N}$, $\lambda_n F(x_n) \in \overline{U} \cap C$, which imply that
$x_n = r(\lambda_n F(x_n)) = \lambda_n F(x_n)$, thus $P_U(\lambda_n F(x_n)) \leq 1$ by Lemma 2.2.
Note that
$$P_U(F(x_n)- x_n)=P_U(F(x_n)- x_n)=P_U(F(x_n)- \lambda_n F(x_n)) =P_U(\frac{(1-\lambda_n) \lambda_n F(x_n)}{\lambda_n})$$
$$ \leq (\frac{1-\lambda_n}{\lambda_n})^p P_U(\lambda_n F(x_n)) \leq (\frac{1-\lambda_n}{\lambda_n})^p,$$
which implies that
$\lim_{n\rightarrow \infty} P_U(F(x_n)-x_n)=0$. Now for any $V \in \mathbb{U}$, without loss of generality, let $U_0 = V \cap U$. Then we have the following conclusion:
$$P_{U_0}(F(x_n)- x_n)=P_{U_0}(F(x_n)- x_n)=P_{U_0}(F(x_n)- \lambda_n F(x_n)) =P_{U_0}(\frac{(1-\lambda_n) \lambda_n F(x_n)}{\lambda_n})$$
$$\leq (\frac{1-\lambda_n}{\lambda_n})^p P_{U_0}(\lambda_n F(x_n)) \leq (\frac{1-\lambda_n}{\lambda_n})^p,$$
which implies that
$\lim_{n\rightarrow \infty} P_{U_0}(F(x_n)-x_n)=0$, where $P_{U_0}$ is the Minkowski $p$-functional of $U_0$ in $E$.

Now if $F$ satisfies the (H) condition, if follows that the consequence
$\{x_n\}_{n \in \mathbb{N}}$ has a convergent subsequence which converges to $x_0$ such that $x_0 = F(x_0)$.
Without loss of the generality, we assume that $\lim_{n \rightarrow \infty} x_n=x_0$, $x_n=\lambda_n y_n$, and $\lim_{n \rightarrow \infty} \lambda_n= 1$, and as $x_0=\lim_{n \rightarrow \infty} (\lambda_n F(x_n))$, which implies that $F(x_0)=\lim_{n\rightarrow \infty} F(x_n)= x_0$. Thus there exists $x_0 =F(x_0)$, thus we have
$0 = d_p(x_0, F(x_0)) = d(y_0, \overline{U}\cap C) = d_p(F(x_0), \overline{I^p_{\overline{U}}(x_0)} \cap C))$ as indeed $x_0 = F(x_0) \in \overline{U}\cap C \subset \overline{I^p_{\overline{U}}( x_0)} \cap C)$.

If $F$ satisfies the (H1) condition, if follows that there exists $x_0 \in \overline{U} \cap C$ with $x_0 = F(x_0)$.
Then we have $0=P_U (F(x_0) - x_0) = d_P(F(x_0), \overline{U}\cap C) =  d_p(F(x_0), \overline{I^p_{\overline{U}}( x_0)} \cap C)$.

Second, by the case (II) there exists a positive integer $n$ such that $\lambda_n F(x_n) \in C \diagdown \overline{U}$.
Then we have that $P_U(\lambda_n F(x_n))> 1$, and also $P_U(F(x_n))> 1$ as $\lambda_n < 1$.
As $ x_n = r(\lambda_n F(x_n)) = \frac{\lambda_n F(x_n)}{(P_U(\lambda_n F(x_n)))^{\frac{1}{p}}}$, which implies that $P_U(x_n)=1$, thus $x_n \in \partial_C(U)$.
Note that
$$P_U(F(x_n) - x_n)=P_U(\frac{(P_U(F(x_n))^{\frac{1}{p}}-1)F(x_n)}{P_U(F(x_n))^{\frac{1}{p}}})=(P^{\frac{1}{p}}_U(F(x_n))-1)^p. $$
By the assumption, we have $(P^{\frac{1}{p}}_U(F(x_n))-1)^p \leq P_U(F(x_n) -x)$ for $x \in C \cap \partial \overline{U}$,
it follows that
$$P_U(F(x_n))-1 \leq P_U(F(x_n)) - \sup\{P_U(z): z \in C\cap \overline{U}\} \leq \inf\{P_U(F(x_n)- z): z \in C \cap \overline{U}\}= d_p(F(x_n), C \cap \overline{U}).$$
Thus  we have the best approximation: $P_U(F(x_n) - x_n)=d_P(y_n, \overline{U} \cap C) = (P^{\frac{1}{p}}_U(F(x_n)-1)^p  > 0.$

Now we want to show that $P_U(y_n-x_n)= d_P(F(x_n), \overline{U} \cap C) = d_p(F(x_n), \overline{I^p_{\overline{U}}( x_0)} \cap C) > 0.$

By the fact that
$(\overline{U}\cap C) \subset I^p_{\overline{U}}(x_n)\cap C$, let $z \in I^p_{\overline{U}}(x_n)\cap C \diagdown(\overline{U}\cap C)$, we first claim that  $P_U(F(x_n) - x_n) \leq P_U(F(x_n) - z)$. If not, we have $P_U(F(x_n) - x_n) > P_U(F(x_n)-z)$.
As $z \in I^p_{\overline{U}}(x_n) \cap C \diagdown (\overline{U} \cap C)$, there exists $y \in \overline{U}$ and a non-negative number $c$ (actually $c\geq 1$ as shown soon below) with $z = x_n + c (y - x_n)$. Since $z \in C$, but $z \notin \overline{U} \cap C$, it implies that $z \notin \overline{U}$. By the fact that $x_n\in \overline{U}$ and $y \in \overline{U}$, we must have
the constant $c \geq 1$; otherwise, it implies that $z ( = (1- c )x_n + c y) \in \overline{U}$, this is impossible by our assumption, i.e., $z\notin \overline{U}$. Thus we have that $c\geq 1$, which implies that $y =\frac{1}{c} z + (1-\frac{1}{c}) x_n \in C$ (as both $x_n \in C$ and $z\in C$). On the other hand, as $z \in I^p_{\overline{U}}(x_n) \cap C \diagdown (\overline{U} \cap C)$, and $c\geq 1$ with
 $(\frac{1}{c})^p+ (1-\frac{1}{c})^p = 1 $, combing with our assumption that for each $x \in \partial_C \overline{U}$ and
 $y \in F(x_n)\diagdown \overline{U}$,
 $P^{\frac{1}{p}}_U(y)- 1 \leq P^{\frac{1}{p}}_U (y-x)$ for $0< p \leq 1$, it then follows that
$$P_U(F(x_n)- y) = P_U[\frac{1}{c}(F(x_n)- z)+(1-\frac{1}{c})(F(x_n) - x_n)] \leq
[(\frac{1}{c})^{p} P_U(F(x_n) -z)+(1-\frac{1}{c})^p P_U(F(x_n) - x_n)] < P_U(F(x_n)- x_n),$$
which contradicts that $P_U (F(x_n) - x_n) = d_P(F(x_n), \overline{U}\cap C)$ as shown above we know that $y \in \overline{U}\cap C$, we should have $P_U(F(x_n)- x_n)\leq P_U(F(x_n) - y)$! This helps us to complete the claim:
$P_U(F(x_n) - x_n) \leq P_U(F(x_n) - z)$ for any $z \in I^p_{\overline{U}}(x_n)\cap C \diagdown(\overline{U}\cap C)$, which means that the following best approximation of Fan's type (see \cite{fan1969}-\cite{fan1972}) holding:
$$ 0 < d_P(F(x_n), \overline{U}\cap C) = P_U (F(x_n) - x_n) = d_p(F(x_n), I^p_{\overline{U}}(x_n) \cap C).$$
Now by the continuity of $P_U$, it follows that the following best approximation of Fan type is also true:
$$ 0 <  P_U(F(x_n) - x_n) = d_P(F(x_n), \overline{U}\cap C) = d_p(F(x_n), I^p_{\overline{U}}(x_n) \cap C)
= d_p(F(x_n), \overline{I^p_{\overline{U}}(x_n)} \cap C);$$
and we have the conclusion below due to that $\lim_{n \rightarrow \infty}x_n=x_0$ and the continuity of $F$ (actually $x_0 \neq F(x_0)$):
$$  P_U(F(x_0) - x_0) = d_P(F(x_0), \overline{U}\cap C) = d_p(F(x_0), I^p_{\overline{U}}(x_0) \cap C)
= d_p(F(x_0), \overline{I^p_{\overline{U}}(x_0)} \cap C)=(P^{\frac{1}{p}}_U(F(x_0))-1)^p  > 0.$$
This completes the proof. $\square$

\vskip.1in
\noindent
{\bf Remark 5.2.} We note that Theorem 5.2 also improves the corresponding best approximation for 1-set contractive mappings given by
Li et al.\cite{lixuduan2006}, Liu \cite{liu2001}, Xu \cite{xu2007}, Xu et al.\cite{xujiali2006}, and results from the references therein; and 3): When $p$=1, we have the similar best approximation result for the mapping $F$ in the locally convex spaces with outward set boundary condition below (see Theorem 3 of Park \cite{park1995b} and related discussion by references therein).

\vskip.1in
Though the main focus of this paper studies best approximation, fixed point theorems for single-valued mappings, when a $p$-vector space $E$ (for $p=1$) being a locally convex space (LCS), we can also have the following best approximation for upper semicontinous set-valued mappings by applying Theorem 4.6 with arguments used by Theorem 5.1 and Theorem 5.2 above (see also more discussion given by Yuan \cite{yuan2022} and references wherein).

\vskip.1in
\noindent
{\bf Theorem 5.3 (Best approximation for USC set-valued mappings in LCS).}
Let $U$ be a bounded open convex subset of a locally convex space $E$ (i.e., $p=1$)  with zero $0 \in intU=U$ (the interior $intU=U$ as $U$ is open), and $C$ a closed $p$-convex subset of $E$ with also zero $0\in C$. Assume that $F: \overline{U}\cap C \rightarrow 2^C$ is a 1-set-contractive upper semicontinuous mapping, and satisfying the condition (H) or (H1) above.
Then there exist $x_0 \in \overline{U} \cap X$ and $y_0 \in F(x_0)$ such that
$ P_U (y_0 - x_0) = d_P(y_0, \overline{U}\cap C) =  d_p(y_0, \overline{I_{\overline{U}}( x_0)} \cap C),$
where $P_U$ is the Minkowski $p$-functional of $U$. More precisely, we have the following either (I) or (II) holding:

(I)  $F$ has a fixed point $x_0 \in U \cap C$, i.e., $x_0 \in F(x_0)$ (so that $P_U (y_0- x_0) = P_U (y_0 - x_0) = d_P(y_0, \overline{U}\cap C)
=  d_p(y_0, \overline{I_{\overline{U}}( x_0)} \cap C))=0$), or

(II)  there exist $x_0 \in \partial_C(U)$ and $y_0 \in F(x_0)$ with $y_0 \notin \overline{U}$ with
$$ P_U (y_0 - x_0) = d_P(y_0, \overline{U}\cap C) = d_p(y_0, I_{\overline{U}}(x_0)\cap C)
=d_p(y_0, \overline{I_{\overline{U}}(x_0)} \cap C) > 0.$$

\noindent
{\bf Proof.} By following the proof used in Theorem 5.1 and Theorem 5.2, then applying Theorem 4.6 for $p=1$, the conclusion follows.
This completes the proof. $\square$

Now by the application of Theorem 5.2 with Remark 5.2 and the argument used in Theorem 5.2, we have the
the following general principle for the existence of solutions for Birkhoff-Kellogg problems in $p$-seminorm spaces for locally $p$-convex spaces, where $0 < p \leq 1$.

\vskip.1in
\noindent
{\bf Theorem 5.4 (Principle of Birkhoff-Kellogg alternative).}
Let $U$ be a bounded open $p$-convex subset of a locally $p$-convex space $E$ ($0 \leq p \leq 1)$  with zero $0 \in intU=(U)$ (the interior $intU$ as $U$ is open), and $C$ a closed $p$-convex subset of $E$ with also zero $0\in C$. Assume that $F: \overline{U}\cap C \rightarrow C$ is a 1-set-contractive continuous mapping, and satisfying the (H) or (H1) condition above.
Then $F$ has at least one of the following two properties:

(I) $F$ has a fixed point $x_0 \in U \cap C$ such that $x_0 = F(x_0)$,

(II) there exist $x_0 \in \partial_C(U)$,  $F(x_0) \notin \overline{U}$,  and $\lambda = \frac{1}{(P_U(F(x_0))^{\frac{1}{p}}} \in (0, 1)$ such that
$x_0 = \lambda F(x_0)$; In addition if for each $x \in \partial_C U$, $P^{\frac{1}{p}}_U(F(x))- 1 \leq P^{\frac{1}{p}}_U (F(x)-x)$ for $0< p \leq 1$ (this is trivial when $p=1$), then the best approximation between $\{x_0\}$ and $F(x_0)$ given by $$ P_U (F(x_0) - x_0) = d_P(F(x_0), \overline{U}\cap C) = d_p(F(x_0), \overline{I^p_{\overline{U}}(x_0)} \cap C) = (P^{\frac{1}{p}}_U(F(x_0))-1)^p > 0.$$

\vskip.1in
{\bf Proof.} If (I) is  not the case, then (II) is proved by the Remark 5.2 and by following the proof in Theorem 5.2 for the case
ii): $F(x_0)\notin \overline{U}$ with $F(x_0)=f(x_0)$, where $f$ is the restriction of the continuous retraction $r$ respect to the set $U$ in $E$ defined in the proof of Theorem 5.2 above.
Indeed, as $F(x_0) \notin \overline{U}$, it follows that $P_U(F(x_0)) > 1$, and $x_0=f(F(x_0)) =F(x_0)\frac{1}{(P_U(F(x_0))^{\frac{1}{p}}}$.
Now let $\lambda = \frac{1}{(P_U(F(x_0))^{\frac{1}{p}}}$, we have $\lambda < 1$ and  $x_0 = \lambda F(x_0)$.
Finally, the additionally assumption in (II) allows us to have the  best approximation between $x_0$ and $F(x_0)$ obtained by following the proof of Theorem 5.2 as $P_U (F(x_0) - x_0) = d_P(F(x_0), \overline{U}\cap C) = d_p(F(x_0), \overline{I^p_{\overline{U}}( x_0)} \cap C) > 0$. This completes the proof. $\square$

\vskip.1in
As an application of Theorem 5.3 for the non-self upper semicontinuous set-valued mappings discussed in Theorem 5.4, we have the
following general principle of Birkhoff-Kellogg alternative in locally convex spaces.

\vskip.1in
\noindent
{\bf Theorem 5.5 (Principle of Birkhoff-Kellogg alternative in LCS).}
Let $U$ be a bounded open $p$-convex subset of a LCS $E$  with the zero $0 \in U$, and $C$ a closed convex subset of $E$ with also zero $0\in C$. Assume the
$F: \overline{U}\cap C \rightarrow 2^C$ is a 1-set contractive and upper semicontinuous mapping, and satisfying the (H) or (H1) condition (H) above. Then it has at least one of the following two properties:

(I) $F$ has a fixed point $x_0 \in U \cap C$ such that $x_0 \in F(x_0)$,

(II) there exists $x_0 \in \partial_C(U)$ and $y_0 \in F(x_0)$ with $y_0 \notin \overline{U}$  and $\lambda \in (0, 1)$ such that $x_0 = \lambda y_0$, and the best approximation between $\{x_0\}$ and $F(x_0)$ is given by
$ P_U (y_0 - x_0) = d_P(y_0, \overline{U}\cap C) = d_p(y_0, \overline{I^p_{\overline{U}}( x_0)} \cap C) > 0.$

\vskip.1in
On the other hand, by the Proof of Theorems 5.2, we note that for case (II) of Theorem 5.2, the assumption $``$each $x \in \partial_C U$ with $P^{\frac{1}{p}}_U(F(x)- 1 \leq P^{\frac{1}{p}}_U (F(x)-x)$" is only used to guarantee the best approximation $``P_U (F(x_0) - x_0) =
d_P(F(x_0), \overline{U}\cap C) = d_p(F(x_0), \overline{I^p_{\overline{U}}( x_0)} \cap C) > 0$", thus
we have the following Leray-Schauder alternative in $p$-vector spaces, which, of course, includes the corresponding results in locally convex spaces as special cases.

\vskip.1in
\noindent
{\bf Theorem 5.6 (The Leray-Schauder Nonlinear Alternative).} Let $C$ a closed $p$-convex subset of $p$-seminorm space $E$ with $0 \leq p \leq 1$  and the zero $0 \in C$. Assume the $F: C \rightarrow C$ is a 1-set contractive and continuous mapping, and satisfying the (H) or (H1) condition above.
Let $\varepsilon(F): =\{x \in C: x = \lambda F(x), \mbox{ for some } 0 < \lambda < 1\}$. Then either $F$ has a fixed point in $C$ or the set $\varepsilon(F)$ is unbounded.

\noindent
{\bf Proof.} We prove the conclusion by assuming that $F$ has no fixed point, then we claim that the set $\varepsilon(F)$ is unbounded. Otherwise, assume the set $\varepsilon(F)$ is bounded. and assume $P$ is the continuous $p$-seminorm for $E$, then there exists $r>0$ such that the set $B(0, r):=\{x \in E: P(x) < r\}$ , which contains the set $\varepsilon(F)$, i.e., $\varepsilon(F) \subset B(0, r)$, which means for any $x \in \varepsilon(F)$, $P(x) < r$. Then $B(0. r)$ is an open $p$-convex subset of $E$ and the zero $0 \in B(0, r)$ by Lemma 2.2 and Remark 2.4.
Now let $U:=B(0, r)$ in Theorem 5.4, it follows that for the mapping $F: B(0, r) \cap C \rightarrow C$ satisfies all general conditions of Theorem 5.4, and we have that any $x_0 \in \partial_C B(0, r)$, no any $\lambda \in (0, 1)$ such that
$x_0=\lambda F(x_0)$. Indeed, for any $x \in \varepsilon(F)$, it follows that $P(x) < r$ as $\varepsilon(F) \subset B(0, r)$, but for any $x_0 \in \partial_C B(0, r)$, we have $P(x_0)=r$, thus the conclusion (II) of Theorem 5.4 does not have hold. By Theorem 5.4 again, $F$ must have a fixed point, but this contradicts with our assumption that $F$ is fixed point free. This completes the proof. $\square$

\vskip.1in
Now assume a given $p$-vector space $E$ equipped with the $P$-seminorm (by assuming it is continuous at zero) for $0< p \leq 1$, then we know that
$P: E \rightarrow \mathbb{R}^+$, $P^{-1}(0)=0$,  $P(\lambda x) = |\lambda|^p P(x)$ for any  $x\in E$ and $\lambda \in \mathbb{R}$. Then we have the following useful result for fixed points due to Rothe and Altman types in $p$-vector spaces, in particular, for locally $p$-convex spaces,  which plays important roles for optimization problem, variational inequality, complementarity problems (see isac \cite{isac}, or Yuan \cite{yuan1999} and references therein for related study in details).

\vskip.1in
\noindent
{\bf Corollary 5.1.} Let $U$ be a bounded open $p$-convex subset of a locally $p$-convex space $E$ and zero $0 \in U$, plus $C$ is a closed $p$-convex subset of $E$ with $U \subset C$, where $0< p \leq 1$. Assume that $F: \overline{U} \rightarrow C$ is a
1-set contractive continuous mapping, and satisfying the (H) or (H1) condition above. If one of the following is satisfied,

(1) (Rothe type condition): $P_U(F(x)) \leq P_U(x)$ for $x \in \partial U$;

(2) (Petryshyn type condition): $P_U(F(x)) \leq P_U(F(x)-x)$ for $x \in \partial U$;

(3) (Altman type condition): $|P_U(F(x))|^{\frac{2}{p}} \leq [P_U(F(x))- x)]^{\frac{2}{p}} + [P_U(x)]^{\frac{2}{p}}$ for $x \in \partial U$;

\noindent
then $F$ has at least one fixed point.

\noindent
{\bf Proof.} By the conditions (1), (2) and (3), it follows that the conclusion of (II) in Theorem 5.4 $``$there exist $x_0 \in \partial_C(U)$ and $\lambda \in (0, 1)$ such that $x_0 \neq  \lambda F(x_0)$" does not hold, thus by the alternative of Theorem 5.4, $F$ has a fixed point. This completes the proof. $\square$.

\vskip.1in
By the fact that for $p=1$, when a $p$-vector space is a locally convex space,  we have the following classical Fan's best approximation (see \cite{fan1969}), which is a powerful tool for nonlinear functional analysis in the supporting on the study in the optimization, mathematical programming, games theory, and mathematical economics, and others related topics in applied mathematics.

\vskip.1in
\noindent
{\bf Corollary 5.2 (Fan's best approximation in LCS).} Let $U$ be a bounded open convex subset of a locally convex space $E$ with
the zero $0 \in U$, and $C$ a closed convex subset of $E$ with also zero $0\in C$, and assume $F: \overline{U}\cap C \rightarrow C$ is a 1-set contractive and continuous mapping, and satisfying the (H) or (H1) condition above. Assume $P_U$ being the Minkowski $p$-functional of $U$ in $E$.
Then there exist $x_0 \in \overline{U} \cap X$  such that
$ P_U(F(x_0) - x_0) = d_P(F(x_0), \overline{U}\cap C) =  d_p(F(x_0), \overline{I_{\overline{U}}(x_0)} \cap C).$
More precisely, we have the following either (I) or (II) holding:

(I)  $F$ has a fixed point $x_0 \in U \cap C$, i.e., $x_0=F(x_0)$ (so that $0= P_U (F(x_0) - x_0) = d_P(F(x_0), \overline{U}\cap C) =
 d_p(F(x_0), \overline{I_{\overline{U}}( x_0)} \cap C))$;

(II)  there exists $x_0 \in \partial_C(U)$ and $F(x_0) \notin \overline{U}$ with
$$ P_U (F(x_0) - x_0) = d_P(F(x_0), \overline{U}\cap C) = d_p(F(x_0), \overline{I_{\overline{U}}( x_0)} \cap C) = P_U(F(x_0)) - 1 > 0.$$

\noindent
{\bf Proof.} When $p=1$, then it automatically satisfies that the inequality:
$ P^{\frac{1}{p}}_U((x))- 1 \leq P^{\frac{1}{p}}_U (F(x)-x)$. Now if $F$ has no fixed points, by Theorem 5.4, indeed we have that for $x_0\in \partial_C(U)$, $P_U (F(x_0) - x_0) = d_P(F(x_0), \overline{U}\cap C) = d_p(F(x_0), \overline{I_{\overline{U}}( x_0)} \cap C)=  P_U(F(x_0)-1$. The conclusions are given by Theorem 5.2 (or Theorem 5.3). The proof is complete. $\square$

\vskip.1in
We like to point out the similar results on Rothe and Leray-Schauder alternative have been developed by Isac \cite{isac}, Park \cite{park1995}, Potter \cite{potter1972}, Shahzad \cite{shahzad2006}-\cite{shahzad2004}, Xiao and Zhu \cite{xiaozhu2011}, Yuan \cite{yuan2022}, and related references therein as tools of nonlinear analysis in $p$-vector spaces.

\vskip.1in
\section{Nonlinear Alternatives Principle for the Class of 1-Set Class Contractive Mappings}

\vskip.1in
As applications of results in Section 5 above, we new establish general results for the existence of solutions for Birkhoff-Kellogg problem, and the principle of Leray-Schauder alternatives in locally $p$-convex spaces for $0 < p \leq 1$.

\noindent
{\bf Theorem 6.1 (Birkhoff-Kellogg alternative in locally $p$-convex spaces).} Let $U$ be a bounded open  $p$-convex subset of a locally $p$-convex space $E$ (where, $0 \leq p \leq 1)$ with the zero $0 \in U$, and $C$ a closed $p$-convex subset of $E$ with also zero $0\in C$, and assume $F: \overline{U}\cap C \rightarrow C$ is a 1-set contractive and continuous mapping, and satisfying the condition (H) or (H1) above. In addition, for each $x \in \partial_C(U)$, $P^{\frac{1}{p}}_U(F(x)- 1 \leq P^{\frac{1}{p}}_U (F(x)-x)$ for $0< p \leq 1$ (this is trivial when $p=1$), where $P_U$ is the Minkowski $p$-functional of $U$.
Then we have that either (I) or (II) holding below:

(I) there exists $x_0 \in  \overline{U}\cap C$ such that; or

(II) there exists $x_0 \in \partial_C(U)$ with $F(x_0)\notin \overline{U}$ and $\lambda >1$ such that
$\lambda x_0 = F(x_0)$, i.e., $F(x_0) \in \{\lambda x_0: \lambda > 1 \} \neq \emptyset$.

\noindent
{\bf Proof.} By following the argument and symbols used in the proof of Theorem 5.2, we have that either

(1) $F$ has a fixed point $x_0 \in U \cap C$; or

(2) there exists $x_0 \in \partial_C(U)$ and $x_0=f(F(x_0))$ such that
$$ P_U (F(x_0) - x_0) = d_P(F(x_0), \overline{U}\cap C) = d_p(F(x_0), \overline{I_{\overline{U}}(x_0)} \cap C) = P_U(F(x_0) - 1 > 0,$$
where $\partial_C(U)$ denotes the boundary of $U$ relative to $C$ in $E$, and $f$ is the restriction of the continuous retraction $r$ respect to the set $U$ in $E$ defined in the proof of Theorem 5.2 above.

If $F$ has no fixed point, then above (2) holds and $x_0 \neq F(x_0)$.
As given by the proof of Theorem 5.2, we have that $F(x_0)\notin \overline{U}$, thus $P_U(F(x_0)) > 1$ and $x_0= f(F(x_0))=\frac{F(x_0)}{(P_U(F(x_0))^{\frac{1}{p}}}$, which means $F(x_0) =(P_U(F(x_0))^{\frac{1}{p}} x_0$.
Let $\lambda = (P_U(F(x_0)))^{\frac{1}{p}}$, then $\lambda > 1$ and we have $ \lambda x_0 = F(x_0)$.
This completes the proof. $\square$.

\vskip.1in
\noindent
{\bf Theorem 6.2 (Birkhoff-Kellogg alternative in LCS).} Let $U$ be a bounded open convex subset of a locally convex space $E$ with
the zero $0 \in U$, and $C$ a closed convex subset of $E$ with also zero $0\in C$, and assume $F: \overline{U}\cap C \rightarrow C$ is a 1-set contractive and continuous mapping, and satisfying the condition (H) or (H1) above. Then we have the following either (I) or (II) holding:

(I) there exists $x_0 \in \overline{U}\cap C$ such that $x_0= F(x_0)$; or

(II) there exists $x_0 \in \partial_C(U)$ with $F(x_0)\notin \overline{U}$ and $\lambda >1$ such that
$\lambda x_0 = F(x_0)$, i.e., $F(x_0) \in  \{\lambda x_0: \lambda > 1 \} \neq \emptyset$.

\noindent
{\bf Proof.} When $p=1$, then it automatically satisfies that the inequality:
$ P^{\frac{1}{p}}_U(F(x))- 1 \leq P^{\frac{1}{p}}_U (F(x_0)-x)$, and indeed we have that for $x_0\in \partial_C(U)$,
we have $P_U (F(x_0) - x_0) = d_P(F(x_0), \overline{U}\cap C) = d_p(F(x_0), \overline{W_{\overline{U}}( x_0)} \cap C)=  P_U(F(x_0))-1$. The conclusions are given by 5.4. The proof is complete. $\square$

\vskip.1in
Indeed, we have the following fixed points for non-self mappings in $p$-vector spaces for $0 < p \leq 1$ under different boundary conditions in
locally $p$-convex spaces.

\vskip.1in
\noindent
{\bf Theorem 6.3 (Fixed Points of non-self mappings in locally $p$-convex space).} Let $U$ be a bounded open $p$-convex subset of a locally $p$-convex space $E$ (where, $0 \leq p \leq 1)$ with the zero $0 \in U$, and $C$ a closed $p$-convex subset of $E$ with also zero $0\in C$, and assume
$F: \overline{U}\cap C \rightarrow C$ is a 1-set contractive and continuous mapping, and satisfying the condition (H) or (H1) above. In addition, for each $x \in \partial_C(U)$, $P^{\frac{1}{p}}_U(F(x))- 1 \leq P^{\frac{1}{p}}_U (F(x)-x)$ for $0< p \leq 1$ (this is trivial when $p=1$), where $P_U$ is the Minkowski $p$-functional of $U$.
If $F$ satisfies any one of the following conditions for any $x \in \partial_C(U) \diagdown F(x)$:

(i)  $P_U(F(x)-z) < P_U(F(x)-x)$ for some $z \in \overline{I_{\overline{U}}(x)}\cap C$;

(ii) There exists $\lambda $ with $|\lambda| < 1$ such that $\lambda x + (1-\lambda)F(x)
 \in \overline{I_{\overline{U}}(x)}\cap C$;

(iii) $F(x) \in  \overline{I_{\overline{U}}(x)}\cap C$;

(iv) $F(x) \in \{\lambda x: \lambda > 1 \} =\emptyset$;

(v) $F(\partial U) \subset \overline{U} \cap C$;

(vi) $P_U(F(x)-x) \neq ((P_U(F(x))^{\frac{1}{p}}-1)^p$;


then $F$ must has a fixed point.

\noindent
{\bf Proof.} By following the argument and symbols used in the proof of Theorem 5.2 (see also Theorem 5.4), we have that either

(1) $F$ has a fixed point $x_0 \in U \cap C$; or

(2) there exists $x_0 \in \partial_C(U)$, with $x_0=f(F(x_0))$ such that
$$ P_U (F(x_0) - x_0) = d_P(F(x_0), \overline{U}\cap C) = d_p(F(x_0), \overline{I_{\overline{U}}( x_0)} \cap C) = P_U(F(x_0)) - 1 > 0,$$
where $\partial_C(U)$ denotes the boundary of $U$ relative to $C$ in $E$, and $f$ is the restriction of the continuous retraction $r$ respect to the set $U$ in $E$.

First, suppose that $F$ satisfies the condition (i), if $F$ has no fixed point, then above (2) holds and $x_0 \neq F(x_0)$. Then by the condition (i), it follows that $P_U(F(x_0)-z) < P_U(F(x_0)-x_0)$ for some $z \in \overline{I_{\overline{U}}(x)}\cap C$, this contradicts with the best approximation equations given by (2) above, thus $F$ mush have a fixed pint.

Second, suppose that $F$ satisfies the condition (ii), if $F$ has no fixed point, then above (2) holds and $x_0 \neq F(x_0)$. Then by
the condition (ii), there exists $\lambda >1$ such that $\lambda  x_0 + (1 - \lambda ) F(x_0) \in
\overline{I_{\overline{U}}(x)}\cap C$. It follows that
$$P_U(F(x_0)- x_0) \leq P_U(F(x_0)- (\lambda x_0 + (1-\lambda F(x_0))=P_U(\lambda (F(x_0) - x_0))=|\lambda|^pP_U(F(x_0)-x_0)<P_U(F(x_0)-x_0)$$
this is impossible and thus $F$ must have a fixed point in $\overline{U}\cap C$.

Third, suppose that $F$ satisfies the condition (iii), i.e., $F(x) \in \overline{I_{\overline{U}}(x)} \cap C$;, then the (2), we have that
$P_U (F(x_0) - x_0)$ and thus $x_0=  F(x_0)$, which means $F$ has a fixed point.

Forth, suppose that $F$ satisfies the condition (iv),  if  if $F$ has no fixed point, then above (2) holds and $x_0 \neq  F(x_0)$.
As given by the proof of Theorem 5.2, we have that $F(x_0) \notin \overline{U}$, thus $P_U(F(x_0)) > 1$ and $x_0= f(F(x_0))=\frac{F(x_0)}{(P_U(F(x_0)))^{\frac{1}{p}}}$, which means $F(x_0)=(P_U(F(x_0)))^{\frac{1}{p}} x_0$,
where $(P_U(F(x_0)))^{\frac{1}{p}} > 1$, this contradicts with the assumption (iv), thus $F$ must have a fixed point in $\overline{U} \cap C$.

Fifth, suppose that $F$ satisfies the condition (v), then $x_0 \neq F(x_0)$.
As $x_0 \in \partial_C{U}$, now by the condition (v), we have that $F(\partial U) \subset \overline{U} \cap C$, it follows that for
$F(x_0)$, we have $F(x_0)\in \overline{U}\cap C$, thus $F(x_0) \notin \overline{U} \diagdown \cap C$, which implies that
$0 < P_U(F(x_0)- x_0) = d_P(F(x_0), \overline{U}\cap C) = 0$, this is impossible, thus $F$ must have a fixed point.
Here, like pointed out by Remark 5.2, we know that  based on the condition (v), the mapping $F$ has a fixed point by applying $F(\partial U) \subset \overline{U} \cap C$ is enough, not needing the general hypothesis: $``$for each $x \in \partial_C(U)$, $P^{\frac{1}{p}}_U(F(x))- 1 \leq P^{\frac{1}{p}}_U (F(x)-x)$ for $0< p \leq 1$".

Finally, suppose that $F$ satisfies the condition (vi), if $F$ has no fixed point, then above (2) holds and $x_0 \neq F(x_0)$. Then the condition (v) implies that $P_U(F(x_0)- x_0) \neq ((P_U(F(x_0))^{\frac{1}{p}}-1)^p$, but the our proof in Theorem 5.2 shows that
$P_U(F(x_0)- x_0)=((P_U(F(x_0)))^{\frac{1}{p}}-1)^p$, this is impossible, thus $F$ must have a fixed point.
Then the proof is complete. $\square$

\vskip.1in
Now by taking the set $C$ in Theorem 6.1 as the whole locally $p$-convex space $E$ itself, we have the following general results for non-self
continuous mappings which include results of Rothe, Petryshyn, Altman and Leray-Schauder types' fixed points as special cases in locally convex spaces.

\vskip.1in
Taking $p=1$ and $C =E$ in Theorem 6.3, we have the following fixed points for non-self single-valued mappings
in locally convex spaces (LCS), and the corresponding results for upper semicontinuous set-valued mappings are discussed by Yuan \cite{yuan2022} and related references wherein.

\vskip.1in
\noindent
{\bf Theorem 6.4 (Fixed Points of non-self mappings with boundary conditions).} Let $U$ be a bounded open convex subset of the LCS $E$ with the zero $0 \in U$, and assume $F: \overline{U} \rightarrow E$ is a 1-set contractive and continuous mapping, and satisfying the condition (H) or (H1) above. If $F$ satisfies any one of the following conditions for any $x \in \partial(U) \diagdown F(x)$

(i) $P_U(F(x)-z) < P_U(F(x)-x)$ for some $z \in \overline{I_{\overline{U}}(x)}$; 

(ii) there exists $\lambda $ with $|\lambda| < 1$ such that $\lambda x + (1-\lambda)F(x)
 \in \overline{I_{\overline{U}}(x)}$; 

(iii) $F(x) \in \overline{I_{\overline{U}}(x)}$; 

(iv) $F(x) \in \{\lambda x: \lambda > 1 \} =\emptyset$;

(v) $F(\partial(U) \subset \overline{U}$;

(vi) $P_U(F(x)-x) \neq P_U(F(x))-1$;


then $F$ must has a fixed point.

\vskip.1in
In what follow, based on the best approximation theorem in $p$-seminorm space, we will also give some fixed point theorems for non-self mappings
with various boundary conditions which are related to the study for the existence of solutions for PDE and differential equations with boundary problems (see, Browder \cite{bro1968}, Petryshyn \cite{petryshyn1966}-\cite{petryshyn1973tams}, Reich \cite{reich}),  which would play roles in nonlinear analysis for $p$-seminorm space as shown below.

First, as discussed by Remark 5.2, the proof of Theorem 5.2, with the strongly boundary condition $``$$F(\partial(U)) \subset \overline{U} \cap C$" only, we can prove that  $F$ has a fixed point, thus we have the following fixed point theorem of Rothe type
in $p$-vector spaces.

\vskip.1in
\noindent
{\bf Theorem 6.5 (Rothe Type).}
Let $U$ be a bounded open $p$-convex subset of a locally $p$-convex space $E$ (where, $ 0 \leq p \leq 1)$ with the zero $0 \in U$.
Assume $F: \overline{U}\rightarrow E$ is a 1-set contractive and continuous mapping, satisfying the condition (H) or (H1) above, and such that $F(\partial(U)) \subset \overline{U}$, then $F$ must has a fixed point.

\vskip.1in
Now as applications of Theorem 6.5, we give the following Leray-Schauder alternative in locally $p$-convex spaces for non-self mappings associated with the boundary condition  which often appear in the applications (see Isac \cite{isac} and references therein for the study of complementary problems and related topics in optimization).

\vskip.1in
\noindent
{\bf Theorem 6.6 (Leray-Schauder Alternative in locally $p$-Convex Spaces).} Let $E$ be a locally $p$-convex space $E$, where $0 < p \leq 1$, $B \subset E$ a bounded closed $p$-convex such that $0 \in int B$. Let $F: [0, 1] \times B  \rightarrow E$  be 1-set contractive and continuous, satisfying the condition (H) or (H1) above, and such that
the set $F([0, 1] \times B)$ be relatively compact in $E$. If the following assumptions are satisfied:

(1) $x \neq  F(t, x)$ for all $x \notin \partial B$ and $t \in [0, 1]$,

(2) $F(\{0\} \times \partial B) \subset B$,

\noindent
then there is an element $x^* \in B$ such that $x^* = F(1, x^*)$.

\noindent
{\bf Proof.} For $n \in \mathbb{N}$, we consider the mapping
\begin{equation}\label{eq10}
F_n(x)=\left\{
\begin{aligned}
&F(\frac{1-P_B(x)}{\epsilon_n}, \frac{x}{P_B(x)}), &\mbox{ if }&  1-\epsilon \leq P_B(x)\leq 1,\\
&F(1, \frac{X}{1-\epsilon_n}),                     &\mbox{ if }&  P_B(x) < 1 - \epsilon_n, \\
\end{aligned}
\right.
\end{equation}
where $P_B$ is the Minkowski $p$-functional of $B$ and $\{\epsilon_n\}_{n \in N}$ is  a sequence of
real numbers such that $\lim_{n \rightarrow \infty} \epsilon_n=0$ and $0 < \epsilon_n < \frac{1}{2}$ for any $n \in \mathbb{N}$, and we also  observe that, the mapping $F_n$ is 1-set contractive continuous with non-empty closed $p$-convex values  on $B$. From assumption (2), we have that $F_n(\partial B) \subset B$, and the assumptions of Theorem 6.5 are satisfied, then for each $n \in \mathbb{N}$, there exists an element $u_n \in B$ such that $u_n = F_n(u_n)$.

We first prove the following statement: ``It is impossible to have an infinite number of the elements $u_n$ satisfy the following inequality: $1 - \epsilon_n  \leq  P_B(u_n) \leq 1$. $"$

If not, we assume to have an infinite number of the elements $u_n$ satisfy the following inequality:
$$1 - \epsilon_n  \leq P_B(u_n) \leq 1.$$
As $F_n(B)$ is relatively compact and by the definition of mappings $F_n$, we have that $\{u_n\}_{n \in \mathbb{N}}$  is contained in a compact set in $E$. Without loss of the generality (indeed, each compact set is also countably compact), we define the sequence $\{t_n\}_{n \in \mathbb{N}}$ by $t_n: =\frac{1-P_B(u_n)}{\epsilon}$ for each $n \in N$.
Then we have that $\{t_n\}_{n\in N}\subset [0, 1]$ and we may assume that $\lim_{n \rightarrow \infty}t_n = t \in [0, 1]$.
The corresponding subsequence of $\{u_n\}_{n \in \mathbb{N}}$ is denoted again by $\{u_n\}_{n \in \mathbb{N}}$ and it also satisfies the
inequality: $1-\epsilon_n \leq P_B(u_n) \leq 1$, which implies that $\lim_{n\rightarrow \infty} P_B (u_n)=1$.

Now let $u^*$ be an accumulation point of $\{u_n\}_{n\in N}$, thus  have $\lim_{n \rightarrow \infty}(t_n,\frac{u_n}{P_B(u_n)}, u_n) = (t, u^*, u^*)$.
By the fact that $F$ is compact, we have assume that $u_n = F(t_n, \frac{u_n}{P_B(u_n)})$ for each $n \in \mathbb{N}$,
it follows that $u^* = F(t, u^*)$, this contradicts with the assumption (1) as we have $\lim_{n \rightarrow \infty}P_B(u_n)=1$ (which means that $u^* \in \partial B$, this is impossible).

Thus it is impossible to have that $``$to have an infinite number of elements $u_n$ satisfy the inequality:
         $1 - \epsilon_n \leq P_B(u_n) \leq 1$" ,
which means that there is only a finite number of elements of sequence $\{u_n\}_{n \in N}$ satisfying the inequality:  $ 1 - \epsilon_n \leq P_B(u_n) \leq 1$. Now, without loss of the generality, for $n \in \mathbb{N}$, we have the following inequality:
$$ P_B(u_n) < 1 - \epsilon_n.$$
By the fact that $\lim_{n \rightarrow} (1-\epsilon_n)=1$, $u_n \in F(1, \frac{u_n}{1-\epsilon})$ for all $n \in \mathbb{N}$ and
assume that $\lim_{n\rightarrow} u_n = u^*$, then the continuity of $F$ with non-empty closed values implies that by
$u_n = F(1, \frac{u_n}{1-\epsilon})$ for each $n \in \mathbb{N}$, it implies that  $u^* = F(1, u^*)$. This completes the proof. $\square$

\vskip.1in
As a special case of Theorem 6.6, we have the following principle for the implicit form of Leray-Schauder type alternative  in locally $p$-convex spaces for $0< p \leq 1.$

\noindent
{\bf Corollary 6.1 (The Implicit Leray-Schauder Alternative).} Let $E$ be  a locally $p$-convex space $E$,
where $0 < p \leq 1$, $B \subset E$ a bounded closed $p$-convex such that $ 0 \in int B$. Let $F: [0, 1] \times B \rightarrow E$  be 1-set contractive and continuous, satisfying the condition (H) or (H1) above, and the set $F([0, 1] \times B)$ be relatively compact in $E$. If the following assumptions are satisfied:

(1) $F(\{0\} \times \partial B) \subset B$,

(2) $x \neq F(0, x)$ for all $x \in \partial B$,

\noindent
then at least one of the following properties is satisfied:

(i) there exists $x^* \in  B$ such that $x^* = F(1, x^*)$; or

(ii) there exists  $(\lambda^*, x^*) \in (0, 1) \times \partial B$ such that $x^* = F(\lambda^*, x^*)$.

\noindent
{\bf Proof.}  The result is an immediate consequence of Theorem 6.6, this completes the proof. $\square$
\vskip.1in

We like to point out the similar results on Rothe and Leray-Schauder alternative have been developed by Furi and Pera \cite{furipera},
Granas and Dugundji \cite{granas}, G\'{o}rniewicz \cite{gorniewicz}, G\'{o}rniewicz et al.\cite{gorniewiczetal},
Isac \cite{isac}, Li et al.\cite{lixuduan2006}, Liu \cite{liu2001}, Park \cite{park1995}, Potter \cite{potter1972}, Shahzad \cite{shahzad2006}-\cite{shahzad2004}, Xu \cite{xu2007}, Xu et al.\cite{xujiali2006}, and related references therein as tools of nonlinear analysis in the Banach space setting and applications to the boundary value problems for ordinary differential equations in noncompact problems, a general class of mappings for nonlinear alternative of Leray-Schauder type in normal topological spaces, and some Birkhoff-Kellogg type theorems for general class mappings in topological vector spaces are also established by Agarwal et al.\cite{Agarwal}, Agarwal and O'Regan \cite{agarwalorgan2003}-\cite{agarwalorgan2004}, Park \cite{park1997}, and references therein for more in detail; and in particular, recently O'Regan \cite{oregan2019} uses the Leray-Schauder type coincidence theory to establish some Birkhoff-Kellogg problem, Furi-Pera type results for a general class of 1-set contractive mappings.

Before closing this section,  we like to share with readers that as the application of the best approximation result for 1-set contractive mappings, we just establish some fixed point theorems and general principle of Leray-Schauder alternative for non-self mappings, which seem would play important roles for the nonlinear analysis under the framework of locally $p$-convex (seminorm) spaces, as the achievement of nonlinear analysis under the framework for underling being locally topological vector spaces, normed spaces, or in Banach spaces.

\vskip.1in
\section{Fixed Points for the class of 1-Set Contractive Mappings}

In this section, based on the best approximation Theorem 5.2 for classes of 1-set contractive mappings developed in section 5, we will show how it can be used as a useful tool to establish fixed point theorems for  non-self upper semi-continuous mappings in locally $p$-convex (seminorm)
spaces for $p \in (0, 1]$, which include norm spaces, uniformly convex Banach spaces as special classes.

By following Browder \cite{bro1968}, Li \cite{li1988}, Goebel and Kirk \cite{goebelkirk1990}, Petryshyn \cite{petryshyn1966}-\cite{petryshyn1973tams}, Tan and Yuan \cite{tanyuan1994}, Xu \cite{xu2007} and references therein, we recall some definitions as follows for $p$-seminorm spaces, where $p \in (0, 1]$.

\vskip.1in
\noindent
{\bf Definition 7.1}. Let $D$ be a non-empty (bounded) closed subset of locally $p$-convex spaces $(E, \|\cdot\|_p)$, where $p \in (0, 1]$. Suppose $f: D \rightarrow X$ is a (single-valued) mapping, then: (1) $f$  is said to be nonexpansive if for each $x, y \in D$, we have $\|f(x)- f(y)\|_p \leq \|x-y\|_p$;
(2) $f$ (actually, $(I-f)$) is said to be demiclosed (see Borwder \cite{bro1968}) at $y \in X$ if for any sequence $\{x_n\}_{n \in \mathbb{N}}$ in $D$, the conditions $x_n \rightarrow x_0\in D$ weakly, and $(I-f)(x_n) \rightarrow y_0$  strongly imply that
 $(I-f)(x_0)=y_0$, where $I$ is the identity mapping;
(3) $f$ is said to be hemicompact (see p.379 of Tan and Yuan \cite{tanyuan1994}) if each sequence $\{x_n\}_{n \in \mathbb{N}}$ in $D$ has a convergent subsequence with the limit $x_0$ such that $x_0 = f(x_0)$,  whenever $\lim_{n \rightarrow \infty}d_p(x_n, f(x_n))=0$, here  $d_P(x_n, f(x_n)):=\inf\{P_U(x_n- z): z \in f(x_n)\}$, and $P_U$ is the Minkowski $p$-functional for any $U \in \mathfrak{U}$, which is the family of all non-empty open $p$-convex subset containing the zero in $E$;
(4) $f$ is said to to be demicompact (by  Petryshyn \cite{petryshyn1966})
if each sequence $\{x_n\}_{n \in \mathbb{N}}$ in $D$  has
a convergent subsequence whenever $\{x_n -f(x_n)\}_{n \in \mathbb{N}}$  is a convergent sequence in $X$;
(5) $f$ is said to be a semi-closed 1-set contractive mapping if $f$ is 1-set contractive mapping, and
$(I-f)$ is closed, where $I$ is identity mapping (by Li \cite{li1988}); and
(6) $f$ is said to be semicontractive (see Petryshyn \cite{petryshyn1973tams} and Browder \cite{bro1968}) if there exists a
mapping $V: D \times D \rightarrow 2^X$  such that $f(x) = V(x, x)$ for each $x \in D$, with (a) for each fixed
$x \in D$, $V(\cdot, x)$ is nonexpansive from $D$ to $X$; and (b) for each fixed $x\in D$, $V(x, \cdot)$  is completely continuous
from $D$ to $X$, uniformly for $u$ in a bounded subset of $D$ (which means if $v_j$ converges weakly to $v$ in
$D$ and $u_j$ is a bounded sequence in $D$, then $V(u_j, v_j) - V(u_j, v) \rightarrow 0$, strongly in $D$).

\vskip.1in
From the definition above, we first observe that the definitions (1) to (6) for set-valued mappings can be given by the similar way with the Hausdorff metric $H$ (we omit their definitions here in details by saving spaces); Secondly, if $f$ is a continuous demicompact mapping, then $(I - f)$ is closed, where $I$ is the identity mapping on $X$. it is also clear from definitions that every demicompact map is hemicompact in seminorm spaces, but the converse is not true by the example in p.380 by Tan and Yuan \cite{tanyuan1994}.
It is evident that if $f$ is demicompact, then $I-f$ is demiclosed. It is know that for each condensing mapping $f$, when $D$ or $f(D)$ is bounded, then $f$ is hemicompact; and also $f$ is demicompact in metric spaces by Lemma 2.1 and Lemma 2.2 of Tan and Yuan \cite{tanyuan1994}, respectively. In addtion, it is known that every nonexpansive map is a 1-set-contractive mapping; and also if $f$ is a hemicompact 1-set-contractive mapping, then $f$ is a 1-set-contractive mapping satisfying the following {\bf (H1) condition} (which is the same as the $``$condition (H1)" in Section 5, but slightly different from the condition (H) used there in the Section 5):

{\bf (H1) condition}: Let $D$ be a nonempty bounded subset of a space $E$ and assume $F: \overline{D} \rightarrow 2^E$ a set-valued mapping. If $\{x_n\}_{n \in \mathbb{N}}$ is any sequence in $D$ such that for each $x_n$, there exists $y_n \in F(x_n)$ with $\lim_{n \rightarrow \infty} (x_n- y_n)=0$, then there exists a point
$x\in \overline{D}$ such that $x \in F(x)$.

\vskip.1in
We first note that the $``$(H1) Condition" above is actually the same one as the $``$Condition (C)"  used by Theorem 1 of  Petryshyn \cite{petryshyn1973tams}. Secondly, it was shown by Browder  \cite{bro1968} that indeed the nonexpansive mapping in a uniformly convex Banach $X$ enjoys the condition (H1) as shown below.

\vskip.1in
\noindent
{\bf Lemma 7.1}. Let $D$ be a nonempty bonded convex subset of a uniformly convex Banach space $E$. Assume $F: \overline{D} \rightarrow E$ is a  nonexpansive (single-valued) mapping, then the mapping $P: = I - F$  defined by $P(x): = (x-F(x)) $  for each $x \in  \overline{D}$ is demiclosed, and in particular, the $``$(H1) Condition" holds.

\noindent
{\bf Proof.}By following the argument given in p.329 (see the proof of Theorem 2.2 and Corollary 2.1) by Petryshyn \cite{petryshyn1973tams}, the mapping $F$ is demiclosed (which actually is called Browder's demiclosednedd principle), which says that by the assumption of (H1) condition, If $\{x_n\}_{n \in \mathbb{N}}$ is any sequence in $D$ such that for each $x_n$, there exists $y_n \in F(x_n)$ with $\lim_{n \rightarrow \infty} (x_n- y_n)=0$, then we have $0 \in (I - F) (\overline{D})$, which means that there exists $x_0 \in \overline{D}$ with $0 \in (I-F)(x_0)$, this implies that  $x_0 \in F(x_0)$. The proof is complete. $\square$.

\vskip.1in
\noindent
{\bf Remark 7.1.} When a $p$-vector space $E$ is with a $p$-norm, then $``$(H) condition"  satisfies the $``$(H1) condition".
The (H1) condition mainly supported by the so-called demiclosedness principle after the work by Browder \cite{bro1968}.

\vskip.1in
By applying Theorem 5.2, we have the following result for non-self mappings in $p$-seminorm spaces for $p \in (0, 1]$.

\vskip.1in
\noindent
{\bf Theorem 7.1}. Let $U$ be a bounded open $p$-convex subset of a locally $p$-convex (or seminorm)  space $E$ ($0 < p \leq 1)$  the zero $0 \in U$.  Assume $F: \overline{U} \rightarrow E$ is a 1-set contractive and continuous mapping, satisfying the condition (H) or (H1) above. In addition, for any $x\in \partial \overline{U}$, we have $\lambda x \neq F(x)$ for any $\lambda > 1$ (i.e., the $``$Leray-Schauder boundary condition"), then $F$ has at least one fixed point.

\noindent
{\bf Proof.} By Theorem 5.2 with $C= E$, it follows that we have the following either (I) or (II) holding:

(I)  $F$ has a fixed point $x_0 \in U $, i.e., $P_U (F(x_0) - x_0) = 0$,

(II)  there exists $x_0 \in \partial(U)$ with $ P_U (F(x_0) - x_0) = (P^{\frac{1}{p}}_U(F(x_0))-1)^p  > 0.$

If $F$ has no fixed point, then above (II) holds and $x_0 \neq  F(x_0)$. By the proof of Theorem 5.2, we have that $x_0=f(F(x_0))$ and
$F(x_0) \notin \overline{U}$. Thus $P_U(F(x_0)) > 1$ and $x_0= f(F(x_0))=\frac{F(x_0)}{(P_U(F(x_0))^{\frac{1}{p}}}$, which means $F(x_0)=(P_U(F(x_0)))^{\frac{1}{p}} x_0$, where $(P_U(F(x_0)))^{\frac{1}{p}} > 1$, this contradicts with the assumption. Thus $F$ must have a fixed point. The proof is complete. $\square$

\vskip.1in
By following the idea used and developed by Browder \cite{bro1968}, Li \cite{li1988}, Li et al.\cite{lixuduan2006}, Goebel and Kirk \cite{goebelkirk1990}, Petryshyn \cite{petryshyn1966}-\cite{petryshyn1973tams}, Tan and Yuan \cite{tanyuan1994}, Xu \cite{xu2007}, Xu et al.\cite{xujiali2006} and references therein,  we have the following a number of existence theorems for the principle of Leray-Schauder type alternatives in locally $p$-convex spaces, or $p$-seminorm spaces  $(E, \|\cdot \|_p)$  for $p \in (0, 1]$.

\vskip.1in
\noindent
{\bf Theorem 7.2}. Let $U$ be a bounded open $p$-convex subset of a $p$-seminorm space $(E, \|\cdot \|_p)$ ($0 < p \leq 1)$  the zero $0 \in U$.  Assume $F: \overline{U} \rightarrow E$ is a 1-set contractive and continuous mapping, satisfying the condition (H) or (H1) above.
In addition, there exist $\alpha >1$, $\beta \geq 0$, such that for each $x \in \partial \overline{U}$, we have that
for any $y \in F(x)$, $\|y- x\|_p^{\alpha/p}\geq \|y\|_p^{(\alpha+\beta)/p}\|x\|_p^{-\beta/p}  - \|x\|_p^{\alpha/p}.$
Then $F$ has at least one fixed point.

\noindent
{\bf Proof.} We prove the conclusion by showing the Leray-Schauder boundary condition in Theorem 7.1 does not hold.
If we assume $F$ has no fixed point, by the boundary condition of Theorem 7.1,
there exist $x_0\in \partial \overline{U}$, $\lambda_0 >1$ such that $F(x_0) = \lambda_0 x_0$.

Now, consider the function $f$ defined by $f(t): =(t-1)^{\alpha} - t^{\alpha + \beta}+1$ for $t\geq 1$. We observe that
$f$ is a strictly decreasing function for $t \in [1, \infty)$ as the derivative of $f$\'{}$(t) =\alpha (t-1)^{\alpha-1} - (\alpha +\beta) t^{\alpha +\beta -1} < 0$  by the differentiation, thus we have
$t^{\alpha + \beta} -1 > (t-1)^{\alpha}$ for $t \in (1, \infty)$. By combining the  boundary condition, we have that
$\|F(x_0)-x_0\|_p^{\alpha/p}=\|\lambda_0 x_0-x_0\|_p^{\alpha/p}=(\lambda_0-1)^{\alpha}\|x_0\|_p^{\alpha/p} < (\lambda_0^{\alpha+\beta}-1)\|x_0\|_p^{(\alpha+\beta)/p}\|x_0\|_p^{-\beta/p}=\|F(x_0)\|_p^{(\alpha+\beta)/p}\|x_0\|_p^{-\beta/p}- \|x_0\|_p^{\alpha/p}$, which contradicts the boundary condition given by Theorem 7.2. Thus, the conclusion follows and the proof is complete. $\square$

\vskip.1in
\noindent
{\bf Theorem 7.3}. Let $U$ be a bounded open $p$-convex subset of a $p$-seminorm space $(E, \|\cdot \|_p)$ ($0 < p \leq 1)$  the zero $0 \in U$.  Assume $F: \overline{U} \rightarrow E$ is a 1-set contractive and continuous mapping, satisfying the condition (H) or (H1) above.
In addition, there exist $\alpha >1$, $\beta \geq 0$, such that for each $x \in \partial \overline{U}$, we have that
 $\|F(x) + x\|_p^{(\alpha+\beta)/p} \leq \|F(x)\|_p^{\alpha/p}\|x\|_p^{\beta/p}  + \|x\|_p^{(\alpha+\beta)/p}.$
Then $F$ has at least one fixed point.

\noindent
{\bf Proof.} We prove the conclusion by showing the Leray-Schauder boundary condition in Theorem 7.1 does not hold.
If we assume $F$ has no fixed point, by the boundary condition of Theorem 7.1,
there exist $x_0\in \partial \overline{U}$ and $\lambda_0 >1$ such that $F(x_0)= \lambda_0 x_0$.

Now, consider the function $f$ defined by $f(t): =(t+1)^{\alpha+\beta} - t^{\alpha} - 1 $ for $t\geq 1$. We then can show that
$f$ is a strictly increasing function for $t \in [1, \infty)$,  thus we have
$t^{\alpha}+1 < (t + 1)^{\alpha +\beta}$ for $t \in (1, \infty)$.
By the  boundary condition given in Theorem 7.3, we have that
$$\|F(x_0) + x_0\|_p^{(\alpha+\beta)/p}=(\lambda_0 +1)^{\alpha+\beta}\|x_0\|_p^{(\alpha+\beta)/p} > (\lambda_0^{\alpha}+1)\|x_0\|_p^{(\alpha+\beta)/p}=\|F(x_0)\|_p^{\alpha/p}\|x_0\|_p^{\beta/p}+ \|x_0\|_p^{\alpha/p},$$
which contradicts the boundary condition given by Theorem 7.3. Thus, the conclusion follows and the proof is complete. $\square$

\vskip.1in
\noindent
{\bf Theorem 7.4}. Let $U$ be a bounded open $p$-convex subset of a $p$-seminorm space $(E, \|\cdot \|_p)$ ($0 < p \leq 1)$  the zero $0 \in U$.  Assume $F: \overline{U} \rightarrow E$ is a 1-set contractive and continuous mapping, satisfying the condition (H) or (H1) above.
In addition, there exist $\alpha >1$, $\beta \geq 0$ (or alternatively, $\alpha >1$, $\beta \geq 0$) such that for each $x \in \partial \overline{U}$, we have that
 $\|F(x) - x\|_p^{\alpha/p} \|x\|_p^{\beta/p} \geq \|F(x)\|_p^{\alpha/p}\|F(x)+x\|_p^{\beta/p} -\|x\|_p^{(\alpha+\beta)/p}.$
 Then $F$ has at least one fixed point.

\noindent
{\bf Proof.} The same as above, we prove the conclusion by showing the Leray-Schauder boundary condition in Theorem 7.1 does not hold. If we assume $F$ has no fixed point, by the boundary condition of Theorem 7.1,
there exist $x_0\in \partial \overline{U}$, and $\lambda_0 >1$ such that $F(x_0) = \lambda_0 x_0$.

Now, consider the function $f$ defined by $f(t): =(t-1)^{\alpha} - t^{\alpha}(t-1)^{\beta}+1$ for $t\geq 1$. We then can show that
$f$ is a strictly decreasing function for $t \in [1, \infty)$,  thus we have
$(t-1)^{\alpha} < t^{\alpha} (t+1)^{\beta}-1$ for $t \in (1, \infty)$.

By the  boundary condition given in Theorem 7.4, we have that
$$\|F(x_0)-x_0\|_p^{\alpha/p}\|x_0\|_p^{\beta/p}=(\lambda_0-1)^{\alpha}\|x_0\|_p^{(\alpha+\beta)/p} <  (\lambda_0^{\alpha}(\lambda_0+1)^{\beta}-1)\|x_0\|_p^{(\alpha+\beta)/p}=\|F(x_0)\|_p^{\alpha/p}\|F(x_0)+x_0\|_p^{\beta/p}- \|x_0\|_p^{(\alpha+\beta)/p},$$
which contradicts the boundary condition given by Theorem 7.4. Thus, the conclusion follows and the proof is complete. $\square$

\vskip.1in
\noindent
{\bf Theorem 7.5}. Let $U$ be a bounded open $p$-convex subset of a $p$-seminorm space $(E, \|\cdot \|_p)$ ($0 < p \leq 1)$  the zero $0 \in U$.  Assume $F: \overline{U} \rightarrow E$ is a 1-set contractive and continuous mapping, satisfying the condition (H) or (H1) above.
In addition, there exist $\alpha >1$, $\beta \geq 0$, we have that
$\|F(x) + x\|_p^{(\alpha+\beta)/p} \leq \|F(x)-x\|_p^{\alpha/p}\|x\|_p^{\beta/p} +\|F(x)\|_p^{\beta/p} \|x\|^{\alpha/p}.$
Then $F$ has at least one fixed point.

\noindent
{\bf Proof.} The same as above, we prove the conclusion by showing the Leray-Schauder boundary condition in Theorem 7.1 does not hold. If we assume $F$ has no fixed point, by the boundary condition of Theorem 7.1,
there exist $x_0\in \partial \overline{U}$, and $\lambda_0 >1$ such that $F(x_0) = \lambda_0 x_0$.

Now, consider the function $f$ defined by $f(t): =(t+1)^{\alpha+\beta} - (t-1)^{\alpha}-t^{\beta}$ for $t\geq 1$. We then can show that $f$ is a strictly increasing function for $t \in [1, \infty)$,  thus we have
$(t+1)^{\alpha+\beta} >  (t-1)^{\alpha} +t^{\beta}$ for $t \in (1, \infty)$.

By the  boundary condition given in Theorem 7.5, we have that
$\|F(x_0) +x_0\|_p^{(\alpha+\beta)/p}=(\lambda_0 +1)^{\alpha+\beta}\|x_0\|_p^{(\alpha+\beta)/p} >
 ((\lambda_0-1)^{\alpha}+ \lambda_0^{\beta})\|x_0\|_p^{(\alpha+\beta)/p}=\|\lambda_0 x_0- x_0\|_p^{\alpha/p}\|x_0\|_p^{\beta/p} + \|\lambda_0 x_0\|_p^{\beta/p}\|x_0\|_p^{\alpha/p} = \|F(x_0)-x_0\|_p^{\beta/p}\|x_0\|_p^{\alpha/p} +\|F(x_0)\|_p^{\beta/p}\|x_9\|^{\alpha/p},$
which implies that
$$ \|F(x_0) +x_0\|_p^{(\alpha+\beta)/p} >  \|F(x_0)-x_0\|_p^{\beta/p}\|x_0\|_p^{\alpha/p} +\|F(x_0)\|_p^{\beta/p}\|x_9\|^{\alpha/p},$$
this contradicts the boundary condition given by Theorem 7.5. Thus, the conclusion follows and the proof is complete. $\square$

As an application of Theorems 7.1 by testing the Leray-Schauder boundary condition, we have the following conclusion for each special
case, and thus we omit their proofs in details here.

\vskip.1in
\noindent
{\bf Corollary 7.1}. Let $U$ be a bounded open $p$-convex subset of a $p$-seminorm space $(E, \|\cdot \|_p)$ ($0 < p \leq 1)$  the zero $0 \in U$.  Assume $F: \overline{U} \rightarrow E$ is a 1-set contractive and continuous mapping, satisfying the condition (H) or (H1) above. Then $F$ has at least one fixed point if one of the the following conditions holds for $x \in \partial \overline{U}$:

(i) $\|F(x)\|_p \leq \|x\|_p$,

(ii) $\|F(x)\|_p \leq \|F(x)-x\|_p$,

(iii) $\|F(x)+x||_p \leq \|F(x)\|_p$,

(iv) $\|F(x)+ x\|_p \leq \|x\|_p$,

(v) $\|F(x)+x\|_p \leq \|F(x)- x\|_p$,

(vi) $\|F(x)\|_p \cdot \|F(x)+x\|_p \leq \|x\|_p^2$,

(vii) $\|F(x)\|_p \cdot \|F(x) +x\|_p \leq \|F(x)- x\|_p \cdot \|x\|_p$.

\vskip.1in
If the $p$-seminorm space $E$ is a uniformly convex Banach space $(E, \| \cdot \|)$ (for $p$-norm space with $p=1$),  then we have the following general existence result (which actually is true for non-expansive set-valued mappings).

\vskip.1in
\noindent
{\bf Theorem 7.6}. Let $U$ be a bounded open convex subset of a uniformly convex Banach space $(E, \|\cdot \|)$ (with $p=1$)
with zero $0 \in U$. Assume $F: \overline{U} \rightarrow E$ is a semi-contractive and continuous single-valued mapping with non-empty  values. In addition, for any $x\in \partial \overline{U}$, we have $\lambda x \neq F(x)$ for any $\lambda > 1$ (i.e., the $``$Leray-Schauder boundary condition"). Then $F$ has at least one fixed point.

\noindent
{\bf Proof.} By the assumption that $F$  is a semi-contractive and continuous single-valued mapping with non-empty values, it follows by Lemma 3.2 in p.338 of Petryshyn \cite{petryshyn1973tams}, $f$ is a 1-set contractive single-valued mapping.
Moreover, by the assumption that $E$ is a uniformly convex Banach, indeed $(I-F)$ is closed at zero, i.e., $F$ is semiclosed (see Browder \cite{bro1968}, or Goebel and Kirk \cite{goebelkirk1990}). Thus all assumptions of Theorem 7.1 are satisfied with  the
(H1) condition. The conclusion follows by Theorem 7.1, and the proof is completes. $\square$

\vskip.1in
Like Lemma 7.1 shows that s single-valued nonexpansive mapping defined in a uniformly convex Banach space (see also Theorem 7.6) satisfied the (H1) condition.  Actually, the nonexpansive set-valued  mappings defined on a special class of Banach spaces with the so-called the $``$ Opial's condition" do not only satisfy the condition (H1), but also belong to the classes of semiclosed 1-set contractive mappings as shown below.

Now let $K(X)$ denote the family of all non-empty compact convex subsets of topological vector space $X$.
The notion of the so-called $``$ Opial's condition" first given by Opial \cite{opial1967}, which says that a Banach space $X$ is said
to satisfy Opial's condition if $\liminf_{n \rightarrow \infty} \| w_n - w \| < \liminf_{n\rightarrow \infty} \|w_n-p\|$ whenever $(w_n)$ is a sequence in $X$ weakly convergent to $w$ and $p\neq w$, we know that Opial's condition plays an important role in the fixed point theory, e.g., see Lami Dozo \cite{lamidozo}, Goebel and Kirk \cite{goebelkirk2008}, Xu \cite{xu2000}
and references where. The following result shows that there nonexpansive set-valued mappings in Banach spaces with Opial's condition (see Lami Dozo \cite{lamidozo} satisfying the condition (H1).

\vskip.1in
\noindent
{\bf Lemma 7.2}. Let $C$ be a convex weakly compact of a Banach space $X$ which satisfies Opial's condition.
Let $T: C \rightarrow K(C)$ be a non-expansive set-valued mapping with non-empty compact values.
Then the graph of $(I-T)$ is closed in $(X, \sigma(X, X^*) \times (X, \|\cdot\|))$, thus $T$ satisfies the $``$(H1) condition", where, $I$ denotes the identity on $X$, $\sigma(X, X^*)$ the weak topology, and $\|\cdot\|$  the  norm (or strong) topology.

\noindent
{\bf Proof.} By following Theorem 3.1 of Lami Dozo \cite{lamidozo}, it follows that the mapping $T$ is demiclosed, thus
$T$ satisfies the $``$(H1) condition". The proof is complete. $\square$

As an application of Lemma 7.2, we have the following results for non-expansive mappings.

\vskip.1in
\noindent
{\bf Theorem 7.7}.
Let $C$ is a nonempty convex weakly compact subset of a Banach space $X$ which satisfies Opial's condition and $0 \in intC$.
Let $T:  C \rightarrow K(X)$  be a nonexpansive set-valued mapping with non-empty compact convex values.
In addition, for any $x\in \partial \overline{C}$, we have $\lambda x \neq F(x)$ for any $\lambda > 1$ (i.e., the $``$Leray-Schauder boundary condition"). Then $F$ has at least one fixed point.

\noindent
{\bf Proof.} As $T$ is nonexpansive, it is 1-set contractive, By Lemma 7.1, it is then semi-contractive and continuous.
Then the (H1) condition of Theorem 7.1 is satisfied. The conclusion follows by Theorem 7.1, and the proof is complete. $\square$.

\vskip.1in
Before the end of this section, by considering the $p$-seminorm space $(E, \|\cdot\|)$ is a seminorm space with $p=1$, the following result is a special case of corresponding results from Theorem 7.2 to Theorem 7.5, and thus we omit its proof.

\vskip.1in
\noindent
{\bf Corollary 7.2}. Let $U$ be a bounded open convex subset of a norm space $(E, \|\cdot\|)$.
 Assume $F: \overline{U} \rightarrow E$ is a 1-set contractive and continuous mapping, satisfying the condition (H) or (H1) above. Then $F$ has at least one fixed point if there exist $\alpha >1$, $\beta \geq 0$, such that any one of the following conditions satisfied

(i) for each $x \in \partial \overline{U}$, $\|F(x)- x\|^{\alpha}\geq \|F(x)\|^{(\alpha+\beta)}\|x\|^{-\beta}  - \|x\|^{\alpha}$,

(ii) for each $x \in \partial \overline{U}$, $\|F(x) + x\|^{(\alpha+\beta)} \leq \|F(x)\|^{\alpha}\|x\|^{\beta}  + \|x\|^{(\alpha+\beta)}$,

(iii) for each $x \in \partial \overline{U}$, $\|F(x) - x\|^{\alpha} \|x\|^{\beta} \geq \|F(x)\|^{\alpha}\|F(x)+x\|^{\beta} -\|x\|^{(\alpha+\beta)}$,

(iv) for each $x \in \partial \overline{U}$, $\|F(x) + x\|^{(\alpha+\beta)} \leq \|F(x)-x\|^{\alpha}\|x\|^{\beta} +\|F(x)\|^{\beta} \|x\|^{\alpha}$.

\vskip.1in
\noindent
{\bf Remark 7.2.} As discussed by Lemma 7.1 and the proof of Theorem 7.6, when the $p$-vector space is a uniformly convex Banach space, the semi-contractive or nonexpansive mappings automatically satisfy the condition (H) or (H1).
Moreover, our results from Theorem 7.1 to Theorem 7.6, Corollary 7.1 and Corollary 7.2 also improve or unify corresponding results given by Browder \cite{bro1968}, Li \cite{li1988}, Li et al.\cite{lixuduan2006}, Goebel and Kirk \cite{goebelkirk1990}, Petryshyn \cite{petryshyn1966}-\cite{petryshyn1973tams}, Reich \cite{reich}, Tan and Yuan \cite{tanyuan1994}, Xu \cite{xu1991}, Xu \cite{xu2007}, Xu et al.\cite{xujiali2006}, and results from the reference therein by extending the non-self mappings to the classes of 1-set contractive set-valued mappings in $p$-seminorm spaces with $p \in (0.1]$ (including the normed space or Banach space when $p=1$, and for $p$-seminorm spaces).

\vskip.1in
\section{Fixed Points for the class of Semiclosed 1-Set Contractive Mappings in $p$-seminorm spaces}

In order to study the fixed point theory for a class of semiclosed 1-set contractive mappings in $p$-seminorm spaces, we first introduce the following definition which is a set-valued generalization of single-value semiclosed 1-set mappings first discussed by Li \cite{li1988}, Xu \cite{xu2007} (see also Li et al.\cite{lixuduan2006}, Xu et al.\cite{xujiali2006} and references therein).

\vskip.1in
\noindent
{\bf Definition 8.1}. Let $D$ be a non-empty (bounded) closed subset of $p$-vector spaces $(E, \|\cdot\|_p)$ with $p$-seminorm for $p$-vector spaces, where $p \in (0, 1]$ (which include norm space, or Banach spaces as special classes), and suppose $T: D \rightarrow X$ is a set-valued mapping. Then $F$ is said to be a semiclosed 1-set contraction mapping if $T$ is 1-set contraction, and $(I-T)$ is closed, which means that for a given net $\{x_n\}_{i \in I}$, for each $i \in I$, there exists $y_i \in T(x_i)$ with $\lim_{i \in I} (x_i - y_i)=0$, then $0 \in (I-T)(\overline{D})$, i.e., there exists $x_0 \in \overline{D}$ such that $x_0 \in T(x_0)$.

\vskip.1in
\noindent
{\bf Remark 8.1}. By Lemma 7.1 and Lemma 7.2 above, it follows that each non-expansive (single-valued) mapping defined on a subset of uniformly convex Banach spaces, and nonexpansive set-valued mappings defined on a subset of Banach spaces satisfying Opial's condition are semiclosed 1-set contractive mapping (see also Goebel \cite{goebel}, Goebel and Kirk \cite{goebelkirk1990},
Petrusel et al.\cite{petrusel2014}, Xu \cite{xu2000}, Yangai \cite{yangai1980} for related discussion by the reference therein).
In particular, under the setting of metric spaces or Banach spaces with certain property, it is clear that each semiclosed 1-set contractive mapping satisfies the condition (H1) above.

Though we know that compared to the single-valued case, based on the  study in the literature about the approximation of fixed points for multi-valued mappings, a well-known counterexample due to Pietramala \cite{pietramala} (see also  Muglia and Marino
\cite{mugliamarino}) proved in 1991 that Browder approximation Theorem 1 given by Browder \cite{bro1967} cannot be extended to the genuine multivalued case even on a finite dimensional space $\mathbb{R}^2$. Moreover, if a Banach space $X$ satisfies Opial's property (see Opial \cite{opial1967}) that is, if $x_n$ weekly converges to $x$, then we have that, $\limsup \|x_n-x\| < \limsup \|x_n -y\|$ for all $x \in X$ and $y \neq x$), then $I - f$ is demiclosed at $0$ (see Lami Dozo \cite{lamidozo}, Yanagi \cite{yangai1980} and related references therein) provided $f: C: \rightarrow K(C)$ is non-expansive (here $K(C)$ denotes the family of nonempty compact subsets of $C$). We know that all Hilbert spaces and $L^p$ spaces $ p \in (1, \infty)$ have Opial's property, but it seems
that whether $I-f$ is demiclosed at zero $0$ if $f$ is a nonexpansive set-valued mapping defined on the space $X$ which is uniformly convex (e.g., $L[0, 1]$,  $1 < p < \infty$, $\neq 2$) and $f: C \rightarrow K(C)$ is nonexpansive. Here we remark that for a single-valued nonexpansive mapping $f$ is yes, which is the famous theorem of Browder \cite{bro1965}.
A remarkable fixed point theorem for multi-valued mappings is Lim's result in \cite{lim1974} which says that: If $C$ is a nonempty closed bounded convex subset of a uniformly convex Banach
space $X$ and $f: C \rightarrow K(C)$ is nonexpansive, then $f$  has a fixed point.

\vskip.1in
Now based on the concept for the semiclosed 1-set contractive mappings, we give the existence results for their best approximation, fixed points and related nonlinear alterative under the framework of $p$-seminorm spaces for $p \in (0, 1]$.

\vskip.1in
\noindent
{\bf Theorem 8.1 (Schauder Fixed Point Theorem for semiclosed 1-set contractive mappings).}
Let $U$ be  a non-empty bounded open subset of a (Hausdorf) locally $p$-convex space $E$ and its zero $0 \in U$, and $C \subset E$  be a closed $p$-convex subset of $E$ such that $0 \in C$, with $0 < p \leq 1$.
If $F: C \cap \overline{U} \rightarrow {C \cap \overline{U}}$ is continuous and semiclosed 1-set contractive. Then $T$ has at least one fixed point in $C \cap \overline{U}$.

\noindent
{\bf Proof.} As the mapping $T$ is 1-set contractive, taking an increasing sequence $\{\lambda_n\}$ such that $0 < \lambda_n < 1$ and  $\lim_{n \rightarrow \infty} \lambda_n =1$, where $n \in \mathbb{N}$.
Now we define a mapping $F_n: C \rightarrow C$ by
$F_n(x): = \lambda_n F(x)$ for each $x \in  C$ and $n\in \mathbb{N}$. Then it follows that $F_n$ is a $\lambda_n$-set-contractive mapping with $ 0 < \lambda_n < 1$. By Theorem 4.5 on the condensing mapping $F_n$ in $p$-vector space with $p$-seminorm $P_U$ for each $n \in \mathbb{N}$, there exists $x_n \in C $ such that $x_n \in F_n(x_n)=\lambda_n F(x_n)$.
Thus we have $x_n=\lambda_n F(x_n)$.
Let $P_U$ is the Minkowski $p$-functional of $U$ in $E$, it follows that $P_U$ is continuous as $0 \in int(U)=U$.
Note that for each $n \in \mathbb{N}$, $\lambda_n x_n \in \overline{U} \cap C$, which imply that
$x_n = r(\lambda_n F(x_n)) = \lambda_n F(x_n)$, thus $P_U(\lambda_n F(x_n)) \leq 1$ by Lemma 2.2.
Note that
$$P_U(F(x_n)- x_n)=P_U(F(x_n)- \lambda_n F(x_n)) =P_U(\frac{(1-\lambda_n) \lambda_n F(x_n)}{\lambda_n})
\leq (\frac{1-\lambda_n}{\lambda_n})^p P_U(\lambda_n F(x_n)) \leq (\frac{1-\lambda_n}{\lambda_n})^p,$$
which implies that
$\lim_{n\rightarrow \infty} P_U(F(x_n)-x_n)=0$. Now by the assumption that $F$ is semiclosed, which means that $(I-F)$ is closed at zero, thus there exists one point $x_0 \in \overline{C}$ such that $0 \in (I-F)(\overline{C})$, thus we have that $x_0=F(x_0)$.

Indeed, without loss of the generality, we assume that $\lim_{n \rightarrow \infty} x_n = x_0$, with $x_n=\lambda_n F(x_n)$,
and $\lim_{n \rightarrow \infty} \lambda_n=1$, it implies that
$x_0=\lim_{n \rightarrow \infty} (\lambda_n F(x_n))$, which means  $F(x_0):=\lim_{n\rightarrow \infty} F(x_n)= x_0$, thus $x_0= F(x_0)$.
We complete the proof. $\square$


\vskip.1in
\noindent
{\bf Theorem 8.2 (Best approximation for semiclosed 1-set contractive mappings).}
Let $U$ be a bounded open $p$-convex subset of a locally $p$-convex space $E$ ($0 \leq p \leq 1)$  the zero $0 \in U$, and $C$ a (bounded) closed $p$-convex subset of $E$ with also zero $0\in C$. Assume $F: \overline{U}\cap C \rightarrow C$ is a
is semiclosed 1-set contractive and continuous mapping, and for each $x \in \partial_C U$ with $F(x) \notin \overline{U}$,
$(P^{\frac{1}{p}}_U(F(x))- 1)^p \leq P_U (F(x)-x)$ for $0< p \leq 1$ (this is trivial when $p=1$). Then we have that there exist
 $x_0 \in C \cap \overline{U}$ and $F(x_0)$  such that
$ P_U (F(x_0) - x_0) = d_P(F(x_0), \overline{U}\cap C) =  d_p(F(x_0), \overline{I^p_{\overline{U}}( x_0)} \cap C)$, where $P_U$ is the Minkowski $p$-functional of $U$. More precisely, we have the following either (I) or (II) holding:

(I)  $F$ has a fixed point $x_0 \in U \cap C$, i.e., $x_0=F(x_0)$ (so that $0=P_U (F(x_0) - x_0) = d_P(F(x_0), \overline{U}\cap C) =  d_p(F(x_0), \overline{I^p_{\overline{U}}( x_0)} \cap C)$);

(II)  there exists $x_0 \in \partial_C(U)$  and $F(x_0) \notin  \overline{U}$ with
$$ P_U (F(x_0) - x_0) = d_P(F(x_0), \overline{U}\cap C) = d_p(F(x_0), \overline{I^p_{\overline{U}}( x_0)} \cap C)=(P^{\frac{1}{p}}_U(F(x_0))-1)^p  > 0.$$

\noindent
{\bf Proof.} Let $r: E \rightarrow U$ be a retraction mapping defined by $r(x): = \frac{x}{\max\{ 1, (P_U(x))^{\frac{1}{p}}\}}$
for each $x \in E$, where $P_U$ is the Minkowski $p$-functional of $U$.
Since the space $E$'s zero $0 \in U(=intU$ as $U$ is open), it follows that $r$ is continuous by Lemma 2.2.
As the mapping $F$ is 1-set contractive, taking an increasing sequence $\{\lambda_n\}$ such that $0 < \lambda_n < 1$ and  $\lim_{n \rightarrow \infty} \lambda_n =1$, where $n \in \mathbb{N}$.
Now we define a mapping $F_n: C \cap \overline{U} \rightarrow C$ by
$F_n(x): = \lambda_n F \circ r(x)$ for each $x \in  C \cap \overline{U}$ and $n\in \mathbb{N}$. Then it follows that $F_n$ is a $\lambda_n$-set-contractive mapping with $ 0 < \lambda_n < 1$  for each $n \in \mathbb{N}$.
As $C$ and $\overline{U}$ are $p$-convex, we have $r(C) \subset C$ and $r(\overline{U}) \subset \overline{U}$, so $r( C \cap \overline{U}) \subset C \cap \overline{U}$. thus $F_n$ is a self-mapping defined on $C \cap \overline{U}$.
By Theorem 4.5 for condensing mapping $F_n$, for each $n \in \mathbb{N}$, there exists $z_n \in C \cap \overline{U}$ such that $z_n \in F_n(z_n)=\lambda_n F \circ r(z_n)$. Let $x_n= r(z_n)$, then we have $x_n \in  C\cap \overline{U}$
with $x_n = r(\lambda_n F(x_n))$ such that the following (1) or (2) holding for each $n \in \mathbb{N}$:

 (1): $\lambda_n F(x_n) \in C\cap \overline{U}$;  or (2): $\lambda_n F(x_n) \in C \diagdown \overline{U}$.

\noindent
Now we prove the conclusion by considering the following two cases:

Case (I): For each $n \in N$, $\lambda_n F(x_n) \in C \cap \overline{U}$; or

Case (II): There there exists a positive integer $n$ such that $\lambda_n F(x_n) \in C \diagdown \overline{U}$.

First, by the case (I), for each $n \in \mathbb{N}$, $\lambda_n F(x_n) \in \overline{U} \cap C$, which imply that
$x_n = r(\lambda_n F(x_n)) = \lambda_n F(x_n)$, thus $P_U(\lambda_n F(x_n)) \leq 1$ by Lemma 2.2.
Note that
$$P_U(F(x_n)- x_n)=P_U(F(x_n)- \lambda_n F(x_n)) =P_U(\frac{(1-\lambda_n) \lambda_n F(x_n)}{\lambda_n})
\leq (\frac{1-\lambda_n}{\lambda_n})^p P_U(\lambda_n F(x_n)) \leq (\frac{1-\lambda_n}{\lambda_n})^p,$$
which implies that $\lim_{n\rightarrow \infty} P_U(F(x_n)-x_n)=0$. Now by the facet that $F$ is semiclosed, it implies that there exists a point $x_0 \in \overline{U}$ (i.e., the consequence $\{x_n\}_{n \in \mathbb{N}}$ has a convergent subsequence with the limit $x_0$) such that  $x_0 = F(x_0)$.  Indeed, without the loss of the generality, we assume that $\lim_{n \rightarrow \infty} x_n=x_0$, with  $x_n=\lambda_n F(x_n)$, and $\lim_{n \rightarrow \infty} \lambda_n=1$, and as $x_0=\lim_{n \rightarrow \infty} (\lambda_n F(x_n))$, which implies that $F(x_0)=\lim_{n\rightarrow \infty} F(x_n)= x_0$. Thus there exists $F(x_0)= x_0$, we have
$0 = d_p(x_0, F(x_0)) = d(F(x_0), \overline{U}\cap C) = d_p(F(x_0), \overline{I^p_{\overline{U}}(x_0)} \cap C))$ as indeed $x_0 =F(x_0) \in \overline{U}\cap C \subset \overline{I^p_{\overline{U}}( x_0)} \cap C)$.

Second, by the case (II) there exists a positive integer $n$ such that $\lambda_n F(x_n) \in C \diagdown \overline{U}$.
Then we have that $P_U(\lambda_n F(x_n))> 1$, and also $P_U(F(x_n))> 1$ as $\lambda_n < 1$.
As $ x_n = r(\lambda_n F(x_n)) = \frac{\lambda_n F(x_n)}{(P_U(\lambda_n F(x_n)))^{\frac{1}{p}}}$, which implies that $P_U(x_n)=1$, thus $x_n \in \partial_C(U)$.
Note that
$$P_U(F(x_n) - x_n)=P_U(\frac{(P_U(F(x_n))^{\frac{1}{p}}-1)F(x_n)}{P_U(F(x_n))^{\frac{1}{p}}})=(P^{\frac{1}{p}}_U(F(x_n))-1)^p. $$
By the assumption, we have $(P^{\frac{1}{p}}_U(F(x_n))-1)^p \leq P_U(F(x_n) -x)$ for $x \in C \cap \partial \overline{U}$,
it follows that
$$P_U(F(x_n))-1 \leq P_U(F(x_n)) - \sup\{P_U(z): z \in C\cap \overline{U}\} \leq \inf\{P_U(F(x_n)- z): z \in C \cap \overline{U}\}= d_p(F(x_n), C \cap \overline{U}).$$
Thus  we have the best approximation: $P_U(F(x_n) - x_n)=d_P(F(x_n), \overline{U} \cap C) = (P^{\frac{1}{p}}_U(F(x_n))-1)^p  > 0.$

Now we want to show that $P_U(F(x_n)-x_n)= d_P(F(x_n), \overline{U} \cap C) = d_p(F(x_n), \overline{I^p_{\overline{U}}( x_0)} \cap C) > 0.$

By the fact that
$(\overline{U}\cap C) \subset I^p_{\overline{U}}(x_n)\cap C$, let $z \in I^p_{\overline{U}}(x_n)\cap C \diagdown(\overline{U}\cap C)$, we first claim that
 $P_U(F(x_n) - x_n) \leq P_U(F(x_n)-z)$. If not, we have $P_U(F(x_n) - x_n) > P_U(F(x_n)-z)$.
As $z \in I^p_{\overline{U}}(x_n) \cap C \diagdown (\overline{U} \cap C)$, there exists $y \in \overline{U}$ and a non-negative number $c$ (actually $c\geq 1$ as shown soon below) with $z = x_n + c (y - x_n)$. Since $z \in C$, but $z \notin \overline{U} \cap C$, it implies that $z \notin \overline{U}$. By the fact that $x_n\in \overline{U}$ and $y \in \overline{U}$, we must have
the constant $c \geq 1$; otherwise, it implies that $z ( = (1- c )x_n + c y) \in \overline{U}$, this is impossible by our assumption, i.e., $z\notin \overline{U}$. Thus we have that $c\geq 1$, which implies that $y =\frac{1}{c} z + (1-\frac{1}{c}) x_n \in C$ (as both $x_n \in C$ and $z\in C$). On the other hand, as $z \in I^p_{\overline{U}}(x_n) \cap C \diagdown (\overline{U} \cap C)$, and $c\geq 1$ with
 $(\frac{1}{c})^p+ (1-\frac{1}{c})^p = 1 $, combing with our assumption that for each $x \in \partial_C \overline{U}$ and
 $F(x)\notin \overline{U}$,
 $P^{\frac{1}{p}}_U(F(x))- 1 \leq P^{\frac{1}{p}}_U (F(x)-x)$ for $0< p \leq 1$, it then follows that
$$P_U(F(x_n)- y) = P_U[\frac{1}{c}(F(x_n)- z)+(1-\frac{1}{c})(F(x_n) - x_n)] \leq
[(\frac{1}{c})^{p} P_U(F(x_n) -z)+(1-\frac{1}{c})^p P_U(F(x_n)- x_n)] < P_U(F(x_n)- x_n),$$
which contradicts that $P_U (F(x_n) - x_n) = d_P(F(x_n), \overline{U}\cap C)$ as shown above we know that $y \in \overline{U}\cap C$, we should have $P_U(F(x_n)- x_n)\leq P_U(F(x_n) - y)$! This helps us to complete the claim:
$P_U(F(x_n) - x_n) \leq P_U(F(x_n) - z)$ for any $z \in I^p_{\overline{U}}(x_n)\cap C \diagdown(\overline{U}\cap C)$, which means that the following best approximation of Fan's type (see \cite{fan1969}-\cite{fan1972}) holding:
$$ 0 < d_P(F(x_n), \overline{U}\cap C) = P_U (F(x_n) - x_n) = d_p(F(x_n), I^p_{\overline{U}}(x_n) \cap C).$$
Now by the continuity of $P_U$, it follows that the following best approximation of Fan type is also true:
$$ 0 <  P_U(F(x_n) - x_n) = d_P(F(x_n), \overline{U}\cap C) = d_p(F(x_n), I^p_{\overline{U}}(x_n) \cap C)
= d_p(F(x_n), \overline{I^p_{\overline{U}}(x_n)} \cap C),$$
and we have that
$$  P_U(F(x_0) - x_0) = d_P(F(x_0), \overline{U}\cap C) = d_p(F(x_0), I^p_{\overline{U}}(x_0) \cap C >0.$$
The proof is complete. $\square$

\vskip.1in
For a $p$-vector space when $p=1$, we have the following best approximation for in LCS.

\vskip.1in
\noindent
{\bf Theorem 8.3 (Best approximation for LCS)).}
Let $U$ be a bounded open convex subset of a locally convex space $E$ (i.e., $p=1$)  with zero $0 \in intU=U$ (the interior $intU=U$ as $U$ is open), and $C$ a closed $p$-convex subset of $E$ with also zero $0 \in C$. Assume that $F: \overline{U}\cap C \rightarrow C$ is a semiclosed 1-set-contractive continuous mapping. Then there exist $x_0 \in \overline{U} \cap X$ such that
$ P_U (F(x_0) - x_0) = d_P(F(x_0), \overline{U}\cap C) =  d_p(F(x_0), \overline{I_{\overline{U}}( x_0)} \cap C),$
where $P_U$ is the Minkowski $p$-functional of $U$. More precisely, we have the following either (I) or (II) holding:

(I)  $F$ has a fixed point $x_0 \in U \cap C$, i.e., $x_0=F(x_0)$ (so that $P_U (F(x_0) - x_0) = d_P(F(x_0), \overline{U}\cap C)
=  d_p(F(x_0), \overline{I_{\overline{U}}( x_0)} \cap C))=0$);

(II)  there exists $x_0 \in \partial_C(U)$ and $F(x_0) \notin \overline{U}$ with
$$ P_U (F(x_0) - x_0) = d_P(F(x_0), \overline{U}\cap C) = d_p(F(x_0), I_{\overline{U}}(x_0)\cap C)
=d_p(F(x_0), \overline{I_{\overline{U}}( x_0)} \cap C) > 0.$$

\noindent
{\bf Proof.} By applying Theorem 5.2 and the same argument used by Theorem 8.2, the conclusion follows. This completes the proof. $\square$

Now by the application of Theorem 8.2 and Theorem 8.3, we have the
the following general principle for the existence of solutions for Birkhoff-Kellogg Problems in $p$-seminorm spaces, where $(0 < p \leq 1)$.

\vskip.1in
\noindent
{\bf Theorem 8.4 (Principle of Birkhoff-Kellogg alternative).}
Let $U$ be a bounded open $p$-convex subset of a locally $p$-convex space $E$ ($0 \leq p \leq 1)$  with zero $0 \in intU=(U)$ (the interior $intU$ as $U$ is open), and $C$ a closed $p$-convex subset of $E$ with also zero $0\in C$. Assume that $F: \overline{U}\cap C \rightarrow C$ is a semiclosed 1-set-contractive continuous mapping. Then $F$ has at least one of the following two properties:

(I) $F$ has a fixed point $x_0 \in U \cap C$ such that $x_0 = F(x_0)$, or

(II) there exist $x_0 \in \partial_C(U)$ and $F(x_0)\notin \overline{U}$,  and $\lambda = \frac{1}{(P_U(F(x_0)))^{\frac{1}{p}}} \in (0, 1)$ such that
$x_0 = \lambda F(x_0)$. In addition if for each $x \in \partial_C U$, $P^{\frac{1}{p}}_U(F(x))- 1 \leq P^{\frac{1}{p}}_U (F(x)-x)$ for $0< p \leq 1$ (this is trivial when $p=1$), then the best approximation between $x_0$ and $F(x_0)$ given by $$ P_U (F(x_0) - x_0) = d_P(F(x_0), \overline{U}\cap C) = d_p(F(x_0), \overline{I^p_{\overline{U}}(x_0)} \cap C) = (P^{\frac{1}{p}}_U(F(x_0))-1)^p > 0.$$

\vskip.1in
{\bf Proof.} If (I) is  not the case, then (II) is proved by the Remark 5.2 and by following the proof in Theorem 8.2 for the case (ii): $F(x_0)\in C \diagdown \overline{U}$,  with $y_0= f(F(x_0))$, where $f$ is the restriction of the continuous mapping $r$ restriction to the subset $U$ in $E$. Indeed, as $y_0 \notin \overline{U}$, it follows that $P_U(y_0) > 1$, and $x_0= f(y_0) = F(x_0) \frac{1}{(P_U(F(x_0)))^{\frac{1}{p}}}$.
Now let $\lambda = \frac{1}{(P_U(F(x_0)))^{\frac{1}{p}}}$, we have $\lambda < 1$ and  $x_0 = \lambda F(x_0)$. Finally, the additionally assumption in (II) allows us to have the  best approximation between $x_0$ and $F(x_0)$ obtained by following the proof of Theorem 8.2 as $P_U (F(x_0) - x_0) = d_P(F(x_0), \overline{U}\cap C) = d_p(F(x_0), \overline{I^p_{\overline{U}}(x_0)} \cap C) > 0$. This completes the proof. $\square$

\vskip.1in
As an application of Theorem 8.2 for the non-self mappings,  we have the following general principle of Birkhoff-Kellogg alternative in TVS.

\vskip.1in
\noindent
{\bf Theorem 8.5 (Principle of Birkhoff-Kellogg alternative in LCS).}
Let $U$ be a bounded open $p$-convex subset of the LCS $E$  with the zero $0 \in U$, and $C$ a closed convex subset of $E$ with also zero $0\in C$. Assume the
$F: \overline{U}\cap C \rightarrow C$ is a semiclosed 1-set contractive and continuous mapping. Then it has at least one of the following two properties:

(I) $F$ has a fixed point $x_0 \in U \cap C$ such that $x_0 = F(x_0)$; or

(II) there exists $x_0 \in \partial_C(U)$ and $F(x_0) \notin \overline{U}$  and $\lambda \in (0, 1)$ such that $x_0 = \lambda F(x_0)$, and the best approximation between $\{x_0\}$ and $F(x_0)$ is given by
$ P_U (F(x_0) - x_0) = d_P(F(x_0), \overline{U}\cap C) = d_p(F(x_0), \overline{I^p_{\overline{U}}( x_0)} \cap C) > 0.$


\vskip.1in
On the other hand, by the Proof of Theorems 8.2, we note that for case (II) of Theorem 8.2, the assumption $``$each $x \in \partial_C U$ with $y \in F(x)$, $P^{\frac{1}{p}}_U(y)- 1 \leq P^{\frac{1}{p}}_U (y-x)$" is only used to guarantee the best approximation $``P_U (y_0 - x_0) = d_P(y_0, \overline{U}\cap C) = d_p(y_0, \overline{I^p_{\overline{U}}( x_0)} \cap C) > 0$", thus
we have the following Leray-Schauder alternative in $p$-vector spaces, which, of course, includes the corresponding results in locally convex spaces as special cases.

\vskip.1in
\noindent
{\bf Theorem 8.6 (The Leray-Schauder Nonlinear Alternative).} Let $C$ a closed $p$-convex subset of $P$-seminorm space $E$ with $0 \leq p \leq 1$  and the zero $0 \in C$. Assume the $F: C \rightarrow C$ is a semiclosed 1-set contractive and continuous mapping. Let $\varepsilon(F): =\{x \in C: x\in \lambda F(x), \mbox{ for some } 0 < \lambda < 1\}$. Then either $F$ has a fixed point in $C$ or the set $\varepsilon(F)$ is unbounded.

\noindent
{\bf Proof.} By assuming the case (I) is not true, i.e., $F$ has no fixed point, then we claim that the set $\varepsilon(F)$ is unbounded. Otherwise, assume the set $\varepsilon(F)$ is bounded. and assume $P$ is the continuous $p$-seminorm for $E$, then there exists $r>0$ such that the set $B(0, r):=\{x \in E: P(x) < r\}$ , which contains the set $\varepsilon(F)$, i.e., $\varepsilon(F) \subset B(0, r)$, which means for any $x \in \varepsilon(F)$, $P(x) < r$. Then $B(0. r)$ is an open $p$-convex subset of $E$ and the zero $0 \in B(0, r)$ by Lemma 2.2 and Remark 2.4.
Now let $U: =B(0, r)$ in Theorem 8.4, it follows that for the mapping $F: B(0, r) \cap C \rightarrow 2^C$ satisfies all general conditions of Theorem 8.4, and we have that any $x_0 \in \partial_C B(0, r)$, no any $\lambda \in (0, 1)$ such that
$x_0=\lambda y_0$, where $y_0 \in F(x_0)$. Indeed, for any $x \in \varepsilon(F)$, it follows that $P(x) < r$ as $\varepsilon(F) \subset B(0, r)$, but for any $x_0 \in \partial_C B(0, r)$, we have $P(x_0)=r$, thus the conclusion (II) of Theorem 8.4 does not have hold. By Theorem 8.4 again, $F$ must have a fixed point, but this contradicts with our assumption that $F$ is fixed point free. This completes the proof. $\square$

\vskip.1in
Now assume a given $p$-vector space $E$ equipped with the $P$-seminorm (by assuming it is continuous at zero) for $0< p \leq 1$, then we know that $P: E \rightarrow \mathbb{R}^+$, $P^{-1}(0)=0$,  $P(\lambda x) = |\lambda|^p P(x)$ for any  $x\in E$ and
$\lambda \in \mathbb{R}$. Then we have the following useful result for fixed points due to Rothe and Altman types in $p$-vector spaces, which plays important roles for optimization problem, variational inequality, complementarity problems.

\vskip.1in
\noindent
{\bf Corollary 8.1.} Let $U$ be a bounded open $p$-convex subset of a locally $p$-convex space $E$ and zero $0 \in U$, plus $C$ is a closed $p$-convex subset of $E$ with $U \subset C$, where $0< p \leq 1$. Assume that $F: \overline{U} \rightarrow C$ is a semiclosed 1-set contractive continuous mapping. if one of the following is satisfied,

(1) (Rothe type condition): $P_U(F(X)) \leq P_U(x)$ for any $x \in \partial U$;

(2) (Petryshyn type condition): $P_U(F(X)) \leq P_U(F(X)-x)$ for any $x \in \partial U$;

(3) (Altman type condition): $|P_U(F(X))|^{\frac{2}{p}} \leq [P_U(F(X))- x)]^{\frac{2}{p}} + [P_U(x)]^{\frac{2}{p}}$ for any $x \in \partial U$;

then $F$ has at least one fixed point.

\noindent
{\bf Proof.} By the conditions (1), (2) and (3), it follows that the conclusion of (II) in Theorem 8.4 $``$there exist $x_0 \in \partial_C(U)$ and $\lambda \in (0, 1)$ such that $x_0 \neq F(x_0)$" does not hold, thus by the alternative of Theorem 8.4, $F$ has a fixed point. This completes the proof. $\square$.

\vskip.1in
By the fact that when  $p=1$ in $p$-vector space being a LCS, we have the following classical Fan's best approximation (see \cite{fan1969}) as a powerful tool for the study in the optimization, mathematical programming, games theory, and mathematical economics, and others related topics in applied mathematics.

\vskip.1in
\noindent
{\bf Corollary 8.2 (Fan's best approximation).} Let $U$ be a bounded open convex subset of a locally convex space $E$ with
the zero $0 \in U$, and $C$ a closed convex subset of $E$ with also zero $0\in C$, and assume $F: \overline{U}\cap C \rightarrow C$ is a semiclosed 1-set contractive and continuous mapping.
Then there exist $x_0 \in \overline{U} \cap X$ such that $ P_U (F(x_0) - x_0) = d_P(F(x_0, \overline{U}\cap C) =
 d_p(F(x_0), \overline{I_{\overline{U}}( x_0)} \cap C)$,
where $P_U$ being the Minkowski $p$-functional of $U$ in $E$. More precisely, we have the following either (I) or (II) holding, where $W_{\overline{U}}(x_0)$ is either inward set $I_{\overline{U}}(x_0)$, or the outward set $O_{\overline{U}}(x_0)$:

(I)  $F$ has a fixed point $x_0 \in U \cap C$, i.e., $x_0=F(x_0)$;

(II)  there exists $x_0 \in \partial_C(U)$  with $ F(x_0) \notin \overline{U}$ such that
$$ P_U (F(x_0) - x_0) = d_P(F(x_0), \overline{U}\cap C) = d_p(F(x_0), \overline{i_{\overline{U}}( x_0)} \cap C) = P_U(F(x_0)) - 1 > 0.$$

\noindent
{\bf Proof.} When $p=1$, then it automatically satisfies that the inequality:
$ P^{\frac{1}{p}}_U(F(x))- 1 \leq P^{\frac{1}{p}}_U (F(x_)-x)$ for each $x \in \overline{U}\cap C$.
Indeed we have that for $x_0 \in \partial_C(U)$, we have $P_U (F(x_0) - x_0) = d_P(F(x_0), \overline{U}\cap C) =
d_p(F(x_0), \overline{I_{\overline{U}}(x_0)} \cap C)=  P_U(F(x_0))-1$. The conclusions are given by Theorem 8.2 (or Theorem 8.3). The proof is complete. $\square$

\vskip.1in
We like to point out the similar results on Rothe and Leray-Schauder alternative have been developed by Isac \cite{isac}, Park \cite{park1995}, Potter \cite{potter1972}, Shahzad \cite{shahzad2006}-\cite{shahzad2004}, Xiao and Zhu \cite{xiaozhu2011}, and related references therein as tools of nonlinear analysis in topological vector spaces. As mentioned above, when $p=1$ and $F$ as a continuous mapping, then we can obtain the version of Lerary-Schauder in locally convex spaces, and we omit their statements in details here due to the limit of the space.

\vskip.1in
\section{Nonlinear Alternatives Principle for the Class of Semiclosed 1-Set Contractive Mappings}

\vskip.1in
As applications of results in Section 8 above, we new establish general results for the existence of solutions for Birkhoff-Kellogg problem, and the principle of Leray-Schauder alternatives for semiclosed 1-set contractive mappings for $p$-vector spaces being locally $p$-convex spaces for $0 < p \leq 1$.

\noindent
{\bf Theorem 9.1 (Birkhoff-Kellogg alternative in locally $p$-convex spaces).} Let $U$ be a bounded open  $p$-convex subset of a locally $p$-convex space $E$ (where, $0 \leq p \leq 1)$ with the zero $0 \in U$, and $C$ a closed $p$-convex subset of $E$ with also zero $0\in C$, and assume $F: \overline{U}\cap C \rightarrow 2^C$ is a semiclosed 1-set contractive and continuous mapping, and for each $x \in \partial_C(U)$ with
$P^{\frac{1}{p}}_U(F(x))- 1 \leq P^{\frac{1}{p}}_U (F(x)-x)$ for $0< p \leq 1$ (this is trivial when $p=1$), where $P_U$ is the Minkowski $p$-functional of $U$.
Then we have that either (I) or (II) holding below:

(I) there exists $x_0 \in  \overline{U}\cap C$ such that $x_0 = F(x_0)$;

(II) there exists $x_0 \in \partial_C(U)$ with $F(x_0) \notin \overline{U}$ and $\lambda >1$ such that
$\lambda x_0 = F(x_0)$, i.e., $F(x_0) \in \{\lambda x_0: \lambda > 1 \}$.

\noindent
{\bf Proof.} By following the argument and notations used by Theorem 8.2, we have that either

(1) $F$ has a fixed point $x_0 \in U \cap C$; or

(2) there exists $x_0 \in \partial_C(U)$ with $x_0=f(y_0)$ such that
$$ P_U (F(x_0) - x_0) = d_P(F(x_0), \overline{U}\cap C) = d_p(F(x_0), \overline{I_{\overline{U}}( x_0)} \cap C) = P_U(F(x_0)) - 1 > 0,$$
where $\partial_C(U)$ denotes the boundary of $U$ relative to $C$ in $E$, and $f$ is the restriction of the continuous retraction $r$ respect to the set $U$ in $E$.

If $F$ has no fixed point, then above (2) holds and $x_0 \neq F(x_0)$.
As given by the proof of Theorem 8.2, we have that $F(x_0) \notin \overline{U}$, thus $P_U(F(x_0)) > 1$ and
$x_0= f(y_0)=\frac{F(x_0)}{(P_U(F(x_0)))^{\frac{1}{p}}}$, which means $F(x_0) =(P_U(F(x_0)))^{\frac{1}{p}} x_0$.
Let $\lambda = (P_U(F(x_0)))^{\frac{1}{p}}$, then $\lambda > 1$ and we have $ \lambda x_0 = F(x_0)$.
This completes the proof. $\square$.

\vskip.1in
\noindent
{\bf Theorem 9.2 (Birkhoff-Kellogg alternative in LCS).} Let $U$ be a bounded open convex subset of a locally $p$-convex space $E$ with
the zero $0 \in U$, and $C$ a closed convex subset of $E$ with also zero $0\in C$, and assume $F: \overline{U}\cap C \rightarrow C$ is a semiclosed 1-set contractive and continuous mapping.
Then we have the following either (I) or (II) holding: 

(I) there exists $x_0 \in  \overline{U}\cap C$ such that $x_0 = F(x_0)$; or

(II) there exists $x_0 \in \partial_C(U)$ with $F(x_0) \notin \overline{U}$ and $\lambda >1$ such that
$\lambda x_0 = F(x_0)$, i.e., $F(x_0) \in \{\lambda x_0: \lambda > 1 \}.$

\noindent
{\bf Proof.} When $p=1$, then it automatically satisfies that the inequality:
$ P^{\frac{1}{p}}_U(F(x))- 1 \leq P^{\frac{1}{p}}_U (F(x)-x)$ for all $x \in \overline{U}\cap C$. Indeed we have that for $x_0\in \partial_C(U)$,
we have $P_U (F(x_0) - x_0) = d_P(F(x_0), \overline{U}\cap C) = d_p(F(x_0), \overline{I_{\overline{U}}( x_0)} \cap C)=  P_U(F(x_0))-1$. The conclusions are given by Theorems 8.3 and 8.4. The proof is complete. $\square$

\vskip.1in
Indeed, we have the following fixed points for non-self mappings in $p$-vector spaces for $0 < p \leq 1$ under different boundary conditions.

\vskip.1in
\noindent
{\bf Theorem 9.3 (Fixed Points of non-self mappings).} Let $U$ be a bounded open $p$-convex subset of a locally $p$-convex space $E$ (where, $0 \leq p \leq 1)$ with the zero $0 \in U$, and $C$ a closed $p$-convex subset of $E$ with also zero $0\in C$, and assume
$F: \overline{U}\cap C \rightarrow C$ is a semiclosed 1-set contractive and continuous mapping. In addition, for each $x \in \partial_C(U)$,
$P^{\frac{1}{p}}_U(F(x))- 1 \leq P^{\frac{1}{p}}_U (F(x)-x)$ for $0 < p \leq 1$ (this is trivial when $p = 1 $), where $P_U$ is the Minkowski $p$-functional of $U$.
If $F$ satisfies any one of the following conditions for any $x \in \partial_C(U) \diagdown F(x)$:

(i)  $P_U(F(x)-z) < P_U(F(x)-x)$ for some $z \in \overline{I_{\overline{U}}(x)}\cap C$;

(ii) there exists $\lambda $ with $|\lambda| < 1$ such that $\lambda x + (1-\lambda)F(x)
 \in \overline{I_{\overline{U}}(x)}\cap C$;

(iii) $F(x) \in \overline{I_{\overline{U}}(x)}\cap C$;

(iv) $F(x) \in  \{\lambda x: \lambda > 1 \} =\emptyset$;

(v) $F(\partial U) \subset \overline{U} \cap C$;

(vi) $P_U(F(x)-x) \neq ((P_U(F(x)))^{\frac{1}{p}}-1)^p$;


then $F$ must has a fixed point.

\noindent
{\bf Proof.} By following the argument and symbols used in the proof of Theorem 8.2 (see also Theorem 8.4), we have that either

(1) $F$ has a fixed point $x_0 \in U \cap C$; or

(2) there exists $x_0 \in \partial_C(U)$ with $x_0=f(F(x_0))$ such that
$$ P_U (F(x_0) - x_0) = d_P(F(x_0), \overline{U}\cap C) = d_p(F(x_0), \overline{I_{\overline{U}}( x_0)} \cap C) = P_U(F(x_0)) - 1 > 0,$$
where $\partial_C(U)$ denotes the boundary of $U$ relative to $C$ in $E$, and $f$ is the restriction of the continuous retraction $r$ respect to the set $U$ in $E$.

First, suppose that $F$ satisfies the condition (i), if $F$ has no fixed point, then above (2) holds and $x_0 \neq  F(x_0)$. Then by the condition (i), it follows that $P_U(F(x_0)-z) < P_U(F(x_0)-x_0)$ for some $z \in \overline{I_{\overline{U}}(x)}\cap C$, this contradicts with the best approximation equations given by (2) above, thus $F$ mush have a fixed pint.

Second, suppose that $F$ satisfies the condition (ii), if $F$ has no fixed point, then above (2) holds and $x_0 \neq F(x_0)$. Then by
the condition (ii), there exists $\lambda >1$ such that $\lambda  x_0 + (1 - \lambda ) F(x_0) \in
\overline{I_{\overline{U}}(x)}\cap C$. It follows that
$$P_U(F(x_0)- x_0) \leq P_U(F(x_0)- (\lambda x_0 + (1-\lambda F(x_0)))=P_U(\lambda (F(x_0) - x_0))=|\lambda|^p P_U(F(x_0)-x_0)<P_U(F(x_0)-x_0)$$
this is impossible and thus $F$ must have a fixed point in $\overline{U}\cap C$.

Third, suppose that $F$ satisfies the condition (iii), i.e., $F(x) \in \overline{I_{\overline{U}}(x)} \cap C$;, then the (2), we have that $P_U (F(x_0) - x_0)$ and thus $x_0= F(x_0)$, which means $F$ has a fixed point.

Forth, suppose that $F$ satisfies the condition (iv),  and  if $F$ has no fixed point, then above (2) holds and $x_0 \neq F(x_0)$.
As given by the proof of Theorem 8.2, we have that $F(x_0) \notin \overline{U}$, thus $P_U(F(x_0)) > 1$ and
$x_0= f(F(x_0))=\frac{F(x_0)}{(P_U(F(x_0)))^{\frac{1}{p}}}$, which means $F(x_0)=(P_U(F(x_0)))^{\frac{1}{p}} x_0$,
where $(P_U(F(x_0)))^{\frac{1}{p}} > 1$, this contradicts with the assumption (iv), thus $F$ must have a fixed point in $\overline{U} \cap C$.

Fifth, suppose that $F$ satisfies the condition (v), then $x_0 \neq  F(x_0)$.
As $x_0 \in \partial_C{U}$, now by the condition (v), we have that $F(\partial U) \subset \overline{U} \cap C$, it follows that for any
we have $F(x_0) \in \overline{U}\cap C$, thus $F(x) \notin \overline{U} \diagdown \cap C$, which implies that
$0 < P_U(F(x_0)- x_0) = d_P(F(x_0), \overline{U}\cap C) = 0$, this is impossible, thus $F$ must have a fixed point.
Here, as pointed out by Remark 5.2, we know that  based on the condition (v), the mapping $F$ has a fixed point by applying $F(\partial U) \subset \overline{U} \cap C$ is enough, not needing the general hypothesis: $``$for each $x \in \partial_C(U)$,  $P^{\frac{1}{p}}_U(F(x))- 1 \leq P^{\frac{1}{p}}_U (F(x)-x)$ for $0 < p \leq 1$".

Finally, suppose that $F$ satisfies the condition (vi), if $F$ has no fixed point, then above (2) holds and $x_0 \neq F(x_0)$. Then the condition (v) implies that $P_U(F(x_0)- x_0) \neq ((P_U(F(x))^{\frac{1}{p}}-1)^p$, but the our proof in Theorem 5.2 shows that
$P_U(y_0- x_0)=((P_U(y))^{\frac{1}{p}}-1)^p$, this is impossible, thus $F$ must have a fixed point.
Then the proof is complete. $\square$

\vskip.1in
Now by taking the set $C$ in Theorem 8.1 as the whole locally  $p$-convex space $E$ itself, we have the following general results for non-self upper semi-continuous mappings which include results of Rothe, Petryshyn, Altman and Leray-Schauder types' fixed points as special cases.


\vskip.1in
Taking $p=1$ and $C =E$ in Theorem 9.3, we have the following fixed points for non-self continuous mappings associated with inward or outward sets for locally convex spaces which are locally $p$-convex spaces for $p=1$.

\vskip.1in
\noindent
{\bf Theorem 9.4 (Fixed Points of non-self mappings with boundary conditions).} Let $U$ be a bounded open convex subset of the LCS  $E$ with the zero $0 \in U$, and assume $F: \overline{U} \rightarrow E$ is a semiclosed 1-set contractive and continuous mapping. If $F$ satisfies any one of the following conditions for any $x \in \partial(U) \diagdown F(x)$

(i) $P_U(F(x)-z) < P_U(F(x)-x)$ for some $z \in \overline{I_{\overline{U}}(x)}$; 

(ii) there exists $\lambda $ with $|\lambda| < 1$ such that $\lambda x + (1-\lambda)F(x)
 \in \overline{I_{\overline{U}}(x)}$; 

(iii) $F(x) \in \overline{I_{\overline{U}}(x)}$; 

(iv) $F(x) \in \{\lambda x: \lambda > 1 \} =\emptyset$;

(v) $F(\partial(U) \subset \overline{U}$;

(vi) $P_U(F(x)-x) \neq P_U(F(x))-1$;


then $F$ must has a fixed point.

\vskip.1in
In what follow, based on the best approximation theorem in $p$-seminorm space, we will also give some fixed point theorems for non-self continuous mappings with various boundary conditions which are related to the study for the existence of solutions for PDE and differential equations with boundary problems (see, Browder \cite{bro1968}, Petryshyn \cite{petryshyn1966}-\cite{petryshyn1973tams}, Reich \cite{reich}),  which would play roles in nonlinear analysis for $p$-seminorm space as shown below.

First, as discussed by Remark 5.2, the proof of Theorem 9.2, with the strongly boundary condition $``$$F(\partial(U)) \subset \overline{U} \cap C$"  only, we can prove that  $F$ has a fixed point, thus we have the following fixed point theorem of Rothe type in locally
$p$-convex spaces.

\vskip.1in
\noindent
{\bf Theorem 9.5 (Rothe Type).}
Let $U$ be a bounded open $p$-convex subset of a locally $p$-convex space $E$ (where, $ 0 \leq p \leq 1)$ with the zero $0 \in U$.
Assume $F: \overline{U}\rightarrow E$ is a semi 1-set contractive and continuous mapping, and such that $F(\partial(U)) \subset \overline{U}$, then $F$ must has a fixed point.

\vskip.1in
Now as applications of Theorem 9.5, we give the following Leray-Schauder Alternative in $p$-vector spaces for non-self set-valued mappings associated with the boundary condition  which often appear in the applications (see Isac \cite{isac} and references therein for the study of complementary problems and related topics in optimization).

By using the same argument used in the proof of Theorem 6.6, we have the
following result. 
\vskip.1in
\noindent
{\bf Theorem 9.6 (Leray-Schauder Alternative in locally $p$-convex Spaces).} Let $E$ be a locally $p$-convex space $E$, where $0 < p \leq 1$, $B \subset E$ a bounded closed $p$-convex such that $0 \in int B$. Let $F: [0, 1] \times B  \rightarrow E$  be semiclosed 1-set contractive and continuous
mapping, and such that
the set $F([0, 1] \times B)$ be relatively compact in $E$. If the following assumptions are satisfied:

(1) $x \neq F(t, x)$ for all $x \notin \partial B$ and $t \in [0, 1]$,

(2) $F(\{0\} \times \partial B) \subset B$,

\noindent
then there is an element $x^* \in B$ such that $x^* = F(1, x^*)$.

\noindent
{\bf Proof.} The conclusion is proved by following the same argument used in Theorem 6.6.
The proof is complete. $\square$

\vskip.1in
As a special case of Theorem 9.6, we have the following principle for the implicit form of Leray-Schauder type alternative in
in locally $p$-convex spaces for $0< p \leq 1.$

\noindent
{\bf Corollary 9.1 (The Implicit Leray-Schauder Alternative).} Let $E$ be  a locally $p$-convex space $E$,
where $0 < p \leq 1$, $B \subset E$ a bounded closed $p$-convex such that $ 0 \in int B$. Let $F: [0, 1] \times B \rightarrow E$  be semiclosed 1-set contractive and continuous, and the set $F([0, 1] \times B)$ be relatively compact in $E$. If the following assumptions are satisfied:

(1) $F(\{0\} \times \partial B) \subset B$,

(2) $x \notin F(0, x)$ for all $x \in \partial B$,

\noindent
then at least one of the following properties is satisfied:

(i) there exists $x^* \in  B$ such that $x^* = F(1, x^*)$; or

(ii) there exists  $(\lambda^*, x^*) \in (0, 1) \times \partial B$ such that $x^* = F(\lambda^*, x^*)$.

\noindent
{\bf Proof.}  The result is an immediate consequence of Theorem 9.6, this completes the proof. $\square$
\vskip.1in

We like to point out the similar results on Rothe and Leray-Schauder alternative have been developed by Furi and Pera \cite{furipera}, Granas and Dugundji \cite{granas}, G\'{o}rniewicz \cite{gorniewicz}, G\'{o}rniewicz et al.\cite{gorniewiczetal},
Isac \cite{isac}, Li et al.\cite{lixuduan2006}, Liu \cite{liu2001}, Park \cite{park1995}, Potter \cite{potter1972}, Shahzad \cite{shahzad2006}-\cite{shahzad2004}, Xu \cite{xu2007}, Xu et al.\cite{xujiali2006}, and related references therein as tools of nonlinear analysis in the Banach space setting and applications to the boundary value problems for ordinary differential equations in noncompact problems, a general class of mappings for nonlinear alternative of Leray-Schauder type in normal topological spaces, and some Birkhoff-Kellogg type theorems for general class mappings in topological vector spaces are also established by Agarwal et al.\cite{Agarwal}, Agarwal and O'Regan \cite{agarwalorgan2003}-\cite{agarwalorgan2004}, Park \cite{park1997}, and references therein for more in detail; and in particular, recently O'Regan \cite{oregan2019} uses the Leray-Schauder type coincidence theory to establish some Birkhoff-Kellogg problem, Furi-Pera type results for a general class of mappings.

Before closing this section,  we like to share that as the application of the best approximation result for 1-set contractive mappings, we can establish the fixed point theorems and general principle of Leray-Schauder alternative for non-self mappings, which would seem play important roles for the development of nonlinear analysis for $p$-vector spaces for $0 < p \leq 1$,  as the nature extension and  achievement of nonlinear functional analysis in mathematics for the underling being locally convex vector spaces locally convex spaces, normed spaces, or in Banach spaces.

\vskip.1in
\section{Fixed Points for the class of Semiclosed 1-Set Contractive Mappings}

In this section, based on the best approximation Theorem 8.2 established for the 1-set contractive mappings in Section 8, we will show how it is used as a useful tool for us to develop fixed point theorems for semiclosed 1-set contractive non-self upper semi-continuous mappings in $p$-seminorm spaces (for $p \in (0, 1]$, by including seminorm, norm spaces, and uniformly convex Banach spaces as special cases).

By following the Definition 7.1 above, we first observe that if $f$ is a continuous demicompact mapping, then $(I - f)$ is closed, where $I$ is the identity mapping on $X$. it is also clear from definitions that every demicompact map is hemicompact in seminorm spaces, but the converse is not true in general (e.g., see the example in p.380 by Tan and Yuan \cite{tanyuan1994}).
It is evident that if $f$ is demicompact, then $I-f$ is demiclosed. It is know that for each condensing mapping $f$, when $D$ or $f(D)$ is bounded, then $f$ is hemicompact; and also $f$ is demicompact in metric spaces by Lemma 2.1 and Lemma 2.2 of Tan and Yuan \cite{tanyuan1994}, respectively. In addtion, it is known that every nonexpansive map is a 1-set-contractive map; and also if $f$ is a hemicompact 1-set-contractive mapping, then f is a 1-set-contractive mapping satisfying the following $``${\bf Condition (H1)}" (the same as (H1), and slightly different from the condition (H) used in Section 5):

{\bf (H1) Condition}: Let $D$ be a nonempty bounded subset of a space $E$ and assume $F: \overline{D} \rightarrow 2^E$ a set-valued mapping. If $\{x_n\}_{n \in \mathbb{N}}$ is any sequence in $D$ such that for each $x_n$, there exists $y_n \in F(x_n)$ with $\lim_{n \rightarrow \infty} (x_n- y_n)=0$, then there exists a point
$x\in \overline{D}$ such that $x \in F(x)$.

\vskip.1in
We first note that the $``$(H1) Condition" above is actually the $``$Condition (C)"  used by Theorem 1 of  Petryshyn \cite{petryshyn1973tams}. Indeed, by following Goebel and Kirk \cite{goebelkirk2008} (see also Xu \cite{xu2000} and reference therein), Browder \cite{bro1968} (see also \cite{bro1976}, p.103) proved that if $K$ is a closed and convex subset of a uniformly
convex Banach space $X$, and if $T: K \rightarrow X$ is nonexpansive, then the mapping
$f: = I - T$ is demiclosed on $X$. This result, known as Browder's demiclosedness
principle (Browder's proof, which was inspired by the technique of G\"ohde in \cite{gohde}), is one of the fundamental results in the theory of nonexpansive mappings, which satisfies the $``$(H1) condition".

The following is the Browder's demiclosedness
principle proved by Browder  \cite{bro1968} that says that a nonexpansive mapping in a uniformly convex Banach $X$ enjoys the condition (H1) as shown below.

\vskip.1in
\noindent
{\bf Lemma 10.1.} Let $D$ be a nonempty bonded convex subset of a uniformly convex Banach space $E$. Assume $F: \overline{D} \rightarrow E$ is a non-expansive single-valued mapping, then the mapping $P: =I - F$  defined by $P(x): = (x-F(x)) $  for each $x \in   \overline{D}$ is demiclosed, and in particular, the $``$(H1) Condition" holds.

\noindent
{\bf Proof.} By following the argument given in p.329 (see also the proof of Theorem 2.2 and Corollary 2.1) by Petryshyn \cite{petryshyn1973tams}, by the Browder demiclosedness principle (see Goebel and Kirk \cite{goebelkirk2008} or Xu \cite{xu2000}),
$P=(I-F)$ is closed at zero, thus there exists $x_0 in \overline{U}$ such $0 \in (I-F)(x_0)$, which means that
$x_0 \in F(x_0)$. The proof is complete. $\square$.

\vskip.1in
On the other hand, by following the notion called $``$Opial's condition"  given by Opial \cite{opial1967}, which says that a Banach space $X$ is said to satisfy Opial's condition if $\liminf_{n \rightarrow \infty} \| w_n - w \| < \liminf_{n\rightarrow \infty} \|w_n-p\|$ whenever $(w_n)$ is a sequence in $X$ weakly convergent to $w$ and $p\neq w$, we know that Opial's condition plays an important role in the fixed point theory, e.g., see Lami Dozo \cite{lamidozo}, Goebel and Kirk \cite{goebelkirk2008}, Xu \cite{xu2000}
and references where. Actually, the following result shows that there exists a class of non-expansive set-valued mappings in Banach spaces with Opial's condition (see Lami Dozo \cite{lamidozo} satisfying the $``$(H1) Condition".

\vskip.1in
\noindent
{\bf Lemma 10.2}. Let $C$ is a nonempty convex weakly compact subset of a Banach space $X$ which satisfies Opial's condition. Let $T:  C \rightarrow K(C)$ be a non-expansive set-valued mapping with non-empty compact-values.
Then the graph of $(I-T)$ is closed $(X, \sigma(X, X^*) \times (X, \|\cdot\|))$,
thus $T$ satisfies the $``$(H1) condition", where, $I$ denotes the identity on $X$, $\sigma(X, X^*)$ the weak topology, and $\|\cdot\|$  the  norm (or strong) topology.

\noindent
{\bf Proof.} By following Theorem 3.1 of Lami Dozo \cite{lamidozo}, it follows that the mapping $T$ is demi-closed, thus
$T$ satisfies the $``$(H1) condition". The proof is complete. $\square$

\vskip.1in
By Theorem 3.1 of Lami Dozo \cite{lamidozo}, indeed we have the following statement which is an another version by using the term of $``$distance convergence" for Lemma 10.2.

\vskip.1in
\noindent
{\bf Lemma 10.3}.  Let $C$ be a nonempty closed convex subset of a
Banach space $(X, d)$ which satisfies the Opial condition. Let $T: C \rightarrow K(C)$ be a multi-valued nonexpansive mapping (with the fixed points).
Let $(y_n)_{n \in \mathbb{N}}$ be a bounded sequence, such that $_{n \rightarrow \infty}d(y_, T(y_n))=0$,
then the weak cluster points of $(y_n)$, $n \in \mathbb{N}$ is a fixed point of $T$.

\noindent
{\bf Proof.} It is Theorem 3.1 of Lami Dozo \cite{lamidozo} (see also Lemma 3.2 of Xu and Muglia \cite{xumuglia}).
$\square$

\vskip.1in
We note that another class of set-valued mappings, called $``*$-nonexpansive mappings in Banach spaces
(introduced by Husain and Tarafdar \cite{husaintarafdar}, see also Husain and Latif \cite{husainlatif}) which was proved to hold the demiclosedness principle in reflexive Banach spaces satisfying Opial's condition by Muglia and Marino (i.e., Lemma 3.4 in \cite{mugliamarino}, thus the demiclosedness principle also holds in reflexive Banach spaces with duality mapping that is weakly sequentially continuous since these satisfy Opial's condition.

More precisely, let $C$ be a subset of a Banach space $(X, \|\cdot\|)$, and $K(C)$ be the family of compact subsets of $C$. By
following Husain and Latif \cite{husainlatif}, a mapping $W: C \rightarrow K(C)$ is said to be $*$-nonexpansive
if for all $x, y \in C$, and $x^W \in W(x)$ such that $\|x -x^W\|=d(x, W(x))$,
there exists $y^W \in W(y)$ with $\|y- y^W\|= d(y, W(y))$ such that
$\|x^W-y^W\| \leq \|x-y\|$.

As pointed by Muglia and Marino \cite{mugliamarino}, however, $*$-nonexpansivity and multivalued nonexpansivity are not so
far, By Theorem 3 of L\'{o}pez-Acdeo and Xu \cite{lopezxu}, it is proved that a multivalued mapping $W: C \rightarrow K(C)$
is $*$-nonexpansive if and only if the metric projection $P_W(x);=\{u_x \in W(x): \|x - u_x\|=\inf_{y \in W(x)}\|x-y\|\}$
is nonexpansive.

We now have the following result which is the demiclosedness principle for multivalued $*$-nonexpansive mapping
given by Lemma 3.4 of  Muglia and Marino \cite{mugliamarino}.

\vskip.1in
\noindent
{\bf Lemma 10.4.} Let $X$ be a reflexive space satisfying Opial condition and let $W: X \rightarrow K(X)$ a $*$-nonexpansive multivalued mapping with fixed points (existing) (denoted by $Fix(W))$.
Let $(y_n)_{n \in \mathbb{N}}$ be a bounded sequence such that
$\lim_{n \rightarrow \infty}d(y_n, W(y_n)) \rightarrow 0$. Then the weak cluster points of
$(y_n)_{n \in \mathbb{N}}$ belong to $Fix(W)$.

\noindent
{\bf Proof.} It is Lemma 3.4 of  Muglia and Marino \cite{mugliamarino}. $\square$

\vskip.1in
\noindent
{\bf Remark 10.1.} We like to point out that indeed, Xu \cite{xu1991} proved existence results of fixed points for $*$-nonexpansive on strictly convex Banach spaces, and L\'{o}pezo-Acdeo and Xu in \cite{lopezxu} have obtained existence result in the setting of Banach space satisfying Opial condition, so the assumption on the existence of fixed points for the mapping $W$ in Lemma 10.4 makes sense for the setting under either strictly convex Banach spaces, or Banach space satisfying Opial condition.

\vskip.1in
Let $E$ denote a Hausdroff locally convex topological vector space, and $\mathfrak{F}$ to denote the family of
oontinuous scminorms generating the topology of $E$. Also  $C(E)$  will denote the family of nonempty compact subsets of $E$.
For each $p\in \mathfrak{F}$ and $A, B \in C(E)$, we can define
$\delta(A, B): =  \sup\{p(a - b): a\in A, b \in B\}$. and
 $D_p(A,B):= \max \{\sup_{a \in A}\inf_{b\in B} P( a- b), \sup_{b\in B}\inf_{a \in A} P(a-b) \}$. Though the $P$ is only a seminorrn, $D_p$ is a Hausdroff metric on $C(E)$ (e.g., see Ko and Tsai \cite{kotsai}).

\vskip.1in
\noindent
{\bf Definition 10.1}. Let $K$ he a nonempty subset of $E$. A mapping $T: K \rightarrow C(E)$ is said to be a multi-valued contraction if there exists a constant $k_p \in (0, 1)$ such that $D_p(T(x), T(y)) \leq k_p P(x-y)$.
$T$ is said to be  non-expansive if for any $x, y \in K$, we have $P_p(T(x), T(y)))\leq P(x-y)$.

\vskip.1in
By Chen and Singh \cite{chensingh}, we now have the following definition of Opial's condition in locally convex spaces.

\vskip.1in
\noindent
{\bf Definition 10.2}. The locally convex space $E$ is said to satisfy the Opial's condition if for
each $x \in E$ and every net $(x_{\alpha})$ converging weakly to $x$, then for each  $P \in \mathfrak{F}$, we  have
$\liminf P(x_{\alpha} - y ) > \liminf P(x_{\alpha} -x)$ for any $y\neq x$.

\vskip.1in
Now we have have the following demiclosedness principle for nonexpansive set-valued mappings in (Hausdorff) local convex spaces $E$, which is indeed Theorem 1 of Chen and Singh \cite{chensingh}).

\noindent
{\bf Lemma 10.5}. Let $K$ be a nonempty, weakly, compact and convex subset of $E$. Let $T: K \rightarrow C(E)$ be non-expansive.
If $E$ satisfies the Opial's condition, then  graph $(I-G)$ is closed in $E_w \times E$, where $E_w$ is $E$ with its weak topology and $I$ is the identity mapping.

\noindent
{\bf Proof.} The conclusion follows by Theorem 1 of  Chen and Singh \cite{chensingh}. $\square$

\vskip.1in
\noindent
{\bf Remark 10.2.} When a $p$-vector space $E$ is with a $p$-norm, then both (H1) and (H) conditions for their convergence can be described by the convergence weakly, and strongly  by the weak topology and strong topology induced by $p$-norm for $p \in (0, 1]$. Secondly, if a given $p$-vector space $E$ has a non-empty open $p$-convex subset $U$ containing zero, then any mapping satisfying the $``$(H) condition" is a hemicompact mapping (with respect $P_U$ for a given bounded open $p$-convex subset $U$  containing zero of $p$-vector space $E$), thus satisfying the $``$(H) condition"  used in Theorem 5.1.

\vskip.1in
By the fact that each semiclosed 1-set mappings satisfy the $``$(H1) condition", we have have the existence of fixed points for the class of semiclosed 1-set mappings. First as an application of Theorem 8.2, we have the following result  for non-self mappings in $p$-seminorm spaces for $p \in (0, 1]$.

\vskip.1in
\noindent
{\bf Theorem 10.1}. Let $U$ be a bounded open $p$-convex subset of a $p$-seminorm space $E$ ($0 < p \leq 1)$  the zero $0 \in U$.  Assume $F: \overline{U} \rightarrow E$ is a semiclosed 1-set contractive and continuous mapping. In addition, for any $x\in \partial \overline{U}$,
we have $\lambda x \neq F(x)$ for any $\lambda > 1$ (i.e., the $``$Leray-Schauder boundary condition"). Then $F$ has at least one fixed point.

\noindent
{\bf Proof.} By the proof of Theorem 8.2 with $C= E$, we actually  have the following either (I) or (II) holding:

(I)  $F$ has a fixed point $x_0 \in U $, i.e., $P_U (F(x_0) - x_0) = 0$,

(II)  there exists $x_0 \in \partial(U)$ with $ P_U (F(x_0) - x_0) = (P^{\frac{1}{p}}_U(F(x_0))-1)^p  > 0.$

If $F$ has no fixed point, then above (II) holds and $x_0 \neq F(x_0)$. By the proof of Theorem 8.2, thus $P_U(F(x_0)) > 1$ and $x_0= f(F(x_0))=\frac{F(x_0)}{(P_U(F(x_0)))^{\frac{1}{p}}}$, which means $F(x_0)=(P_U(F(x_0)))^{\frac{1}{p}} x_0$, where $(P_U(F(x_0)))^{\frac{1}{p}} > 1$, this contradicts with the assumption, thus  $F$ must have a fixed point. The proof is complete. $\square$

\vskip.1in
By following the idea used and developed by Browder \cite{bro1968}, Li \cite{li1988}, Li et al.\cite{lixuduan2006}, Goebel and Kirk \cite{goebelkirk1990}, Petryshyn \cite{petryshyn1966}-\cite{petryshyn1973tams}, Tan and Yuan \cite{tanyuan1994}, Xu \cite{xu2007}, Xu et al.\cite{xujiali2006} and references therein,  we have the following existence theorems for the principle of Leray-Schauder type alternatives in $p$-seminorm spaces $(E, \|\cdot \|_p)$  for $p \in (0, 1]$.

\vskip.1in
\noindent
{\bf Theorem 10.2}. Let $U$ be a bounded open $p$-convex subset of a $p$-seminorm space $(E, \|\cdot \|_p)$ ($0 < p \leq 1)$  the zero $0 \in U$.  Assume $F: \overline{U} \rightarrow E$ is a semiclosed 1-set contractive and continuous mapping. In addition, there exist $\alpha >1$, $\beta \geq 0$, such that for each $x \in \partial \overline{U}$, we have
$\|F(x)- x\|_p^{\alpha/p}\geq \|F(x)\|_p^{(\alpha+\beta)/p}\|x\|_p^{-\beta/p}  - \|x\|_p^{\alpha/p}$.
Then $F$ has at least one fixed point.

\noindent
{\bf Proof.} By assuming $F$ has no fixed point, we prove the conclusion by showing the Leray-Schauder boundary condition in Theorem 10.1 does not hold. If we assume $F$ has no fixed point, by the boundary condition of Theorem 10.1,
there exist $x_0\in \partial \overline{U}$, and $\lambda_0 >1$ such that $F(x_0) = \lambda_0 x_0$.

Now, consider the function $f$ defined by $f(t): =(t-1)^{\alpha} - t^{\alpha + \beta}+1$ for $t\geq 1$. We observe that
 $f$ is a strictly decreasing function for $t \in [1, \infty)$ as the derivative of $f$\'{}$(t) =\alpha (t-1)^{\alpha-1} - (\alpha +\beta) t^{\alpha +\beta -1} < 0$  by the differentiation, thus we have
$t^{\alpha + \beta} -1 > (t-1)^{\alpha}$ for $t \in (1, \infty)$. By combining the  boundary condition, we have that
$\|F(x_0)-x_0\|_p^{\alpha/p}=\|\lambda_0 x_0-x_0\|_p^{\alpha/p}=(\lambda_0-1)^{\alpha}\|x_0\|_p^{\alpha/p} < (\lambda_0^{\alpha+\beta}-1)\|x_0\|_p^{(\alpha+\beta)/p}\|x_0\|_p^{-\beta/p}=\|F(x_0)\|_p^{(\alpha+\beta)/p}\|x_0\|_p^{-\beta/p}- \|x_0\|_p^{\alpha/p}$, which contradicts the boundary condition given by Theorem 10.2. Thus, the conclusion follows and the proof is complete. $\square$

\vskip.1in
\noindent
{\bf Theorem 10.3}. Let $U$ be a bounded open $p$-convex subset of a $p$-seminorm space $(E, \|\cdot \|_p)$ ($0 < p \leq 1)$  the zero $0 \in U$.  Assume $F: \overline{U} \rightarrow 2^E$ is a semiclosed 1-set contractive and continuous mapping.
In addition, there exist $\alpha >1$, $\beta \geq 0$, such that for each $x \in \partial \overline{U}$, we have
$\|F(x) + x\|_p^{(\alpha+\beta)/p} \leq \|F(x)\|_p^{\alpha/p}\|x\|_p^{\beta/p}  + \|x\|_p^{(\alpha+\beta)/p}.$
Then $F$ has at least one fixed point.

\noindent
{\bf Proof.} We prove the conclusion by showing the Leray-Schauder boundary condition in Theorem 10.1 does not hold.
If we assume $F$ has no fixed point, by the boundary condition of Theorem 10.1,
there exist $x_0\in \partial \overline{U}$, and $\lambda_0 >1$ such that $F(x_0) = \lambda_0 x_0$.

Now, consider the function $f$ defined by $f(t): =(t+1)^{\alpha+\beta} - t^{\alpha} - 1 $ for $t\geq 1$. We then can show that
$f$ is a strictly increasing function for $t \in [1, \infty)$,  thus we have
$t^{\alpha}+1 < (t + 1)^{\alpha +\beta}$ for $t \in (1, \infty)$.
By the  boundary condition given in Theorem 7.3, we have that
$$\|F(x_0)+x_0\|_p^{(\alpha+\beta)/p}=(\lambda_0 +1)^{\alpha+\beta}\|x_0\|_p^{(\alpha+\beta)/p} > (\lambda_0^{\alpha}+1)\|x_0\|_p^{(\alpha+\beta)/p}=\|F(x_0)\|_p^{\alpha/p}\|x_0\|_p^{\beta/p}+ \|x_0\|_p^{\alpha/p},$$
which contradicts the boundary condition given by Theorem 10.3. Thus, the conclusion follows and the proof is complete. $\square$

\vskip.1in
\noindent
{\bf Theorem 10.4}. Let $U$ be a bounded open $p$-convex subset of a $p$-seminorm space $(E, \|\cdot \|_p)$ ($0 < p \leq 1)$  the zero $0 \in U$.  Assume $F: \overline{U} \rightarrow E$ is a semiclosed 1-set contractive and continuous mapping. In addition, there exist $\alpha >1$, $\beta \geq 0$ (or alternatively, $\alpha >1$, $\beta \geq 0$) such that for each $x \in \partial \overline{U}$, we have that
 $\|F(x) - x\|_p^{\alpha/p} \|x\|_p^{\beta/p} \geq \|F(x)\|_p^{\alpha/p}\|F(x)+x\|_p^{\beta/p} -\|x\|_p^{(\alpha+\beta)/p}.$
Then $F$ has at least one fixed point.

\noindent
{\bf Proof.} The same as above, we prove the conclusion by showing the Leray-Schauder boundary condition in Theorem 10.1 does not hold. If we assume $F$ has no fixed point, by the boundary condition of Theorem 10.1,
there exist $x_0\in \partial \overline{U}$ and $\lambda_0 >1$ such that $F(x_0) = \lambda_0 x_0$.

Now, consider the function $f$ defined by $f(t): =(t-1)^{\alpha} - t^{\alpha}(t-1)^{\beta}+1$ for $t\geq 1$. We then can show that
$f$ is a strictly decreasing function for $t \in [1, \infty)$,  thus we have
$(t-1)^{\alpha} < t^{\alpha} (t+1)^{\beta}-1$ for $t \in (1, \infty)$.

By the boundary condition given in Theorem 10.3, we have that
$$\|F(x_0)-x_0\|_p^{\alpha/p}\|x_0\|_p^{\beta/p}=(\lambda_0-1)^{\alpha}\|x_0\|_p^{(\alpha+\beta)/p} <  (\lambda_0^{\alpha}(\lambda_0+1)^{\beta}-1)\|x_0\|_p^{(\alpha+\beta)/p}=\|F(x_0)\|_p^{\alpha/p}\|y_0+x_0\|_p^{\beta/p}- \|x_0\|_p^{(\alpha+\beta)/p},$$
which contradicts the boundary condition given by Theorem 10.4. Thus, the conclusion follows and the proof is complete. $\square$

\vskip.1in
\noindent
{\bf Theorem 10.5}. Let $U$ be a bounded open $p$-convex subset of a $p$-seminorm space $(E, \|\cdot \|_p)$ ($0 < p \leq 1)$  the zero $0 \in U$.  Assume $F: \overline{U} \rightarrow E$ is a semiclosed 1-set contractive and continuous mapping. In addition, there exist $\alpha >1$, $\beta \geq 0$, we have that $\|F(x) + x\|_p^{(\alpha+\beta)/p} \leq \|F(x)-x\|_p^{\alpha/p}\|x\|_p^{\beta/p} +\|F(x)\|_p^{\beta/p} \|x\|^{\alpha/p}.$
Then $F$ has at least one fixed point.

\noindent
{\bf Proof.} The same as above, we prove the conclusion by showing the Leray-Schauder boundary condition in Theorem 7.1 does not hold. If we assume $F$ has no fixed point, by the boundary condition of Theorem 10.1,
there exist $x_0\in \partial \overline{U}$ and $\lambda_0 >1$ such that $F(x_0) = \lambda_0 x_0$.

Now, consider the function $f$ defined by $f(t): =(t+1)^{\alpha+\beta} - (t-1)^{\alpha}-t^{\beta}$ for $t\geq 1$. We then can show that $f$ is a strictly increasing function for $t \in [1, \infty)$,  thus we have
$(t+1)^{\alpha+\beta} >  (t-1)^{\alpha} +t^{\beta}$ for $t \in (1, \infty)$.

By the  boundary condition given in Theorem 10.3, we have that
$\|F(x_0) +x_0\|_p^{(\alpha+\beta)/p}=(\lambda_0 +1)^{\alpha+\beta}\|x_0\|_p^{(\alpha+\beta)/p} >
 ((\lambda_0-1)^{\alpha}+ \lambda_0^{\beta})\|x_0\|_p^{(\alpha+\beta)/p}=\|\lambda_0 x_0- x_0\|_p^{\alpha/p}\|x_0\|_p^{\beta/p} + \|\lambda_0 x_0\|_p^{\beta/p}\|x_0\|_p^{\alpha/p} = \|F(x_0)-x_0\|_p^{\beta/p}\|x_0\|_p^{\alpha/p} +\|F(x_0)\|_p^{\beta/p}\|x_9\|^{\alpha/p},$
which implies that
$$ \|F(x_0) + x_0\|_p^{(\alpha+\beta)/p} >  \|F(x_0) - x_0\|_p^{\beta/p}\|x_0\|_p^{\alpha/p} +\| F(x_0) \|_p^{\beta/p}\|x_9\|^{\alpha/p},$$
this contradicts the boundary condition given by Theorem 10.5. Thus, the conclusion follows and the proof is complete. $\square$

As an application of Theorems 10.1 by testing the Leray-Schauder boundary condition, we have the following conclusion for each special case, and thus we omit their proofs in details here.

\vskip.1in
\noindent
{\bf Corollary 10.1}. Let $U$ be a bounded open $p$-convex subset of a $p$-seminorm space $(E, \|\cdot \|_p)$ ($0 < p \leq 1)$  the zero $0 \in U$.  Assume $F: \overline{U} \rightarrow E$ is a semiclosed 1-set contractive and continuous mapping. Then $F$ has at least one fixed point if one of the the following (strong)  conditions holds for $x \in \partial \overline{U}$:

(i) $\|F(x)\|_p \leq \|x\|_p$;

(ii) $\|F(x)\|_p \leq \|F(x)-x\|_p$;

(iii) $\|F(x)+x||_p \leq \|F(x)\|_p$;

(iv) $\|F(x) + x\|_p \leq \|x\|_p$;

(v) $\|F(x) + x\|_p \leq \|F(x)- x\|_p$;

(vi) $\|F(x)\|_p \cdot \|F(x)+x\|_p \leq \|x\|_p^2$;

(vii) $\|F(x)\|_p \cdot \|F(x) +x\|_p \leq \|F(x)- x\|_p \cdot \|x\|_p$.

\vskip.1in
If the $p$-seminorm space $E$ is a uniformly convex Banach space $(E, \| \cdot \|)$ (for $p$-norm space with $p=1$),  then we have the following general existence result which can apply to general non-expansive (single-valued) mappings, too.

\vskip.1in
\noindent
{\bf Theorem 10.6}. Let $U$ be a bounded open convex subset of a uniformly convex Banach space $(E, \|\cdot \|)$ (with $p=1$)
with zero $0 \in U$. Assume $F: \overline{U} \rightarrow E$ is a semi-contractive and continuous  (single-valued) mapping.
In addition, for any $x\in \partial \overline{U}$, we have $\lambda x \neq F(x)$ for any $\lambda > 1$ (i.e., the $``$Leray-Schauder boundary condition").
Then $F$ has at least one fixed point.

\noindent
{\bf Proof.} By Lemma 10.1, $F$ is a semiclosed 1-set contractive mapping. Moreover, by the assumption that $E$ is a uniformly convex Banach, the mapping $(I-F)$ is closed at zero, and thus $F$ is semiclosed at zero (see Browder \cite{bro1968}, or Goebel and Kirk \cite{goebelkirk1990}). Thus all assumptions of Theorem 10.2 are satisfied.
The conclusion follows by Theorem 10.2. The proof is completes. $\square$

\vskip.1in
Now we can also have the following result for nonexpansive set-valued mappings (instead of single-valued) in a Banach space $X$ with Opial's condition.

\vskip.1in
\noindent
{\bf Theorem 10.7}.
Let $C$ is a nonempty convex weakly compact subset of a local convex space $X$ which satisfies Opial's condition and $0 \in intC$.
Let $T:  C \rightarrow K(X)$ be a nonexpansive set-valued mapping with non-empty compact convex values.
In addition, for any $x\in \partial \overline{C}$, we have $\lambda x \neq F(x)$ for any $\lambda > 1$ (i.e., the $``$Leray-Schauder boundary condition"). Then $F$ has at least one fixed point.

\noindent
{\bf Proof.} As $T$ is nonexpansive, it is 1-set contractive, By Lemma 10.2,  it is then semi-contractive and continuous. By following the idea of Theorem 10.1, indeed using the proof of Theorem 8.2 (or the similar argument used by Theorem 5.2) by applying Theorem 5.3 (instead of Theorem 5.2) for the fixed point theorem of upper semicontinuous set-valued mappings in locally convex space, the conclusion follows. The proof is complete. $\square$.

\vskip.1in
By using Lemma 10.4, we have the following result in local convex spaces for $*$-nonexpansive single-valued mappings.

\noindent
{\bf Theorem 10.8}.
Let $C$ is a nonempty (bonded) convex closed subset of a Banach space $X$ which is either strictly convex, or with satisfying Opial condition. Let $T:  C \rightarrow X$  be a $*$-nonexpansive and continuous mapping.
In addition, for any $x\in \partial \overline{C}$, we have $\lambda x \neq F(x)$ for any $\lambda > 1$ (i.e., the $``$Leray-Schauder boundary condition"). Then $F$ has at least one fixed point.

\noindent
{\bf Proof.} As $T$ is $*$-nonexpansive, and by the  demiclosednedd principle for $*$-nonexpansive mappings given by Lemma 10.4, it follows that $T$ satisfies the (H1) condition of Theorem 7.1, then all conditions of Theorem 7.1 are satisfied, then the conclusion follows by Theorem 7.1. The proof is complete. $\square$.

\vskip.1in
By considering $p$-seminorm space $(E, \|\cdot\|)$ with a seminorm for $p=1$, the following corollary is a special case of corresponding results from Theorem 10.2 to Theorem 10.5, and thus we omit its proof.

\vskip.1in
\noindent
{\bf Corollary 10.2}. Let $U$ be a bounded open convex subset of a norm space $(E, \|\cdot\|)$.
 Assume $F: \overline{U} \rightarrow E$ is a semiclosed 1-set contractive and continuous mapping. Then $F$ has at least one fixed point if there exist $\alpha >1$, $\beta \geq 0$, such that any one of the following conditions satisfied:

(i) for each $x \in \partial \overline{U}$, $\|F(x)- x\|^{\alpha}\geq \|F(x)\|^{(\alpha+\beta)}\|x\|^{-\beta}  - \|x\|^{\alpha}$;

(ii) for each $x \in \partial \overline{U}$,  $\|F(x) + x\|^{(\alpha+\beta)} \leq \|F(x)\|^{\alpha}\|x\|^{\beta}  + \|x\|^{(\alpha+\beta)}$;

(iii) for each $x \in \partial \overline{U}$, $\|F(x) - x\|^{\alpha} \|x\|^{\beta} \geq \|F(x)\|^{\alpha}\|y+x\|^{\beta} -\|x\|^{(\alpha+\beta)}$;

(iv) for each $x \in \partial \overline{U}$, $\|F(x) + x\|^{(\alpha+\beta)} \leq \|F(x)-x\|^{\alpha}\|x\|^{\beta} +\|F(x)\|^{\beta} \|x\|^{\alpha}$.

\vskip.1in
\noindent
{\bf Remark 10.3.} As discussed by Lemma 10.1 and the proof of Theorem 10.6, when the  $p$-vector space is a uniformly convex Banach space, the semi-contractive or nonexpansive mappings automatically satisfy the conditions (see (H1)) required by Theorem 10.1, that is, the mappings are indeed semiclosed. Moreover, our results from Theorem 10.1 to Theorem 10.6, Corollary 10.1 and Corollary 10.2 also improve or unify corresponding results given by Browder \cite{bro1968}, Li \cite{li1988}, Li et al.\cite{lixuduan2006}, Goebel and Kirk \cite{goebelkirk1990}, Petryshyn \cite{petryshyn1966}-\cite{petryshyn1973tams}, Reich \cite{reich}, Tan and Yuan \cite{tanyuan1994}, Xu \cite{xu1991}, Xu \cite{xu2007}, Xu et al.\cite{xujiali2006}, Yuan \cite{yuan2022} and results from the reference therein by extending the non-self mappings to the classes of semiclosed 1-set contractive set-valued mappings in $p$-seminorm spaces with $p \in (0.1]$ (including the norm space or Banach space when $p=1$ for $p$-seminorm spaces).

\vskip.1in
Before the ending of this paper, we like to share with readers that the main goal of this paper is to develop some new results and tools in the nature way for the category of nonlinear analysis for three classes of mappings which are: 1) condensing; 2) 1-set contractive; and 3) semiclosed mappings under the general framework of locally $p$-convex spaces (where, $(0< p \leq 1)$) for £¨single-valued) continuous mappings, instead of set-valued mappings without the strong condition with closed $p$-convex values! We do also  expect that these new results would become very useful tools for the development of nonlinear functional analysis under the general framework of $p$-vector spaces which include the topological vector spaces as a special classes, and also the related applications for nonlinear problems on
on optimization, nonlinear programming, variational inequality, complementarity, game theory, mathematical economics, and so on.

 Like what mentioned in the beginning of this paper, we do expect that nonlinear results and principles of the best approximation theorem established in this paper would play a very important role for the nonlinear analysis under the general framework of $p$-vector spaces for $(0< p \leq 1)$, as shown by those results given from Sections 6 and 7 for both condensing and 1-set contractive mappings; and general new results in nonlinear analysis from Sections 8, 9 and 10 for semiclosed 1-set contractive mappings for the development of fixed point theorems for non-self mappings, principle of nonlinear alternative, Rothe type, Leray-Schauder alternative, and related topics,  which are not only include corresponding results in the existing literature as special cases, but expected to be important tools for the study of its nonlinear analysis.

 Finally, we like to point out that the work presented by this paper focuses on the development of nonlinear analysis for single-valued (instead of set-valued) mappings for locally $p$-convex spaces, essentially, is very important, and indeed, the continuation of the work given recently by Yuan \cite{yuan2022} therein, the attention is given to establish new results on fixed points, principle of nonlinear alternative for nonlinear mappings mainly on set-valued mappings developed in locally $p$-convex spaces for $0 < p  \leq 1$.
 Though some new results for set-valued mappings in locally $p$-convex spaces have been developed (see Gholizadeh et al.\cite{gholizadeh}, Park \cite{park2016}, Qiu and Rolewicz \cite{qiu}, Xiao and Zhu \cite{xiaozhu2011}-\cite{xiaozhu2018}, Yuan \cite{yuan2022} and others), we still like to emphasize that results obtained for set-valued mappings  for $p$-vector spaces may face some challenging in dealing with true nonlinear problems.  One example is that the assumption used for $``$set-valued mappings with closed $p$-convex values" seems too strong as it always means that the  zero element is a trivial fixed point of the set-valued mappings, and this fact was also discussed in P.40-41 by Yuan \cite{yuan2022} for $0 < p \leq 1$.

\vskip.1in
\noindent
{\bf Acknowledgement}

The author thanks Professor Shih-sen Chang (Shi-Sheng Zhang), Professor K.K. Tan, Professor Hong Ma, Professor Jian Yu, Professor Y.J. Cho, Professor S. Park for their always encouragements in the past more than  two decades; also my thanks go to Professor Hong-Kun Xu, Professor Lishan Liu, Professor Jian-Zhong Xiao, Professor Xiao-Long Qin, Professor Ganshan Yang, Professor Xian Wu, Professor Nanjing Huang, Professor Mohamed Ennassik, Professor Tiexin Guo, and my colleagues and friends across China, Australia, Canada, UK, USA and else where.

\vskip.1in
\noindent
{\bf Funding}

This research is partially supported by the National Natural Science Foundation of China [grant numbers 71971031 and U1811462].

\vskip.1in
\noindent
{\bf Availability of data and materials}

Not applicable.

\vskip.1in
\noindent
{\bf Declarations}

\noindent
{\bf Competing interests}

The authors declare no competing interests.

\noindent
{\bf Author contributions}

Full contribution for research, manuscript writing, reading and approving.

\noindent
{\bf Author details}

1 Business School, Chengdu University, Chengdu 610601, China.

2 College of Mathematics, Sichuan University, Chengdu
610065, China.

3 Business School, Sun Yat-Sen University, Guangzhou 510275, China.

4 Business School, East China
University of Science and Technology, Shanghai 200237, China.

\vskip.1in
\noindent
{\bf Compliance with Ethical Standards}

The author declares that he has no conflict of interest.


\vskip.1in
\noindent
{\bf References}
\begin{thebibliography}{}

\bibitem{Agarwal}R.P. Agarwal, M. Meehan, and D. O'Regan, Fixed Point Theory and Applications, Cambridge Tracts in Mathematics, vol. 141, Cambridge University Press, Cambridge, 2001.

\bibitem{agarwalorgan2003}
R.P. Agarwal, and D. O'Regan, Birkhoff-Kellogg theorems on invariant directions for multimaps, Abstr. Appl. Anal. 2003, no. 7, 435 - 448.

\bibitem{agarwalorgan2004}
R.P. Agarwal, and D. O'Regan, Essential $U_c^k$-type maps and Birkhoff-Kellogg theorems, J. Appl. Math. Stoch. Anal.,
2004, no.1, 1 - 8.

\bibitem{alghamdiporegan}
 M.A. Alghamdi, D. O'Regan, and N. Shahzad, Krasnosel'skii type fixed point theorems for mappings on nonconvex sets,
 Abstr. Appl. Anal. Article, 2020, ID 267531, 1 - 23 (2012).



\bibitem{Askoura}Y. Askoura and C. Godet-Thobie, Fixed points in contractible spaces and convex subsets of topological vector spaces, J. Convex Anal., 2006, 13 (no.2), 193 - 205.

\bibitem{balachandra} V.K. Balachandran, Topological Algebras, vol. 185, Elsevier, Amsterdam, 2000.

\bibitem{bayoumi}A. Bayoumi, Foundations of complex analysis in non locally convex spaces. Function theory without convexity condition, North-Holland Mathematics Studies, Vol. 193. Elsevier Science B.V., Amsterdam, 2003.

\bibitem{bayoumi2015}A. Bayoumi, N. Faried, N., and R. Mostafa, Regularity properties of p-distance transformations in image analysis,
Int. J. Contemp. Math. Sci., 2015,  10, 143 - 157.

\bibitem{berbstein}S. Bernstein, Sur les equations de calcul des variations, Ann. Sci. Ecole Normale Sup.,1912, 29, 431 - 485.

\bibitem{bernues1997}
J. Bernu\'{e}es, and A. Pena, On the shape of p-convex hulls $0 < p < 1$, Acta Math. Hungar., 1997, 74, no.4, 345 - 353.

\bibitem{bk}G.D. Birkhoff, and O.D. Kellogg, Invariant points in function space, Trans. Amer. Math. Soc., 1922, 23, no.1, 96 - 115.

\bibitem{bro1965}
F.E. Browder, Nonexpansive nonlinear operators in a Banach space,  Proc. Nat. Acad. Sci. U.S.A., 1965, 54,  1041 - 1044.

\bibitem{bro1967}
F.E. Browder, Convergence of approximants to fixed points of nonexpansive
nonlinear mappings in Banach spaces, Arch. Ration. Mech. Anal., 1967, 24(1), 82 - 90.

\bibitem{bro1968b}
F.E. Browder, The fixed point theory of multi-valued mappings in topological vector spaces, Math. Ann., 1968, 177, 283 - 301.

\bibitem{bro1968}
F.E. Browder, Semicontractive and semiaccretive nonlinear mappings in Banach spaces, Bull. Amer. Math.
Soc., 1968, 74, 660 - 665.


\bibitem{bro1976} F.E. Browder, Nonlinear Functional Analysis, Proc. Sympos. Pure Math., 18, Part
2, Amer. Math. Soc., Providence, RI, 1976.

\bibitem{bro1983}
F.E. Browder, Fixed point theory and nonlinear problems, Bull. Amer. Math. Soc. (N.S.), 1983, 9,  no.1, 1 - 39.





\bibitem{carbone1991}
A. Carbone, and G. Conti, Multivalued maps and existence of best approxima- tions, Jour. Approx. Theory, 1991, 64,  203 - 208.

\bibitem{cauty2005}R. Cauty, R\'etract\`es absolus de voisinage alg\'ebriques.(French) [Algebraic absolute neighborhood retracts], Serdica Math. J., 2005,  31(4), 309 - 354.

\bibitem{cauty2007}R. Cauty, Le th\'{e}or\.{e}me de Lefschetz - Hopf pour les applications compactes des espaces ULC. (French) [The Lefschetz - Hopf theorem for compact maps of uniformly locally contractible spaces], J. Fixed Point Theory Appl., 2007, 1(1), 123 - 134.

\bibitem{chang1997}
S.S. Chang, Some problems and results in the study of nonlinear analysis, Proceedings of the Second World Congress of Nonlinear Analysts, Part 7 (Athens, 1996), Nonlinear Anal., 1997,  30, no.7, 4197 - 4208.


\bibitem{chang1993}
S.S. Chang, Y.J. Cho, and Y.Zhang, The topological versions of KKM theorem and Fan's matching theorem with applications,
Topol. Methods Nonlinear Anal.,  1993, 1, no. 2, 231 - 245.






\bibitem{chang2001}
T.H. Chang, Y.Y. Huang, and J.C. Jeng, Fixed point theorems for multi- functions in S-KKM class, Nonlinear Anal., 2001, 44, 1007 - 1017.


\bibitem{chang1999}
T.H. Chang, Y.Y. Huang, J.C. Jeng, and K.H. Kuo, On S-KKM property and related topics, Jour. Math. Anal. Appl., 1999, 229, 212 - 227.


\bibitem{chang1996}
T.H. Chang, and C.L. Yen, KKM property and fixed point theorems, Jour. Math. Anal. Appl., 1996,203,224 - 235.

\bibitem{chen}Y.Q. Chen, Fixed points for convex continuous mappings in topological vector spaces, Proc. Amer. Math. Soc., 2001, 129, no. 7, 2157 - 2162.

\bibitem{chensingh}
Y.K. Chen, and  K.L. Singh, Fixed points for nonexpansive multivalued mapping and the Opial's condition,
J\~{n}\={a}n\={a}bha, 1992, 22, 107 - 110.

\bibitem{darbo}
G. Darbo, Punti uniti in trasformazioni a condominio non compatto, Rend. Sem. Mat. Univ. Padova, 1955, 24, 84 - 92.

\bibitem{ding}
G.G. Ding, New Theory in Functional Analysis, Academic Press, Beijing, 2007.




\bibitem{dobrowolski2003}T. Dobrowolski, Revisiting Cauty's proof of the Schauder conjecture, Abstr. Appl. Anal., 2003, 7, 407 - 433.



\bibitem{ennassik2021}M. Ennassik, and M.A. Taoudi, On the conjecture of Schauder., J. Fixed Point Theory Appl., 2021, 23,
no.4, Paper No. 52, 15pp.

\bibitem{ennassik}M. Ennassik, L. Maniar and M.A. Taoudi, Fixed point theorems in r-normed and locally $r$-convex spaces and applications, Fixed Point Theory, 2021, 22(2), 625 - 644.

\bibitem{fan1952}K. Fan, Fixed-point and minimax theorems in locally convex topological linear spaces, Proc. Nat. Acad. Sci. U.S.A., 1952,  38, 121 - 126.

\bibitem{fan1960}K. Fan, A generalization of Tychonoff's fixed point theorem, Math. Ann., 1960/61, 142, 305 - 310.

\bibitem{fan1969}K. Fan, Extensions of two fixed point theorems of F. E. Browder. Math. Z., 1969, 112, 234 - 240.


\bibitem{fan1972}K. Fan, A minimax inequality and applications. Inequalities, III (Proc. Third Sympos., Univ. California, Los Angeles, Calif., 1969; dedicated to the memory of Theodore S. Motzkin), 1972, pp. 103 - 113. Academic Press, New York.

\bibitem{furipera}
M. Furi, and M.P. Pera, A continuation method on locally convex spaces and applications to ordinary differential equations on noncompact intervals, Ann. Polon. Math., 1987, 47, no.3, 331 - 346.

\bibitem{gal}
S.G. Gal, and J.A. Goldstein, Semigroups of linear operators on p-Fr\'echet spaces $0 < p <1$, Acta Math. Hungar., 2007,
114 (1-2), 13 - 36.

\bibitem{gholizadeh}L. Gholizadeh, E.Karapinar, and M.Roohi, Some fixed point theorems in locally $p$-convex spaces, Fixed Point Theory Appl., 2013, 2013:312, 10 pp.


\bibitem{goebel} K. Goebel, On a fixed point theorem for multivalued nonexpansive mappings, Ann. Univ. Mariae Curie-Sk?odowska Sect. A,  1975, 29, 69 - 72.

\bibitem{goebelkirk1990} K. Goebel, W.A. Kirk, Topics in metric fixed point theory, Cambridge Studies in Advanced Mathematics, 28,
Cambridge University Press, Cambridge, 1990.

\bibitem{goebelkirk2008} K. Goebel, and W.A. Kirk, Some problems in metric fixed point theory, J. Fixed Point Theory Appl., 2008, 4, no.1, 13 -  25.

\bibitem{gohde}
D. G\"ohde, Zum Prinzip der kontraktiven Abbildung, Math. Nachr., 1965, 30, 251 - 258.

\bibitem{gorniewicz} L. G\'{o}rniewicz, Topological fixed point theory of multivalued mappings. Mathematics and its Applications, vol. 495, Kluwer Academic Publishers, Dordrecht, 1999.

\bibitem{gorniewiczetal}
L. G\'{o}rniewicz, A. Granas, W. Kryszewski, On the homotopy method in the fixed point index theory of multivalued mappings of compact absolute neighborhood retracts, J. Math. Anal. Appl., 1991, 161, no. 2, 457 - 473.


\bibitem{granas}A.Granas and J.Dugundji, Fixed Point Theory, Springer Monographs in Mathematics. Springer-Verlag, New York, 2003.

\bibitem{hergman1968}
B.R. Halpern and G.H. Bergman, A fixed-point theorem for inward and outward maps, Trans. Amer. Math. Soc., 1965, 130,
 353 - 358.

\bibitem{huang1997}
N.J. Huang, B.S. Lee, and M.K. Kang, Fixed point theorems for compatible mappings with applications to the solutions of functional equations arising in dynamic programmings, Internat. J. Math. Math. Sci., 1997,  20, no. 4, 673 - 680.


\bibitem{husainlatif} T. Husain, and A. Latif, Fixed points of multivalued nonexpansive maps, Math. Japonica, 1988, 33, 385 - 391.

\bibitem{husaintarafdar}T. Husain, and E. Tarafdar,  Fixed point theorems for multivalued mappings of nonexpansive type,
Yokohama Math. J., 1980, 28, no. 1-2, 1 - 6.

\bibitem{jarchow}H. Jarchow, Locally Convex Spaces, Stuttgart: B.G. Teubner, 1981.

\bibitem{isac}G. Isac, Leray-Schauder Type Alternatives, Complementarity Problems and Variational Inequalities. Nonconvex Optimization and Its Applications, vol. 87. Springer, New York, 2006.


\bibitem{kalton1977}N.J. Kalton, Compact p-convex sets, Q.J.Math.Oxf. Ser., 1977, 28, no.2, 301 - 308.

\bibitem{kalton1977b}N.J. Kalton, Universal spaces and universal bases in metric linear spaces, Studia Math., 1977, 61, 161 - 191.


\bibitem{kalton1984}N.J. Kalton, N.T. Peck, and J.W. Roberts, An F-Space Sampler, London Mathematical Society Lecture Note Series, vol. 89. Cambridge University Press, Cambridge, 1984.

\bibitem{kaniok}L. Kaniok, On measures of non compactness in general topological vector spaces, Comment. Math.
Univ. Carol., 1990,  31(3),  479 - 487.


\bibitem{kim2002}I.S. Kim, K. Kim, and S.Park, Leray-Schauder alternatives for approximable maps in topological vector spaces,
Math. Comput. Modelling, 2002,  35,  385 - 391.

\bibitem{kirk2014}
W. Kirk, and N. Shahzad, Fixed Point Theory in Distance Spaces, Springer, Cham, 2014.


\bibitem{klee1960}V. Klee, Convexity of Chevyshev sets. Math. Ann., 1960/61, 142, 292 - 304.

\bibitem{kkm}H. Knaster, C. Kuratowski, and S.Mazurkiwiecz, Ein beweis des fixpunktsatzes f\"ur n-dimensional simplexe, Fund. Math., 1929, 63, 132 - 137.

\bibitem{kotsai}
H.M. Ko, and Y.H. Tsai, Fixed point theorems for ptiint to set mappings in locally convex spaces and a characterization of
complete metric spaces, Bull. Academia Sinica, 1979, 7, no.4, 461 - 470.

\bibitem{kozlovetal}
V. Kozlov, J. Thim, and B. Turesson, A fixed point theorem in locally convex spaces, Collect. Math., 2010, 61, no. 2, 223 - 239.

\bibitem{kuratowski}K. Kuratowski, Sur les espaces complets, Fund. Math., 1930, 15, 301 - 309.

\bibitem{lamidozo} E. Lami Dozo, Multivalued nonexpansive mappings and Opial's condition, Proc. Amer. Math. Soc., 1973, 38, 286 - 292.

\bibitem{lerayschauder}
J. Leray, and J. Schauder,Topologie et equations fonctionnelles, Ann. Sci. Ecole Normale Sup.,1934,51,45-78.















\bibitem{li1988}
G.Z. Li, The fixed point index and the fixed point theorems of 1-set-contraction mappings, Proc. Amer. Math. Soc.,
1988,  104, 1163 - 1170.

\bibitem{lixuduan2006}
G.Z. Li, S.Y. Xu, and H.G. Duan, Fixed point theorems of 1-set-contractive operators in Banach spaces,
Appl. Math. Lett., 2006, 19, no. 5, 403 - 412.

\bibitem{li}J.L. Li, An extension of Tychonoff's fixed point theorem to pseudonorm adjoint topological vector spaces, Optimization, 2021, 70, no.5-6, 1217 - 1229.

\bibitem{lim1974}
T.C. Lim, A fixed point theorem for multivalued nonexpansive mappings in a uniformly
convex Banach space, Bull. Amer. Math. Soc., 1974, 80, 1123 - 1126.

\bibitem{liu2001}L.S. Liu, Approximation theorems and fixed point theorems for various classes of 1-set-contractive mappings in Banach spaces, Acta Math. Sin. (Engl. Ser.), 2001, 17, no.1, 103 - 112.

\bibitem{lopezxu}G. L\'{o}pez-Acedo, H.K. Xu, Remarks on multivalued nonexpansive mappings, Soochow J. Math., 1995, 21(1), 107 - 115.

\bibitem{machrafioubbi} N. Machrafi, and L. Oubbi, Real-valued non compactness measures in topological vector spaces and a pplications. [Corrected title: Real-valued non compactness measures in topological vector spaces and applications],  Banach J. Math. Anal., 2020, 14, no.4, 1305 - 1325.

\bibitem{manka}R. Ma\'{n}ka, The topological fixed point property - an elementary continuum-theoretic approach. Fixed point theory and its applications, 2007, 183 - 200, Banach Center Publ., 77, Polish Acad. Sci. Inst. Math., Warsaw,



\bibitem{Mauldin}R.D. Mauldin, The Scottish book, Mathematics from the Scottish Caf\'{e} with selected problems from the new Scottish book, Second Edition, Birkhauser, 2015.
\bibitem{mugliamarino}
L. Muglia, and G. Marino, Some results on the approximation of solutions of variational inequalities for multivalued maps on Banach spaces, Mediterr. J. Math., 2021, 18,  no.4, Paper No. 157, 19 pp.


\bibitem{nhu1996}N.T. Nhu, The fixed point property for weakly admissible compact convex sets: searching for a solution to Schauder's conjecture, Topology Appl., 1996, 68, no. 1, 1 - 12.

\bibitem{nussbaum}R.D. Nussbaum, The fixed point index and asymptotic fixed point theorems for k-set-contractions,
Bull. Amer. Math. Soc., 1969, 75, 490 - 495.

\bibitem{okon}T. Okon, The Kakutani fixed point theorem for Robert spaces, Topology Appl.,2002, 123, no.3, 461 - 470.

\bibitem{opial1967}
Z. Opial, Weak convergence of the sequence of successive approximations for nonexpansive
mappings, Bull. Amer. Math. Soc., 1967, 73, 595 - 597.

\bibitem{oregan2021}
D. O'Regan, Continuation theorems for Monch countable compactness-type set-valued maps,  Appl. Anal., 2021, 100,
 no.7, 1432 - 1439.

\bibitem{oregan2019}
D. O'Regan, Abstract Leray-Schauder type alternatives and extensions,
An.tiin. Univ. $``$Ovidius" Constana Ser. Mat., 2019, 27, no.1, 233 - 243.

\bibitem{OP}D. O'Regan, and R. Precup, Theorems of Leray-Schauder Type and Applications, Gordon and Breach Science
Publishers, 2001.

%

\bibitem{oubbi}
L. Oubbi, Algebras of Gelfand-continuous functions into Arens-Michael algebras, Commun. Korean Math. Soc., 2019, 34(2), 585 - 602.

\bibitem{park1991}S. Park, Some coincidence theorems on acyclic multifunctions and applications to KKM theory,
In, Fixed point theory and applications (Halifax, NS, 1991), 1992, 248 - 277, World Sci. Publ., River Edge, NJ.


\bibitem{park1995}
S. Park, Generalized Leray-Schauder principles for compact admissible multifunctions, Topol. Methods Nonlinear Anal., 1995, 5,
 no. 2, 271 - 277.

\bibitem{park1995b}
S. Park, Acyclic maps, minimax inequalities and fixed points. Nonlinear Anal., 24 (1995), no. 11, 1549 - 1554.

\bibitem{park1997}S. Park, Generalized Leray-Schauder principles for condensing admissible multifunctions, Ann. Mat. Pura Appl., 1997, 172, no.4, 65 - 85.

\bibitem{park2010}S. Park, The KKM principle in abstract convex spaces: equivalent formulations and applications, Nonlinear Anal., 2010, 73, no. 4, 1028 - 1042.

\bibitem{park2016}S. Park, On the KKM theory of locally $p$-convex spaces (Nonlinear Analysis and Convex Analysis), Institute of Mathematical Research(Kyoto University), 2016, Vol. 2011, 70 - 77 (http://hdl.handle.net/2433/231597), Kyoto University, Japan.


\bibitem{park100}S. Park, One hundred years of the Brouwer fixed point theorem,
J. Nat. Acad. Sci., ROK, Nat. Sci. Ser., 2021, 60(1), 1 - 77.

\bibitem{park2022}S. Park, Some new equivalents of the Brouwer fixed point theorem, Advances in the Theory of Nonlinear Analysis and its Applications, 2022, 6, No.3, 300 - 309 (https://doi.org/10.31197/atnaa.1086232.)


\bibitem{petryshyn1966}
W.V. Petryshyn, Construction of fixed points of demicompact mappings in Hilbert space, J. Math. Anal.
Appl., 1966, 14, 276 - 284.

\bibitem{petryshyn1973tams}
W.V. Petryshyn, Fixed point theorems for various classes of 1-set-contractive and 1-ball-contractive mappings in Banach spaces,
Trans. Amer. Math. Soc., 1973, 182, 323 - 352.

\bibitem{petrusel2014} A. Petrusel, I.A. Rus, M.A. Serban, Basic problems of the metric fixed point theory and the relevance of a metric fixed point theorem for a multivalued operator, J. Nonlinear Convex Anal., 2014, 15, no.3, 493 - 513.

\bibitem{pietramala}
P. Pietramala, Convergence of approximating fixed points sets for multivalued
nonexpansive mappings, Comment. Math. Univ. Carol., 1991, 32(4), 697 - 701.

\bibitem{poincare2}
H. Poincare, Sur un theoreme de geometric, Rend. Circ. Mat. Palermo, 1912, 33,  357 -  407.

\bibitem{potter1972}
A.J.B. Potter, An elementary version of the Leray-Schauder theorem, J. London Math. Soc., 1972, 5(2), 414 - 416.

\bibitem{qiu}J. Qiu, and S. Rolewicz,
       Ekeland's variational principle in locally p-convex spaces and related results, Stud. Math., 2008, 186, No.3, 219 - 235.

\bibitem{reich}S. Reich, Fixed points in locally convex spaces. Math. Z., 1972, 125, 17 - 31.


\bibitem{robertson} L.B. Robertson,
Topological vector spaces, Publ. Inst. Math., 1971, 12 (26), 19 - 21.

\bibitem{roberts1977}J.W. Roberts, A compact convex set with no extreme points, Studia Math., 1977, 60, no.3, 255 - 266.

\bibitem{rolewicz}S. Rolewicz, Metric Linear Spaces, Warszawa: PWN-Polish Scientific Publishers, 1985.


\bibitem{rothe1981}
E.H. Rothe, Some homotopy theorems concerning Leray-Schauder maps. Dynamical systems, II (Gainesville, Fla.,
1981), 1982, 327 - 348, Academic Press, New York-London.


\bibitem{rothe1986}
E.H. Rothe, Introduction to various aspects of degree theory in Banach spaces, Mathematical Surveys and Monographs, 23,  American Mathematical Society, Providence, RI, 1986.


\bibitem{sadovskii}B.N. Sadovskii, On a fxed point principle [in Russian], Funkts. Analiz Prilozh., 1967,  1(2), 74 - 76.

\bibitem{zezer2021}
S. Sezer, Z. Eken, G. Tinaztepe, and G. Adilov, $p$-convex functions and some of their properties,
Numer. Funct. Anal. Optim.,  2021, 42, no.4, 443 - 459.


\bibitem{shahzad2004a}
N. Shahzad, Fixed point and approximation results for multimaps in $S-KKM$ class,  Nonlinear Anal., 2004, 56, no.6, 905 - 918.

\bibitem{shahzad2004}
N. Shahzad, Approximation and Leray-Schauder type results for $U_{c}^{k}$ maps, Topol. Methods Nonlinear Anal., 2004, 24, no.2, 337 - 346.

\bibitem{shahzad2006}
N. Shahzad, Approximation and Leray-Schauder type results for multimaps in the S-KKM class,
Bull. Belg. Math. Soc., 2006,  13, no.1, 113 - 121.


\bibitem{simons}S. Simons, Boundness in linear topological spaces, Trans. Amer. Math. Soc., 1964, 113, 169 - 180.

\bibitem{singh1997}S.P. Singh, B. Watson, and F. Srivastava, Fixed Point Theory and Best Approximation: the KKM-map principle,
 Mathematics and its Applications, vol. 424, Kluwer Academic Publishers, Dordrecht, 1997.






























%




\bibitem{schauder}
J. Schauder, Der Fixpunktsatz in Funktionalraumen, Stud. Math., 1930, 2, 171 - 180.



\bibitem{silva} E.B. Silva, D.L. Fernandez, and L. Nikolova, Generalized quasi-Banach sequence spaces and measures of noncompactness, An. Acad. Bras. Cie., 2013, 85(2), 443 - 456.

\bibitem{smart}D.R. Smart, Fixed Point Theorems, Cambridge University Press, Cambridge, 1980.



\bibitem{tabor}
J.A. Tabor, J.O. Tabor, M. Idak, Stability of isometries in $p$-Banach spaces, Funct. Approx., 2008,  38, 109 - 119.

\bibitem{tan}
D.N. Tan, On extension of isometries on the unit spheres of $L^p$ - spaces for $0 < p \leq 1$,  Nonlinear Anal., 2011, 74,
6981 - 6987.

\bibitem{tanyuan1994}
K.K. Tan, and X.Z. Yuan, Random fixed-point theorems and approximation in cones, J. Math. Anal. Appl., 1994, 185, 378 - 390.

\bibitem{tych1935}A. Tychonoff, Ein Fixpunktsatz, Math. Ann., 1935, 111, 767 - 776.



\bibitem{wang}J.Y. Wang, An Introduction to Locally $p$-Convex Spaces, 2013, pp. 26 - 64, Academic Press, Beijing.


\bibitem{weber1991}
H. Weber, Compact convex sets in non-locally convex linear spaces, Schauder-Tychonoff fixed point theorem, In, Topology, measures, and fractals (Warnemunde, 1991), 1992,  37 - 40, Math. Res., 66, Akademie-Verlag, Berlin.

\bibitem{weber1992}
H. Weber, Compact convex sets in non-locally convex linear spaces. Dedicated to the memory of Professor Gottfried K$\ddot{\mbox{o}}$the,  Note Mat., 1992, 12, 271 -  289.


\bibitem{xiaolu}J.Z. Xiao, and Y. Lu, Some fixed point theorems for $s$-convex subsets in $p$-normed spaces based on measures ofnoncompactness, J. Fixed Point Theory Appl., 2018, 20, no.2, Paper No.83, 22 pp.

\bibitem{xiaozhu2018}
J.Z. Xiao, and X.H. Zhu, The Chebyshev selections and fixed points of set-valued mappings in Banach spaces with some uniform convexity,
Math. Comput. Modelling,  2011, 54, no.5-6, 1576 -  1583.


\bibitem{xiaozhu2011}J.Z. Xiao, and X.H. Zhu, Some fixed point theorems for $s$-convex subsets in $p$-normed spaces,
Nonlinear Anal., 2011, 74, no.5, 1738 - 1748.



\bibitem{xu1991}
H.K. Xu, Inequalities in Banach spaces with applications. Nonlinear Anal., 1991,  16, no.12, 1127 - 1138.

\bibitem{xu1991}
H.K. Xu, On weakly nonexpansive and ?-nonexpansive multivalued mappings, Math. Jpn., 1991,  36(3), 441 - 445.


\bibitem{xu2000}
H.K. Xu, Metric fixed point theory for multivalued mappings, Dissertationes Math. (Rozprawy Mat.), 2000, 389, 39 pp.

\bibitem{xumuglia} H.K. Xu, and L. Muglia, On solving variational inequalities defined on fixed point sets of multivalued mappings in Banach spaces, J. Fixed Point Theory Appl., 2020, 22, no.4, Paper No.79, 17 pp.

\bibitem{xu2007}
S.Y. Xu, New fixed point theorems for 1-set-contractive operators in Banach spaces, Nonlinear Anal.,
 2007, 67, no.3, 938 - 944.

\bibitem{xujiali2006}
S.Y. Xu, B.G. Jia, and G.Z. Li, Fixed points for weakly inward mappings in Banach spaces,
J. Math. Anal. Appl., 2006, 319, no. 2, 863 - 873.




\bibitem{yangai1980} K. Yanagi, On some fixed point theorems for multivalued mappings,
Pacific J. Math., 1980, 87, no. 1, 233 - 240.

\bibitem{yuan1998}
G.X. Yuan, The study of minimax inequalities and applications to economies and variational inequalities, Mem. Amer. Math. Soc., 1998, 132, no. 625.

\bibitem{yuan1999}
G.X. Yuan, KKM Theory and Applications in Nonlinear Analysis, Monographs and Textbooks in Pure and Applied Mathematics, vol. 218. Marcel Dekker, Inc., New York, 1999.


\bibitem{yuan2022}
G.X. Yuan, Nonlinear analysis by applying best approximation method in p-vector spaces, Fixed Point Theory Algorithms Sci Eng, vol. 20 (2022), pp.1-45.
https://doi.org/10.1186/s13663-022-00730-x, 2022.






\bibitem{Zeidler}
E. Zeidler, Nonlinear functional analysis and its applications, vol. I, Fixed-Point Theorems, Springer Verlag, New York, 1986.







\end{thebibliography}
\end{document}